%% file: ex_article_v3.tex
\documentclass[review]{siamart0216}

\input{ex_shared}

\ifpdf
\hypersetup{
  pdftitle={\TheTitle},
  pdfauthor={\TheAuthors}
}
\fi

\usepackage{amssymb,amsmath,bm}
\usepackage{graphics}
\usepackage{epsfig}
\usepackage{algorithm} 
\usepackage{multirow} 
\usepackage{xcolor}
\usepackage{makecell}

\theoremstyle{assumption}
\newtheorem{assumption}[theorem]{Assumption}
\theoremstyle{assumptions}

\usepackage{setspace}
\usepackage[normalem]{ulem}
\newtheorem{remark}[theorem]{Remark}

\usepackage{setspace}
\newenvironment{citeRef}
{\vskip 0.1cm \begin{singlespace*}\textbf{References:}\begin{small}}
		{\end{small}\end{singlespace*}}

\DeclareMathOperator*{\argmax}{\textup{arg\,max}}

\usepackage[]{algorithm}
\usepackage{algorithmic} 
\usepackage{makecell}
\makeatletter
\newcommand{\xRightarrow}[2][]{\ext@arrow 0359\Rightarrowfill@{#1}{#2}}
\makeatother

\begin{document}

\maketitle

\begin{abstract}
In this paper, we propose an infinite-dimensional version of the Stein variational gradient descent (iSVGD) method for solving  Bayesian inverse problems. The method can generate approximate samples from posteriors efficiently. Based on the concepts of operator-valued kernels and vector-valued reproducing kernel Hilbert spaces, a rigorous definition is given for the infinite-dimensional objects, e.g., the Stein operator, which are proved to be the limit of finite-dimensional ones. Moreover, a more efficient iSVGD with preconditioning operators is constructed by generalizing the change of variables formula and introducing a regularity parameter. The proposed algorithms are applied to an inverse problem of the steady state Darcy flow equation. Numerical results confirm our theoretical findings and demonstrate the potential applications of the proposed approach in the posterior sampling of large-scale nonlinear statistical inverse problems.
\end{abstract}



\section{Introduction}

Driven by rapid algorithmic development and a steady increase of computer power, the Bayesian approach has enjoyed great popularity for solving inverse problems over the last decade. By transforming inverse problems into statistical inference problems, the approach provides a general framework to quantify uncertainties \cite{Arridge2019AN}. The posterior distribution automatically delivers an estimate of the statistical uncertainty in the reconstruction, and hence suggests ``confidence'' intervals that allow to reject or accept scientific hypotheses \cite{Nickl2020JEMS}. It has been widely used in many applications, e.g., artifact detecting in medical imaging \cite{Zhou2020SIIMS}.

The approach begins with establishing an appropriate Bayes model. When the parameters are in a finite-dimensional space, the finite-dimensional Bayesian method can be employed \cite{Tarantola1982JG}. A comprehensive account of the finite-dimensional theory can be found in \cite{Kaipio2004Book}. When the inferred parameters are in the infinite-dimensional space, the problems are more challenging since the Lebesgue measure cannot be defined rigorously in this case \cite{DaPrato1996Book}. Recently, some attempts have been made to handle the issue. For example, a general framework was designed  for the Bayesian formula and the general theory was applied to inverse problems of fluid mechanic equations \cite{Cotter2009IP}. A survey can be found in  
\cite{Stuart2010AN} on the basic framework of the infinite-dimensional Bayes' approach for solving inverse problems. Inverse problems of partial differential equations (PDEs) often involve infinite-dimensional spaces, and the infinite-dimensional Bayes' theory has recently attracted more attention \cite{Burger2014IP,Cotter2013SS,Helin2015IP,Pinski2015SIAMSC,Pinski2015SIAMMA}.

As pointed out in \cite{Arridge2019AN}, one of the challenges for the Bayesian approach is how to effectively extract information encoded in the posterior probability measure. To overcome the difficulty, the two main strategies are the point estimate method and the sampling method. The former is to find the maximum a posteriori (MAP) estimate which is equivalent to solve an  optimization problem \cite{Burger2014IP,Helin2015IP}. In some situations, the MAP estimates are more desirable and computationally feasible than the entire posterior distribution \cite{Jia2019IP,Tarantola2005Book}. However, the point estimates cannot convey uncertainty information and are usually recognized as an incomplete Bayes' method. The sampling type methods, such as the well known Markov chain Monte Carlo (MCMC), are often used to extract posterior information. They are well studied in the finite-dimensional setting \cite{Levin2017Book}. Although the MCMC methods are accurate and effective, they are usually not robust under mesh refinement \cite{Cotter2013SS}. Multiple dimension-independent MCMC-type algorithms have been
proposed \cite{Cotter2013SS,Cui2016JCP,Feng2018SISC,Spantini2015SISC}. However, these MCMC-type algorithms are computationally too expensive to be adopted in such an application as seismic exploration \cite{Fichtner2011Book}.

The finite-dimensional problems have been extensively studied and many efficient algorithms have been developed to quantify uncertainties effectively. In particular, the variational inference (VI) methods have been broadly investigated in machine learning \cite{Bishop2006PRML,Matthews2016PHD,Zhang2018IEEE,Zhao2015IEEE}. Under the mean-field assumption, the linear inverse problems were examined in \cite{Jin2010JCP,Jin2012JCP} by using a hierarchical formulation with Gaussian and centered-t noise distribution. The skewed-t noise distribution was considered for a similar setting in \cite{Guha2015JCP}. A new type of variational inference algorithm, called the Stein variational gradient descent (SVGD), was proposed in \cite{Liu2016NIPS}. The method can achieve reliable uncertainty estimation by efficiently using an interacting repulsive mechanism. The SVGD has shown to be a fast and flexible method for solving challenging machine learning problems and inverse problems of PDEs \cite{Chen2020SISC,Chen2019NIPS}. 

Compared with the finite-dimensional problems, the infinite-dimensional problems are much less studied for the variational inference (VI). When the approximate measures are restricted to be Gaussian, the novel Robbins--Monro algorithm was developed 
in \cite{Pinski2015SIAMSC,Pinski2015SIAMMA} from a calculus-of-variations viewpoint. It was shown in \cite{Sun2019ICLR} that the Kullback--Leibler (KL) divergence between the stochastic processes is equal to the supremum of the KL divergence between the measures restricted to finite marginals. Meanwhile, they developed a VI method for functions parameterized by Bayesian neural networks. Under the classical mean-field assumption, a general VI framework defined on separable Hilbert spaces was proposed recently in \cite{Jia2020SISC}. A function space particle optimization method including the SVGD was developed in \cite{Wang2019ICLR} to solve the particle optimization directly in the space of functions. The function space algorithm was also employed to solve computer vision problems, e.g., the context of semantic segmentation and depth estimation \cite{Carvalho2020CVPR}. However, the function spaced SVGD assumes that the random functions can be parameterized by a finite number of parameters, e.g., parameterized by some neural networks \cite{Wang2019ICLR}. Hence, the probability measures on functions are implicitly defined through the probability distributions of a finite number of parameters, instead of the expected infinite-dimensional function space.

This work concerns inverse problems of PDEs imposed on infinite dimensional function spaces. Motivated by the preconditioned Crank--Nicolson (pCN) algorithm \cite{Cotter2013SS}, we aim to construct the SVGD on separable Hilbert spaces with random functions. Throughout, the iSVGD stands for SVGD defined on the infinite-dimensional function space. The goal is to develop algorithms defined on Hilbert spaces and lay a foundation for appropriate discretizations. It contains three  contributions:  

\begin{enumerate}
	\item[(1)] We investigate the Bayesian formula in infinite-dimensional spaces. The rigorous definition of the
	SVGD on separable Hilbert spaces is provided, the Stein operator is defined 
	and the corresponding optimization problem on some Hilbert spaces is considered, and the finite-dimensional problem is proved to converge to the infinite-dimensional counterpart; 
	
	\item[(2)] By introducing vector-valued reproducing kernel Hilbert space (RKHS) and operator-valued kernel,
	we improve the iSVGD with precondition information (e.g., Hessian information operator), which can 
	accelerate the iSVGD algorithm significantly. This is the first work on such an iSVGD algorithm with precondition information; 
	
	\item[(3)] Explicit numerical strategies are designed by using the finite-element approach.
	Through theoretical analysis and numerical examples, we demonstrate that the regularity parameter $s$ introduced in the abstract theory (see Assumptions \ref{assumPhi} and \ref{assKernel} in Section \ref{ConstructSVGDSection_const})
	should belong to the interval $(0,0.5)$ and be close to $0.5$.
	The scalability of the algorithm depends only on the scalability of the forward and adjoint PDE solvers.  
	Hence, the algorithm is applicable to solve large-scale inverse problems of PDEs. 
\end{enumerate}

The paper is organized as follows. The SVGD in finite-dimensional spaces is introduced in Section 2. Section 3 is devoted to the construction of the iSVGD. The basic concepts of operator-valued kernels and Hilbert scales are briefly reviewed; the Stein operator is defined on separable Hilbert spaces; it is shown that the infinite-dimensional version is indeed equivalent to the finite-dimensional version in some limit sense; Based on the Stein operator and the theory of reproducing kernel Hilbert space (RKHS), the update direction of the iSVGD is derived; In addition, the change of variables is studied and the iSVGD is constructed with preconditioning operators; a preliminary theoretical study is given for the corresponding continuous equations. 
In Section 4, the algorithm is applied to solve an inverse problem governed by the steady state Darcy flow equation. The paper is concluded with some general remarks and directions for future work in Section 5.

\section{A short review of SVGD}\label{SectionReview}

Let $\mathcal{H}$ be a separable Hilbert space endowed with the Borel  $\sigma$-algebra $\mathcal{B}(\mathcal{H})$.
Denote by $\mathcal{G}$, $u$, and $\bm{d}$ the solution operator of some PDE, the model parameter, and the observation, respectively. We assume that $u\in\mathcal{H}$ and $\bm{d}\in\mathbb{R}^{N_d}$ with $N_d$ being a positive integer. 
The observation $\bm{d}$ is related to $\mathcal{G}(u)$ and the random noise $\bm{\epsilon}$ through some functions \cite{Kaipio2004Book}, e.g., the additive noise model or the multiplicative noise model. 
We refer to Section \ref{ApplicationSection} for a specific example.

For statistical inverse problems, it is usually required to find a probability measure $\mu^{\bm{d}}$ on $\mathcal{H}$,
which is known as the posterior probability measure and is specified by its density with respect to a prior probability measure $\mu_0$. The Bayesian formula on a Hilbert space is defined by
\begin{align}\label{BasicBayesFormula}
	\frac{d\mu^{\bm{d}}}{d\mu_0}(u) = \frac{1}{Z_{\bm{d}}}\exp\Big( -\Phi(u;\bm{d}) \Big),
\end{align}
where $\Phi\in C(\mathcal{H}\times\mathbb{R}^{N_d};\mathbb{R})$ 
and $\exp(-\Phi(u;\bm{d}))$ is integrable with respect to $\mu_0$.
The constant $Z_{\bm{d}}$ is chosen to ensure that $\mu^{\bm{d}}$ is indeed a probability measure. The prior measure $\mu_0:=\mathcal{N}(0,\mathcal{C}_0)$ is assumed to be a Gaussian measure defined on $\mathcal{H}$
with $\mathcal{C}_0$ being a self-adjoint, positive definite, and trace class operator.
Let $(\lambda_k,\varepsilon_k)_{k=1}^{\infty}$ be the eigensystem of $\mathcal{C}_0$ satisfying
$\mathcal{C}_0\varepsilon_k = \lambda_k^2\varepsilon_k$.
Denote by $P^N$ and  $Q^N$ the orthogonal projections of $\mathcal{H}$ onto
$X^N:=\text{span}\{ \varepsilon_1, \varepsilon_2, \ldots, \varepsilon_N \}$
and $X^{\perp}:=\text{span}\{ \varepsilon_{N+1}, \varepsilon_{N+2}, \ldots \}$, respectively. 
Clearly, we have $Q^N=\text{Id}-P^N$. Let $u^N:=P^N u\in X^N$ and $u^{\perp}:=Q^N u \in X^{\perp}$.
Define $\mathcal{C}_0^N=P^N\mathcal{C}_0 P^N$ and let $\mu_0^N=\mathcal{N}(0,\mathcal{C}_0^N)$
be a finite-dimensional Gaussian measure defined on $X^N$. Then an approximate measure $\mu^{dN}$ on $X^N$ can be defined by
\begin{align}\label{ApproBayesTheorem}
	\frac{d\mu^{\bm{d}N}}{d\mu_0^N}(u^N) = \frac{1}{Z_{\bm{d}}^N}\exp\Big( -\Phi(u^N; \bm{d}) \Big),
\end{align}
where
\begin{align*}
	Z_{\bm{d}}^N = \int_{X^N}\exp\Big( -\Phi(u^N; \bm{d}) \Big)\mu_0^N(du^N).
\end{align*}

Some more properties of the above approximate measure can be found in  \cite[Subsection 5.6]{Dashti2017}.
The probability measure $\mu^{\bm{d}N}$ can be written as the pushforward of the posterior measure $\mu^{\bm{d}}$ on $\mathbb{R}^N$, i.e., $\mu^{\bm{d}N} = P^{N}_{\#}\mu^{\bm{d}} := \mu^{\bm{d}}\circ (P^{N})^{-1}$.
Hence the measure $\mu^{\bm{d}N}$ has a Lebesgue density denoted by $p^{\bm{d}N}$ with the following form:
\begin{align}
	p^{\bm{d}N}(u^N) \propto \exp\Big(
	-\Phi(u^N; \bm{d}) - \frac{1}{2}\|u^N\|_{\mathcal{C}_0^N}^2
	\Big),
\end{align}
where $\|\cdot\|_{\mathcal{C}_0^N}$ represents $\|(\mathcal{C}_0^N)^{-1/2}\cdot\|_{\ell^2}$ with $\|\cdot\|_{\ell^2}$ standing for the usual $\ell^2$-norm.
Obviously, the target distribution $\mu^{\bm{d}N}$ is the solution to the optimization problem defined on the set $\mathcal{P}_2(\mathbb{R}^N)$
of probability measures $\nu$ such that $\int\|u^N\|^2d\nu^N(u^N) < \infty$ by:
\begin{align}\label{min1}
	\min_{\nu^N\in\mathcal{P}_2(\mathbb{R}^N)}\text{KL}(\nu^N || \mu^{\bm{d}N}),
\end{align}
where KL denotes the Kullback-Leibler (KL) divergence. 

Now, we present the Stein variational gradient descent (SVGD) algorithm. Denote $\text{KL}(\cdot || \mu^{\bm{d}N}) : \mathcal{P}_2(\mathbb{R}^N) \rightarrow [0,+\infty)$ as the functional $\nu^N \mapsto \text{KL}(\nu^N || \mu^{\bm{d}N})$. 
In order to obtain samples from $\mu^{\bm{d}N}$, the SVGD applies a gradient descent-like algorithm to the functional $\text{KL}(\cdot || \mu^{\bm{d}N})$. The standard gradient descent algorithm in the Wasserstein space applied to $\text{KL}(\cdot||\mu^{\bm{d}N})$, at each iteration $\ell \geq 0$, is 
\begin{align}\label{iter1}
	\nu_{\ell+1}^N = \left( \text{Id} - \epsilon\nabla\log\left( \frac{d\nu_{\ell}^N}{d\mu^{\bm{d}N}} \right) \right)_{\#}\nu_{\ell}^N,
\end{align}
where $\epsilon > 0$ is the step size.
This corresponds to a forward Euler discretization of the gradient flow of $\text{KL}(\cdot || \mu^{\bm{d}N})$
with respect to Stein geometry \cite{Duncan2019arXiv}. 
Instead of the Wasserstein gradient $\nabla\log\left(d\nu_{\ell}^N/d\mu^{\bm{d}N} \right)$ used in (\ref{iter1}), 
the SVGD uses $P_{\nu_{\ell}^N}\nabla\log\left( d\nu_{\ell}^N/d\mu^{\bm{d}N}\right)$ to generate the following iteration:
\begin{align}\label{iter2}
	\nu_{\ell+1}^N = \left( \text{Id} - \epsilon P_{\nu_{\ell}^N}\nabla\log\left( \frac{d\nu_{\ell}^N}{d\mu^{\bm{d}N}} \right) \right)_{\#}\nu_{\ell}^N,
\end{align}
where $P_{\nu_{\ell}^{N}}$ is the same as that in Subsection 3.1 of \cite{Korba2020NIPS}.
Let $\mathcal{H}_{K}^N$ be an $N$-dimensional reproducing kernel Hilbert space (RKHS) \cite{SteinwartBook} with the kernel function $K:\mathbb{R}^N\times \mathbb{R}^N \rightarrow \mathbb{R}$. 
To define $P_{\nu_{\ell}^{N}}$ rigorously, it is necessary to introduce the kernel integral operator based on the kernel function $K$, 
which will not be used in the rest of the paper. Hence, we omit it and refer to \cite{Korba2020NIPS} for the details. The reason for introducing the operator $P_{\nu_{\ell}^{N}}$ is that we have  
\begin{align}\label{inter3}
	P_{\nu_{\ell}^N}\nabla\log\left( \frac{d\nu_{\ell}^N}{d\mu^{\bm{d}N}}\right)(\cdot) = 
	-\mathbb{E}_{u^N\sim \nu_{\ell}^N}\left[ K(u^N,\cdot)\nabla_{u^N}\log p^{\bm{d}N}(u^N) + \nabla_{u^N}K(u^N,\cdot) \right]
\end{align}
under some mild conditions.
For every $\ell\geq 0$, let $u^{N,\ell}$ be distributed according to $\nu_{\ell}^N$. Using (\ref{iter2})--(\ref{inter3}), 
we obtain a particle update scheme
\begin{align}\label{iter4}
	u^{N,\ell+1} = u^{N,\ell} + \epsilon \phi_{\ell}^{N}(u^{N,\ell}),
\end{align}
where
\begin{align}\label{tran2for}
	\phi^{N*}_{\ell}(\cdot) = \mathbb{E}_{u^N\sim q_{\ell}^{N}}\Big[ 
	K(u^N,\cdot)\nabla_{u^N}\log p^{\bm{d}N}(u^N) + \nabla_{u^N}K(u^N,\cdot)
	\Big].
\end{align}

The basic SVGD algorithm is given in Algorithm \ref{finiteSVGD}. Inspired by applications in machine learning, the SVGD type algorithms have been widely studied over the last few years \cite{Detommaso2018NIPS,Duncan2019arXiv,Korba2020NIPS,Liu2017NIPS,Liu2016NIPS,Lu2019SIAM}.

\begin{algorithm}
	\setstretch{1}
	\caption{Finite-dimensional Stein variational gradient descent}
	\label{finiteSVGD}
	\begin{algorithmic}
		\STATE {\textbf{Input: }A target probability measure with density function $p^{\bm{d}N}(u^N)$ and a set of particles
			$\{u^{N,0}_i\}_{i=1}^{m}$.
		}
		\STATE {\textbf{Output: }A set of particles $\{u_i^N\}_{i=1}^m$ that approximates the target probability measure.
		}
		\STATE {\textbf{for }iteration $\ell$ do
			\vskip -0.5 cm
			\begin{align*}
				u_i^{N,\ell+1}\longleftarrow u_i^{N,\ell}+\epsilon_{\ell}\phi^*(u_{i}^{N,\ell}),
			\end{align*}
			\vskip -0.2 cm
			$\qquad$where 
			\begin{align*}
				\phi^*(u^N) = \frac{1}{m}\sum_{j=1}^m \Big[
				K(u_{j}^{N,\ell}, u^N)\nabla_{u_{j}^{N,\ell}}\log p^{\bm{d}N}(u_{j}^{N,\ell}) +
				\nabla_{u_{j}^{N,\ell}}K(u_{j}^{N,\ell}, u^N)
				\Big],
			\end{align*}
			\qquad and $\epsilon_{\ell}$ is the step size at the $\ell$-th iteration. \\
			\textbf{end for}
		}
	\end{algorithmic}
\end{algorithm}

\section{SVGD on separable Hilbert spaces}

This section is devoted to the construction of iSVGD and the preconditioning operators. The corresponding continuity equations are provided for a preliminary theoretical study of the method. 

\subsection{Hilbert scale and vector-valued RKHS}
For constructing iSVGD, we need to characterize the smoothness of functions that belong to some infinite dimensional spaces. The Sobolev spaces are usually employed to characterize the smoothness of functions. However, for presenting a general theory, we introduce the Hilbert scales defined by the prior covariance operator \cite{Engl1996Book}. The reason is that different covariance operators employed in practical problems lead to the same form of Hilbert scales. However, they are related to different Sobolev spaces. Hence, the same form of the general theory can be flexibly adapted to different practical problems.

Let $\mathcal{C}_0:\, \mathcal{H} \rightarrow \mathcal{H}$ be the covariance operator introduced in Section \ref{SectionReview}.
Denote by $\mathcal{D}(\mathcal{C}_0)$ and $\mathcal{R}(\mathcal{C}_0)$ the domain and range of $\mathcal{C}_0$, respectively.
Let $\mathcal{H} = \overline{\mathcal{R}(\mathcal{C}_0)}\oplus\mathcal{R}(\mathcal{C}_0)^{\perp} = \overline{\mathcal{R}(\mathcal{C}_0)}$ (the closure of $\mathcal{R}(\mathcal{C}_0)$). It is clear to note that $\mathcal{C}_0^{-1}$ is a densely defined, unbounded, symmetric and positive-definite operator in $\mathcal{H}$.
Let $\langle\cdot,\cdot\rangle_{\mathcal{H}}$ and $\|\cdot\|_{\mathcal{H}}$ be the inner product and norm defined on the Hilbert space $\mathcal{H}$, respectively.
Define the Hilbert scales $(\mathcal{H}^{t})_{t\in\mathbb{R}}$
with $\mathcal{H}^{t} := \overline{\mathcal{S}_{f}}^{\|\cdot\|_{\mathcal{H}^{t}}}$, where
\begin{align*}
	\mathcal{S}_{f} := \bigcap_{n=0}^{\infty}\mathcal{D}(\mathcal{C}_0^{-n}), \quad
	\langle u, v\rangle_{\mathcal{H}^{t}} := \langle\mathcal{C}_0^{-t/2}u, \mathcal{C}_0^{-t/2}v\rangle_{\mathcal{H}}, \quad
	\|u\|_{\mathcal{H}^{t}} := \left\| \mathcal{C}_0^{-t/2}u \right\|_{\mathcal{H}}.
\end{align*}
The norms defined above possess the following properties (cf. \cite[Proposition 8.19]{Engl1996Book}).

\begin{lemma}\label{L2Inter}
	Let $(\mathcal{H}^{t})_{t\in\mathbb{R}}$ be the Hilbert scale induced by the operator $\mathcal{C}_0$ given above.
	Then the following assertions hold:
	\begin{enumerate}
		\item Let $-\infty < s < t < \infty$. Then the space $\mathcal{H}^{t}$ is densely and continuously embedded into $\mathcal{H}^{s}$.
		\item If $t \geq 0$, then $\mathcal{H}^{t} = \mathcal{D}(\mathcal{C}_0^{-t/2})$, and $\mathcal{H}^{-t}$ is the dual space of $\mathcal{H}^{t}$.
		\item Let $-\infty < q < r < s < \infty$ then the interpolation inequality
		$\|u\|_{\mathcal{H}^{r}} \leq \|u\|_{\mathcal{H}^{q}}^{\frac{s-r}{s-q}} \|u\|_{\mathcal{H}^{s}}^{\frac{r-q}{s-q}}$ holds when $u \in \mathcal{H}^{s}$.
	\end{enumerate}
\end{lemma}

Now, we introduce some basic notations of vector-valued reproducing kernel Hilbert space (RKHS).
The following definition concerns the Hilbert space adjoint opertor \cite{Reed1980FunctionalAnalysis}. 

\begin{definition}\label{DefAdjBanach}
	Let $\mathcal{X}$ and $\mathcal{Y}$ be Banach spaces, and $T$ be a bounded linear operator from $\mathcal{X}$ to $\mathcal{Y}$.
	The \textbf{Banach space adjoint} of $T$, denoted by $T'$, is the bounded linear operator from $\mathcal{Y}^*$ to $\mathcal{X}^*$ and is defined by
	$(T'\ell)(u) = \ell(Tu)$
	for all $\ell\in\mathcal{Y}^*$, $u\in\mathcal{X}$. 
	Let $\mathcal{X}$ and $\mathcal{Y}$ be Hilbert spaces, and 
	$C_1:\mathcal{X}\rightarrow\mathcal{X}^*$ be the map that assigns to each $u\in\mathcal{X}$,
	the bounded linear functional $\langle u,\cdot\rangle_{\mathcal{X}}$ in $\mathcal{X}^*$.
	Let $C_2:\mathcal{Y}\rightarrow\mathcal{Y}^{*}$ be defined similarly as $C_1$.
	Then the \textbf{Hilbert space adjoint} of $T$ is a map $T^*:\mathcal{Y}\rightarrow\mathcal{X}$ given by
	$T^* = C_1^{-1}T'C_2$.
\end{definition}

Next, we introduce operator-valued positive definite kernels, which constitute the framework for specifying
vector-valued RKHS. Following Kadri et al. \cite{Kadri2016JMLR} to avoid topological and measurability issues, we focus on separable Hilbert spaces with reproducing operator-valued kernels whose elements are continuous functions. 
Denote by $\mathcal{X}$ and $\mathcal{Y}$ the separable Hilbert spaces and by $\mathcal{L}(\mathcal{X},\mathcal{Y})$ the set of
bounded linear operators from $\mathcal{X}$ to $\mathcal{Y}$. When $\mathcal{X}=\mathcal{Y}$, we write $\mathcal{L}(\mathcal{Y},\mathcal{Y})$ briefly as $\mathcal{L}(\mathcal{Y})$.

\begin{definition}\label{defOpValKernel}(Operator-valued kernels)
	An $\mathcal{L}(\mathcal{Y})$-valued kernel $\bm{K}$ on $\mathcal{X}\times\mathcal{X}$ is an operator
	$\bm{K}(\cdot,\cdot):\mathcal{X}\times\mathcal{X}\rightarrow\mathcal{L}(\mathcal{Y})$;
	\begin{enumerate}
		\item $\bm{K}$ is Hermitian if $\forall\, u,v\in\mathcal{X}$, $\bm{K}(u,v) = \bm{K}(v,u)^*$;
		\item $\bm{K}$ is nonnegative on $\mathcal{X}$ if it is Hermitian and for every natural number $r$ and all
		$\{(u_i, v_i)_{i=1,\ldots,r}\}\in\mathcal{X}\times\mathcal{Y}$, the matrix with $ij$-th entry
		$\langle \bm{K}(u_i,u_j)v_i, v_j\rangle_{\mathcal{Y}}$ is nonnegative (positive-definite).
	\end{enumerate}
\end{definition}


\begin{definition}\label{defRKHS}(Vector-valued RKHS)
	Let $\mathcal{X}$ and $\mathcal{Y}$ be separable Hilbert spaces.
	A Hilbert space $\mathcal{F}$ of operators from $\mathcal{X}$ to $\mathcal{Y}$ is called a reproducing kernel
	Hilbert space if there is a nonnegative $\mathcal{L}(\mathcal{Y})$-valued kernel
	$\bm{K}$ on $\mathcal{X}\times\mathcal{X}$ such that
	\begin{enumerate}
		\item the operator $v\longmapsto\bm{K}(u,v)g$ belongs to $\mathcal{F}$ for all
		$v, u\in\mathcal{X}$ and $g\in\mathcal{Y}$;
		\item for every $f\in\mathcal{F}$, $u\in\mathcal{X}$ and $g\in\mathcal{Y}$, we have
		$\langle f(u), g\rangle_{\mathcal{Y}} = \langle f(\cdot), \bm{K}(u,\cdot)g\rangle_{\mathcal{F}}$.
	\end{enumerate}
\end{definition}

Throughout the paper, we assume that the kernel $\bm{K}$ is \emph{locally bounded and separately continuous}, which guarantee that $\mathcal{F}$ is a subspace of $C(\mathcal{X},\mathcal{Y})$
(the vector space of continuous operators from $\mathcal{X}$ to $\mathcal{Y}$).
If the kernel $\bm{K}$ is nice enough \cite{Carmeli2006AA,Carmeli2010AA}, 
then it is the reproducing kernel of some Hilbert space $\mathcal{F}$.

Since the kernel is an important part of the SVGD, we provide some intuitive ideas about the operator-valued kernel. 
Let $u,v\in\mathcal{H}$ and $h>0$ be a positive constant. To construct the infinite-dimensional SVGD, 
we may introduce a scalar-valued kernel $K(u,v) := \exp\left( -\frac{1}{h}\|u - v\|_{\mathcal{H}}^2 \right)$ 
and consider the operator-valued kernel 
\begin{align}\label{2kernel1}
	\bm{K}(u, v) = K(u, v)\text{Id}.
\end{align}
For example, we can take $\mathcal{H} = L^2(\Omega)$ with $\Omega$ being a bounded open domain and have  
\begin{align}\label{kernelTerm1}
	\|u-v\|_{\mathcal{H}}^2 = \int_{\Omega} |u(x)-v(x)|^2 dx.
\end{align}
However, for solving inverse problems of PDEs, it is useful to introduce some preconditioning operators which require to consider operator-valued kernels. Here, we illustrate this by a simple example. 
Let the prior measure $\mu_0=\mathcal{N}(0,(\text{Id}-\Delta)^{-2})$, where $\Delta$ is the Dirichlet Laplace operator 
and $\mathcal{H}=L^2(\Omega)$. Intuitively we have $\mathcal{H}^{1} \approx H^{2}(\Omega)$,
where $H^2(\Omega)$ is the usual Sobolev space. 
By the theory of Gaussian measures \cite{Prato2006IDAnalysis}, we approximately have $\mu_0(H^2(\Omega)) = 0$ (not rigorously correct).
Inspired by the pCN algorithm \cite{Cotter2013SS}, we may choose the preconditioning operator $T=\text{Id}-\Delta$.
If we choose the Gaussian kernel as (\ref{2kernel1}), then the transformed kernel function becomes
\begin{align}
	\bm{K}(u, v) = \exp\left( -\frac{1}{h}\|T(u-v)\|_{L^2}^2 \right)T^{-1}(T^{-1})^*,
\end{align}
which is approximately equal to 
\begin{align}
	\bm{K}(u, v) \approx \exp\left( -\frac{1}{h}\|u-v\|_{H^2}^2 \right)(\text{Id}-\Delta)^{-2}.
\end{align}
Obviously, the kernel function equals to zero when $u-v$ does not belong to $H^2(\Omega)$, i.e., $\|u-v\|_{H^2} < \infty$ when
$u-v\in H^2(\Omega)$. Hence, the kernel function takes nonzero values and the algorithms can work only if the differences of any two particles reside in a measure zero set. 
In our opinion, this restriction seems too strong in the infinite-dimensional setting 
to make the particles over concentrated (see our numerical example in Section \ref{ApplicationSection} to demonstrate this in details).

Based on the above discussion, we may introduce a parameter $s$ and have an approximate transformed kernel
\begin{align}
	\bm{K}(u, v) \approx \exp\left( -\frac{1}{h}\|u-v\|_{H^{2-2s}}^2 \right)(\text{Id}-\Delta)^{-2}.
\end{align}
However, to achieve this, we should not choose the original kernel (the kernel is not transformed by the operator $T$) to be the usual scalar-valued kernel. The original kernel may be chosen as
$\bm{K}_0(u, v) = K_0(u,v)(\text{Id}-\Delta)^{-2s}$,
where $K_0(u, v) := e^{-\frac{1}{h}\|u-v\|_{L^2}}$ with $h>0$ being a positive constant. In this setting, the preconditioning operator can be chosen as $T:=(\text{Id}-\Delta)^{1-s}$.
These intuitive ideas indicate that it is necessary to construct the infinite-dimensional SVGD based on the more involved operator-valued kernel theory. 

\subsection{iSVGD}\label{ConstructSVGDSection_const}

In this subsection, we present an infinite-dimensional version of the SVGD, i.e., iSVGD. 
For a function $u$, denote by $D_{u}$ and $D_{u_k}$ the Fr\'{e}chet derivative and the directional derivative in the $k$th direction, respectively. For simplicity of notation, we shall use $D$ and $D_k$ instead of $D_u$ and $D_{u_k}$, and write $\Phi(u; \bm{d})$ as $\Phi(u)$. Let 
\begin{align}\label{defV}
	V(u) = \Phi(u) + \frac{1}{2}\|u\|_{\mathcal{H}^{1}}^{2}, 
\end{align}
where the potential functional $\Phi$ is required to satisfy the following assumptions.  

\begin{assumption}\label{assumPhi}
	Let $\mathcal{X}$ and $\mathcal{H}$ be two separable Hilbert spaces. For $s\in[0,1]$, we assume
	$\mathcal{H}^{1-s}\subset\mathcal{X}\subset\mathcal{H}$.
	Let $M_1\in\mathbb{R}^{+}$ be a positive constant.
	For each $u\in\mathcal{X}\subset\mathcal{H}$, we introduce $D\Phi: \mathcal{X}\rightarrow \mathcal{X}^*$ and 
	$D^2\Phi:\mathcal{X}\rightarrow\mathcal{L}(\mathcal{X}, \mathcal{X}^*)$, then 	
	the functional $\Phi:\mathcal{X}\rightarrow\mathbb{R}$ satisfies
	\begin{align*}
		-M_1 \leq \Phi(u) &\leq M_2 (\|u\|_{\mathcal{X}}), \\
		\|D\Phi(u)\|_{\mathcal{X}^*} &\leq M_3 (\|u\|_{\mathcal{X}}), \\
		\|D^2\Phi(u)\|_{\mathcal{L}(\mathcal{X},\mathcal{X}^*)} &\leq M_4(\|u\|_{\mathcal{X}}),
	\end{align*}
	where $M_2(\cdot)$, $M_3(\cdot)$, and $M_4(\cdot)$ are some monotonic non-decreasing functions.
\end{assumption}

The above assumption is a local version of \cite[Assumption 4]{Dashti2017}, which can be verified for many problems,
e.g., the Darcy flow model (Theorem \ref{ConditionVerifyThm} in Section \ref{ApplicationSection}).
We now optimize $\phi$ in the unit ball of a general vector-valued RKHS $\mathcal{H}_{\bm{K}}$ with an operator
valued kernel $\bm{K}(u, u')\in\mathcal{L}(\mathcal{Y})$:
\begin{align}\label{gSo1}
	\phi_{\bm{K}}^* = \argmax_{\phi\in\mathcal{H}_{\bm{K}}}\left\{
	\mathbb{E}_{u\sim\mu}[\mathcal{S}\phi(u)], \,\, \text{s.t. }\|\phi\|_{\mathcal{H}_{K}}\leq 1
	\text{ and }D\phi:\mathcal{X}\rightarrow\mathcal{L}_1(\mathcal{X},\mathcal{Y})
	\right\},
\end{align}
where $\mathcal{S}$ is the generalized Stein operator defined formally as follows:
\begin{align}\label{gSo2}
	\mathcal{S}\phi(u) = -\langle DV(u), \phi(u) \rangle_{\mathcal{Y}} +
	\sum_{k=1}^{\infty}D_{k}\langle\phi(u), e_k\rangle_{\mathcal{Y}},
\end{align}
and $\mathcal{L}_1(\mathcal{X},\mathcal{Y})$ denotes the set of all trace class operators from $\mathcal{X}$ to $\mathcal{Y}$.
For the convergence of the infinite sum, we illustrate it in Theorem \ref{limitDefineApp}.
Here, $\{e_k\}_{k=1}^{\infty}$ stands for an orthonormal basis of space $\mathcal{Y}$
and $\mu$ is a probability measure defined on $\mathcal{H}$. Moreover,  we assume that $\phi:\mathcal{X}\rightarrow\mathcal{Y}$ is Fr\'{e}chet differentiable, and the derivative is continuous to ensure the validity of (\ref{gSo1}). 

\begin{remark}
	In the finite-dimensional case, the operator $D\phi(u)$ naturally belongs to $\mathcal{L}_1(\mathcal{X},\mathcal{Y})$
	(cf.  \cite[Appendix C]{DaPrato1996Book}).
\end{remark}

The following assumption is also needed for the operator-valued kernels, which include many useful kernels,
e.g., the radial basis function (RBF) kernel.

\begin{assumption}\label{assKernel}
	Let $\mathcal{X}$, $\mathcal{Y}$, and $\mathcal{H}$ be three separable Hilbert spaces. For $s \in [0,1]$, we assume that
	$\mathcal{H}^{-s-1}\subset\mathcal{Y}$ and
	\begin{align}\label{kernelCond}
		\sup_{u\in\mathcal{X}}\|\bm{K}(u,u)\|_{\mathcal{L}(\mathcal{Y})} < \infty.
	\end{align}
\end{assumption}

\begin{remark}\label{necessaryOpKernel}
	We mention that Condition (\ref{kernelCond}) holds for the bounded scalar-valued kernel functionals since a scalar-valued kernel functional can be seen as a scalar-valued kernel functional composite with an identity operator as demonstrated in (\ref{2kernel1}).
\end{remark}

To illustrate (\ref{gSo1}) and (\ref{gSo2}), we prove Theorem \ref{limitDefineApp}.
For each particle $u$, we assume that $u\in\mathcal{H}^{1-s}$,
which is based on the following two considerations:
\begin{itemize}
	\item The SVGD with one particle is an optimization algorithm for finding maximum a posterior (MAP) estimate.
	The MAP estimate belongs to the separable Hilbert space $\mathcal{H}^{1}$.
	\item For the prior probability measure, the space $\mathcal{H}^{1}$ has zero measure \cite{DaPrato1996Book}.
	Intuitively, if all particles belong to $\mathcal{H}^{1}$, the particles tend to concentrate
	around a small set that leads to unreliable estimates of statistical quantities. Hence, we may assume that the particles
	belong to a larger space containing $\mathcal{H}^{1}$.
\end{itemize}
\begin{theorem}\label{limitDefineApp}
	The generalized Stein operator (\ref{gSo2}) defined on $\mathcal{Y}$ can be obtained
	by taking $N\rightarrow\infty$ in the following finite-dimensional Stein operator:
	\begin{align}
		\mathcal{S}^N\phi^N(u^N) = -\langle DV(u^N), \phi^N(u^N) \rangle_{\mathcal{Y}} +
		\sum_{k=1}^N D_{k}\langle \phi^N(u^N), e_k\rangle_{\mathcal{Y}},
	\end{align}
	where $\phi^N := P^N\circ \phi$.
\end{theorem}
\begin{proof}
	By straightforward calculations, we have
	\begin{align}
		\begin{split}
			\mathcal{S}\phi(u) - \mathcal{S}^N\phi^N(u^N) =
			& - \Big(\langle DV(u), \phi(u) \rangle_{\mathcal{Y}} -
			\langle DV(u^N), \phi^N(u^N) \rangle_{\mathcal{Y}} \Big) \\
			& + \left(\sum_{k=1}^{\infty}D_{k}\langle\phi(u), e_k\rangle_{\mathcal{Y}} -
			\sum_{k=1}^N D_{k}\langle \phi^N(u^N), e_k\rangle_{\mathcal{Y}}\right) \\
			= & - \text{I} + \text{II}.
		\end{split}
	\end{align}
	For term $\text{I}$, we have
	\begin{align}
		\begin{split}
			\text{I} = &
			\langle D(V(u)-V(u^N)), \phi^N(u^N) \rangle_{\mathcal{Y}} +
			\langle DV(u), \phi(u)-\phi^N(u^N) \rangle_{\mathcal{Y}}\\
			= & \, \text{I}_1(N) + \text{I}_2(N).
		\end{split}
	\end{align}
	For term $\text{I}_1(N)$, we find that
	\begin{align}\label{limitTheo1}
		\begin{split}
			\text{I}_1(N) &\! = \!\langle D(\Phi(u)\!-\!\Phi(u^N)), \phi^N(u^N) \rangle_{\mathcal{Y}}
			\! + \! \langle \mathcal{C}_0^{-1/2}(u\!-\!u^N), \mathcal{C}_0^{-1/2}\phi^N(u^N) \rangle_{\mathcal{Y}},
		\end{split}
	\end{align}
	where the second term on the right-hand side is understood as the white noise mapping \cite{Prato2006IDAnalysis}.
	According to Assumptions \ref{assumPhi} and \ref{assKernel}, we know that
	\begin{align}
		\begin{split}
			\lim_{N\rightarrow\infty}\|D(\Phi(u) - \Phi(u^N))\|_{\mathcal{Y}} \leq &
			\lim_{N\rightarrow\infty}C\|D(\Phi(u)-\Phi(u^N))\|_{\mathcal{H}^{-1-s}} \\
			\leq &
			\lim_{N\rightarrow\infty}C\|D(\Phi(u)-\Phi(u^N))\|_{\mathcal{H}^{-1+s}} \\
			\leq &
			\lim_{N\rightarrow\infty}CM_4(2\|u\|_{\mathcal{X}})\|u-u^N\|_{\mathcal{H}^{1-s}} = 0,
		\end{split}
	\end{align}
	where $C$ is a generic constant that can be different from line to line.
	Hence, we obtain
	\begin{align}\label{limitTheo2}
		\lim_{N\rightarrow\infty}\langle D(\Phi(u)-\Phi(u^N)), \phi^N(u^N) \rangle_{\mathcal{Y}} = 0.
	\end{align}
	Taking $u_m\in\mathcal{H}^2$ such that $u_m \rightarrow u$ in $\mathcal{H}^{1-s}$, 
	we have
	\begin{align*}
		\langle \mathcal{C}_0^{-1/2}(u-u^N), \mathcal{C}_0^{-1/2}\phi^N(u^N) \rangle_{\mathcal{Y}} = &
		\lim_{m\rightarrow\infty}\langle \mathcal{C}_0^{-1/2}(u_m-u_m^N), \mathcal{C}_0^{-1/2}\phi^N(u^N) \rangle_{\mathcal{Y}}\\
		= & \lim_{m\rightarrow\infty}\langle P^N\mathcal{C}_0^{-1}(u_m-u_m^N), \phi(u^N) \rangle_{\mathcal{Y}} \\
		= & \lim_{m\rightarrow\infty}\langle \phi(\cdot), \bm{K}(u^N,\cdot)P^N\mathcal{C}_0^{-1}(u_m-u_m^N) \rangle_{\mathcal{H}_{\bm{K}}}.
	\end{align*}
	As for the last term in the above equality, we have the following estimates:
	\begin{align*}
		& \langle \phi(\cdot), \bm{K}(u^N,\cdot)P^N\mathcal{C}_0^{-1}(u_m-u_m^N) \rangle_{\mathcal{H}_{\bm{K}}} \leq \\
		& \qquad
		\langle \phi(\cdot), \phi(\cdot) \rangle_{\mathcal{H}_{\bm{K}}}
		\langle \bm{K}(u^N,\cdot)P^N\mathcal{C}_0^{-1}(u_m-u_m^N),\bm{K}(u^N,\cdot)P^N\mathcal{C}_0^{-1}(u_m-u_m^N) \rangle_{\mathcal{H}_{\bm{K}}} \\
		& \quad \leq
		\langle \phi(\cdot), \phi(\cdot) \rangle_{\mathcal{H}_{\bm{K}}}
		\langle P^N\bm{K}(u^N,u^N)P^N\mathcal{C}_0^{-1}(u_m-u_m^N), \mathcal{C}_0^{-1}(u_m-u_m^N) \rangle_{\mathcal{Y}} \\
		& \quad \leq C
		\langle \phi(\cdot), \phi(\cdot) \rangle_{\mathcal{H}_{\bm{K}}}
		\|\mathcal{C}_0^{-1}(u_m-u_m^N)\|_{\mathcal{Y}}^2  \\
		& \quad \leq C
		\langle \phi(\cdot), \phi(\cdot) \rangle_{\mathcal{H}_{\bm{K}}}
		\|\mathcal{C}_0^{-\frac{1-s}{2}}(u_m-u_m^N)\|_{\mathcal{H}}^2.
	\end{align*}
	Replacing $u_m - u_m^N$ by $(u_m - u_m^N) - (u - u^N)$,
	we deduce 
	\begin{align}\label{limitTheo4}
		\langle \mathcal{C}_0^{-1/2}(u-u^N), \mathcal{C}_0^{-1/2}\phi^N(u^N) \rangle_{\mathcal{Y}}
		= & \lim_{m\rightarrow\infty}\langle \phi(\cdot), \bm{K}(u^N,\cdot)P^N\mathcal{C}_0^{-1}(u_m-u_m^N) \rangle_{\mathcal{H}_{\bm{K}}}
		\nonumber \\
		= & \langle \phi(\cdot), \bm{K}(u^N,\cdot)P^N\mathcal{C}_0^{-1}(u-u^N) \rangle_{\mathcal{H}_{\bm{K}}}.
	\end{align}
	Hence, we obtain
	\begin{align}\label{limitTheo3}
		\begin{split}
			& \lim_{N\rightarrow\infty}\langle \mathcal{C}_0^{-1/2}(u-u^N), \mathcal{C}_0^{-1/2}\phi^N(u^N) \rangle_{\mathcal{Y}}  \\
			= & \lim_{N\rightarrow\infty}\langle \phi(\cdot), \bm{K}(u^N,\cdot)P^N\mathcal{C}_0^{-1}(u-u^N) \rangle_{\mathcal{H}_{\bm{K}}} \\
			\leq & \lim_{N\rightarrow\infty}\langle\phi(\cdot),\phi(\cdot)\rangle_{\mathcal{H}_{\bm{K}}}
			\langle P^N\bm{K}(u^N\!,u^N)P^N\mathcal{C}_0^{-1}(u-u^N),
			\mathcal{C}_0^{-1}(u-u^N)\rangle_{\mathcal{Y}}	\\
			\leq & C \langle\phi(\cdot),\phi(\cdot)\rangle_{\mathcal{H}_{\bm{K}}} \lim_{N\rightarrow\infty}
			\|\mathcal{C}_0^{-\frac{1-s}{2}}(u-u^N)\|_{\mathcal{H}}^2 = 0.
		\end{split}
	\end{align}
	Plugging (\ref{limitTheo2}) and (\ref{limitTheo3}) into (\ref{limitTheo1}), we arrive at
	$\lim_{N\rightarrow\infty}I_1(N) = 0$.
	For term $\text{I}_2(N)$, it can be decomposed as follows:
	\begin{align}
		\text{I}_2(N) = \langle D\Phi(u), \phi(u)-\phi^N(u^N)\rangle_{\mathcal{Y}} \! + \!
		\langle \mathcal{C}_0^{-1/2}u, \mathcal{C}_0^{-1/2}(\phi(u)-\phi^N(u^N)) \rangle_{\mathcal{Y}}.
	\end{align}
	It follows from the continuity of $\phi$ that we have $\lim_{N\rightarrow\infty}\langle D\Phi(u), \phi(u)-\phi^N(u^N)\rangle_{\mathcal{Y}}=0$.
	Using similar estimates as those for deriving (\ref{limitTheo4}), we obtain
	\begin{align}\label{limitTheo5}
		\begin{split}
			& \langle \mathcal{C}_0^{-1/2}u, \mathcal{C}_0^{-1/2}(\phi(u)-\phi^N(u^N)) \rangle_{\mathcal{Y}} \\
			& \qquad\qquad
			= \langle \phi(\cdot), \bm{K}(u,\cdot)\mathcal{C}_0^{-1}u\rangle_{\mathcal{H}_{\bm{K}}}
			- \langle \phi(\cdot), \bm{K}(u^N,\cdot)P^N\mathcal{C}_0^{-1}u\rangle_{\mathcal{H}_{\bm{K}}}.
		\end{split}
	\end{align}
	By the continuity of $\bm{K}(\cdot,\cdot)$, we obtain
	\begin{align}
		\lim_{N\rightarrow\infty}\langle \mathcal{C}_0^{-1/2}u, \mathcal{C}_0^{-1/2}(\phi(u)-\phi^N(u^N)) \rangle_{\mathcal{Y}} = 0.
	\end{align}
	Now, we conclude that
	$\lim_{N\rightarrow\infty}\text{I}_2(N) = 0.$
	For term $\text{II}$, we have
	\begin{align}\label{limitTheo6}
		\text{II} = \sum_{k=1}^N D_{k}\langle \phi(u)-\phi(u^N), e_k\rangle_{\mathcal{Y}} +
		\sum_{k=N+1}^{\infty}D_{k}\langle \phi(u), e_k \rangle_{\mathcal{Y}}.
	\end{align}
	Let $\{\varphi_{k}\}_{k=1}^{\infty}$ be an orthonormal basis in $\mathcal{X}$, and then we have
	\begin{align}
		\sum_{k=N+1}^{\infty}D_{k}\langle \phi(u), e_k \rangle_{\mathcal{Y}} =
		\sum_{k=N+1}^{\infty}\langle D\phi(u)\varphi_{k}, e_k\rangle_{\mathcal{Y}}\rightarrow 0 \quad
		\text{as }N\rightarrow\infty,
	\end{align}
	where we use the condition $D\phi(u)\in\mathcal{L}_1(\mathcal{X},\mathcal{Y})$.
	For the first term on the right-hand side of (\ref{limitTheo6}), we find that
	\begin{align}
		\sum_{k=1}^N D_{k}\langle \phi(u)-\phi(u^N), e_k\rangle_{\mathcal{Y}} =
		\sum_{k=1}^N \langle (D\phi(u)-D\phi(u^N))\varphi_{k}, e_k\rangle_{\mathcal{Y}}.
	\end{align}
	Due to the continuity of the Fr\'{e}chet derivative of $\phi$, we know that the above summation goes to $0$ as $N\rightarrow\infty$.
	Combining the estimates of $\text{I}$ and $\text{II}$, we complete the proof.
\end{proof}

The following theorem gives explicitly the iSVGD update directions that are essential for the construction of iSVGD.
\begin{theorem}\label{optimalSolTheorem}
	Let $\bm{K}(\cdot,\cdot):\mathcal{X}^{2}\rightarrow\mathcal{L}(\mathcal{Y})$ be a positive definite kernel that is
	Fr\'{e}chet differentiable on both variables. In addition, we assume that
	\begin{align}\label{optimalAssum}
		\mathbb{E}_{u\sim\mu}\Big[
		D_{u'}\bm{K}(u,u')\mathcal{C}_0^{-1/2}g+\sum_{k=1}^{\infty}D_{u_k}D_{u'}\bm{K}(u,u')e_k
		\Big]
	\end{align}
	belongs to $\mathcal{L}_1(\mathcal{X},\mathcal{Y})$ for each $u'\in\mathcal{X}$ and $g\in\mathcal{H}^{-s}$.
	Then, the optimal $\phi_{\bm{K}}^*$ in (\ref{gSo1}) is
	\begin{align}\label{transITH}
		\phi^*_{\bm{K}}(\cdot)\propto \mathbb{E}_{u\sim\mu}\Big[
		\bm{K}(u,\cdot)(-D\Phi(u)-\mathcal{C}_0^{-1}u)+\sum_{k=1}^{\infty}D_{u_k}\bm{K}(u,\cdot)e_k
		\Big],
	\end{align}
	where $\{e_k\}_{k=1}^{\infty}$ is an orthonormal basis of $\mathcal{Y}$ and the term $\bm{K}(u,\cdot)\mathcal{C}_0^{-1}u$
	is understood in the following limiting sense:
	\begin{align}
		\bm{K}(u,\cdot)\mathcal{C}_0^{-1}u := \lim_{m\rightarrow\infty}\bm{K}(u,\cdot)\mathcal{C}_0^{-1}u_m.
	\end{align}
	Here the limit is taken in $\mathcal{H}_{\bm{K}}$ and $\{u_{m}\}_{m=1}^{\infty}\subset\mathcal{H}^2$ such that
	$\|\mathcal{C}_0^{-\frac{1-s}{2}}(u_m-u)\|_{\mathcal{H}}\rightarrow0$ as $m\rightarrow\infty$.
\end{theorem}
\begin{proof}
	First, by taking $\phi(u)$ as an element in $\mathcal{H}_{\bm{K}}$, we have
	\begin{align}\label{optimalThe1}
		\langle DV(u), \phi(u)\rangle_{\mathcal{Y}} = \langle D\Phi(u), \phi(u) \rangle_{\mathcal{Y}}
		+ \langle \mathcal{C}_0^{-1/2}u, \mathcal{C}_0^{-1/2}\phi(u) \rangle_{\mathcal{Y}} = \text{I} + \text{II},
	\end{align}
	where term II is understood as the white noise mapping. For term I, we have
	\begin{align}\label{optimalThe4}
		\text{I} = \langle \phi(\cdot), \bm{K}(u,\cdot)D\Phi(u) \rangle_{\mathcal{H}_{\bm{K}}},
	\end{align}
	where the proposition (2) in Definition \ref{defRKHS} is employed. For term II, we take $u_m\in\mathcal{H}^2$ such that
	$\lim_{m\rightarrow\infty}\|\mathcal{C}_0^{-\frac{1-s}{2}}(u_m - u)\|_{\mathcal{H}} = 0.$
	It is clear to note that 
	\begin{align}\label{optimalThe2}
		\langle \mathcal{C}_0^{-1/2}u_m, \mathcal{C}_0^{-1/2}\phi(u) \rangle_{\mathcal{Y}} =
		\langle \mathcal{C}_0^{-1}u_m, \phi(u) \rangle_{\mathcal{Y}} =
		\langle \phi(\cdot), \bm{K}(u, \cdot)\mathcal{C}_0^{-1}u_m \rangle_{\mathcal{H}_{\bm{K}}}.
	\end{align}
	Because
	\begin{align*}
		& |\langle \phi(\cdot),\bm{K}(u,\cdot)\mathcal{C}_0^{-1}u_m \rangle_{\mathcal{H}_{\bm{K}}} -
		\langle \phi(\cdot), \bm{K}(u,\cdot)\mathcal{C}_0^{-1}u \rangle_{\mathcal{H}_{\bm{K}}} |^2  \\
		\leq & \langle \phi(\cdot), \phi(\cdot) \rangle_{\mathcal{H}_{\bm{K}}}
		\langle \bm{K}(u,\cdot)\mathcal{C}_0^{-1}(u_m-u),
		\bm{K}(u,\cdot)\mathcal{C}_0^{-1}(u_m-u) \rangle_{\mathcal{H}_{\bm{K}}}	\\
		= & \langle \phi(\cdot), \phi(\cdot) \rangle_{\mathcal{H}_{\bm{K}}}
		\langle \bm{K}(u,u)\mathcal{C}_0^{-1}(u_m-u),\mathcal{C}_0^{-1}(u_m-u) \rangle_{\mathcal{Y}} \\
		\leq & \langle \phi(\cdot), \phi(\cdot) \rangle_{\mathcal{H}_{\bm{K}}}
		\langle \bm{K}(u,u)\mathcal{C}_0^{-1}(u_m-u),
		\mathcal{C}_0^{-1}(u_m-u) \rangle_{\mathcal{Y}}	\\
		\leq & C \langle \phi(\cdot), \phi(\cdot) \rangle_{\mathcal{H}_{\bm{K}}} \|\mathcal{C}_0^{-\frac{1-s}{2}}(u_m-u)\|_{\mathcal{H}}^2,
	\end{align*}
	we find that
	$\lim_{m\rightarrow\infty}\langle \phi(\cdot), \bm{K}(u,\cdot)\mathcal{C}_0^{-1}u_m\rangle_{\mathcal{H}_{\bm{K}}} =
	\langle \phi(\cdot), \bm{K}(u,\cdot)\mathcal{C}_0^{-1}u \rangle_{\mathcal{H}_{\bm{K}}}.$
	Hence, let $m\rightarrow\infty$ in (\ref{optimalThe2}), we have
	\begin{align}\label{optimalThe3}
		\langle \mathcal{C}_0^{-1/2}u, \mathcal{C}_0^{-1/2}\phi(u) \rangle_{\mathcal{Y}} =
		\langle \phi(\cdot), \bm{K}(u,\cdot)\mathcal{C}_0^{-1}u \rangle_{\mathcal{H}_{\bm{K}}}.
	\end{align}
	Plugging (\ref{optimalThe3}) and (\ref{optimalThe4}) into (\ref{optimalThe1}), we obtain
	\begin{align}\label{optimalThe5}
		\begin{split}
			\langle DV(u), \phi(u) \rangle_{\mathcal{Y}} = &
			\langle \phi(\cdot), \bm{K}(u,\cdot)D\Phi(u) + \bm{K}(u,\cdot)\mathcal{C}_0^{-1}u \rangle_{\mathcal{H}_{\bm{K}}} \\
			= & \langle \phi(\cdot), \bm{K}(u,\cdot)DV(u) \rangle_{\mathcal{H}_{\bm{K}}}.
		\end{split}
	\end{align}
	
	Next, let us calculate the second term on the right-hand side of (\ref{gSo2}). A simple calculation yields  
	\begin{align}
		\sum_{k=1}^{\infty}D_{k}\langle \phi(u), e_k \rangle_{\mathcal{Y}} =
		\sum_{k=1}^{\infty}D_{k}\langle \phi(\cdot), \bm{K}(u,\cdot)e_k \rangle_{\mathcal{H}_{\bm{K}}}.
	\end{align}
	Since
	\begin{align}
		\begin{split}
			D_{k}\langle\phi(\cdot), \bm{K}(u,\cdot)e_k\rangle_{\mathcal{H}_{\bm{K}}} & = \lim_{\epsilon\rightarrow 0}
			\frac{1}{\epsilon}\langle\phi(\cdot),\bm{K}(u+\epsilon\varphi_k,\cdot)e_k - \bm{K}(u,\cdot)e_k\rangle_{\mathcal{H}_{\bm{K}}} \\
			& = \langle\phi(\cdot),D_{k}\bm{K}(u,\cdot)e_k\rangle_{\mathcal{H}_{\bm{K}}},
		\end{split}
	\end{align}
	we have
	\begin{align}\label{optimalThe6}
		\sum_{k=1}^{\infty}D_{k}\langle\phi(u),e_k\rangle_{\mathcal{Y}} =
		\Big\langle \phi(\cdot),\sum_{k=1}^{\infty}D_{k}\bm{K}(u,\cdot)e_k \Big\rangle_{\mathcal{H}_{\bm{K}}}.
	\end{align}
	Combining (\ref{optimalThe5}) and (\ref{optimalThe6}) with (\ref{gSo2}), we obtain
	\begin{align}
		\mathcal{S}\phi(u) = \Big\langle \phi(\cdot), -\bm{K}(u,\cdot)DV(u)
		+ \sum_{k=1}^{\infty}D_{k}\bm{K}(u,\cdot)e_k \Big\rangle_{\mathcal{H}_{\bm{K}}}.
	\end{align}
	Thus, the optimization problem (\ref{gSo1}) possesses a solution $\phi^*_{\bm{K}}(\cdot)$ satisfying
	\begin{align}
		\phi^*_{\bm{K}}(\cdot) \propto \mathbb{E}_{u\sim\mu}\Big[
		-\bm{K}(u,\cdot)DV(u) + \sum_{k=1}^{\infty}D_{k}\bm{K}(u,\cdot)e_k
		\Big].
	\end{align}
	Based on condition (\ref{optimalAssum}), we know that $D\phi^*_{\bm{K}}(u)$ belongs to $\mathcal{L}_1(\mathcal{X},\mathcal{Y})$
	for each $u\in\mathcal{X}$, which completes the proof.
\end{proof}

\begin{remark}
	The optimal $\phi_{\bm{K}}^*$ is given in (\ref{transITH}) which is consistent with the finite-dimensional case. 
	Since the first and second terms on the right-hand side of (\ref{transITH}) are similar, we may just focus on the second 
	term which is usually named as the repulsive force term. 
	For each $u, v\in\mathcal{X}$, consider $\bm{K}(u, v):=K(u,v)\text{Id}$ with 
	$K(u,v):=\exp\left( -\frac{1}{h}\|u-v\|_{\mathcal{X}}^2 \right)$. Then, we have 
	\begin{align}\label{compareIDandFD}
		\begin{split}
			\sum_{k=1}^{\infty}D_{u_k}\bm{K}(u,v)e_k = & \sum_{k=1}^{\infty}\langle D_uK(u,v)e_k, \varphi_k \rangle_{\mathcal{X}} \\
			= & \sum_{k=1}^{\infty} -\frac{2}{h}\langle u-v, \varphi_k \rangle_{\mathcal{X}}K(u,v)e_k.
		\end{split}
	\end{align} 
	Projecting (\ref{compareIDandFD}) on one particular coordinate $e_{\ell}$ with $\ell\in\mathbb{N}$, we obtain 
	\begin{align}\label{eachSecondTerm}
		\begin{split}
			\left( \sum_{k=1}^{\infty}D_{u_k}\bm{K}(u,v)e_k \right)_{\ell} = & \left\langle
			\sum_{k=1}^{\infty} -\frac{2}{h}\langle u-v, \varphi_k \rangle_{\mathcal{X}}K(u,v)e_k, e_{\ell} \right\rangle_{\mathcal{Y}} \\
			= & -\frac{2}{h}\langle u-v, \varphi_{\ell} \rangle_{\mathcal{X}}K(u, v),
		\end{split}
	\end{align}
	which is similar to the $\ell$th coordinate of $\nabla_{u^N}K(u^N,v^N)$ appearing in (\ref{tran2for}).
	Additionally, we mention that the assumption (\ref{optimalAssum}) given in Theorem \ref{optimalSolTheorem} can be verified for many useful kernels. Detailed illustrations are provided in the supplementary material.
\end{remark}

By Theorem \ref{optimalSolTheorem}, we can construct a series of transformations as follows:
\begin{align}\label{algformul1}
	T_{\ell}(u) = u + \epsilon_{\ell}\mathbb{E}_{u'\sim\mu_{\ell}}\Big[
	-\bm{K}(u',u)DV(u')+\sum_{k=1}^{\infty}D_{(u')_k}\bm{K}(u',u)e_k
	\Big]
\end{align}
with $\ell=1,2,\ldots$.
In practice, we draw a set of particles $\{u^0_i\}_{i=1}^m$ from some initial measure, and then iteratively update
the particles with an empirical version of the above transformation in which the expectation under $\mu_{\ell}$ is approximated by the empirical mean of particles $\{u_i^{\ell}\}_{i=1}^m$ at the $\ell$-th iteration. The 
iSVGD is summarized in Algorithm \ref{infiniteSVGD}.

\begin{algorithm}
	\setstretch{1}
	\caption{Infinite-dimensional Stein variational gradient descent (iSVGD)}
	\label{infiniteSVGD}
	\begin{algorithmic}
		\STATE {\textbf{Input: }A target probability measure $\mu^{\bm{d}}$ that is absolutely
			continuous w.r.t the Gaussian measure $\mu_0=\mathcal{N}(0,\mathcal{C}_0)$ with
			$\frac{d\mu^{\bm{d}}}{d\mu_0}(u)\propto\exp(-\Phi(u))$
			and a set of particles $\{u^{0}_i\}_{i=1}^{m}$.
		}
		\STATE {\textbf{Output: }A set of particles $\{u_i\}_{i=1}^m$ that approximates the target probability measure.
		}
		\STATE {\textbf{for }iteration $\ell$ do
			\vskip -0.5 cm
			\begin{align*}
				u_i^{\ell+1}\longleftarrow u_i^{\ell}+\epsilon_{\ell}\phi^*(u_{i}^{\ell}),
			\end{align*}
			\vskip -0.2 cm
			$\qquad$where
			\vskip -0.5 cm
			\begin{align*}
				\phi^*(u) = \frac{1}{m}\sum_{j=1}^m \Big[
				\bm{K}(u_{j}^{\ell}, u)(-D\Phi(u_j^{\ell})-\mathcal{C}_0^{-1}u_j^{\ell}) +
				\sum_{k=1}^{\infty}D_{(u_{j}^{\ell})_k}\bm{K}(u_{j}^{\ell},u)e_k
				\Big].
			\end{align*}
			\vskip -0.2 cm
			\textbf{end for}
		}
	\end{algorithmic}
\end{algorithm}

\subsection{iSVGD with precondition information}\label{SVGDPre}

In the supplementary material, the numerical experiments indicate that 
the SVGD without preconditioning operators converges slowly for some inverse problems of PDEs. By the 
finite-dimensional SVGD \cite{Wang2019NIPS}, it may accelerate the convergence and give reliable estimates efficiently
by introducing preconditioning operators.
For constructing the iSVGD with preconditioning operators, let us begin with a theorem concerning the change of variables.

\begin{theorem}\label{ChangeVariableTheorem}
	Let $\mathcal{X}$ and $\mathcal{Y}$ be two separable Hilbert spaces, and let 
	$\mathcal{F}_0$ be a RKHS with a nonnegative $\mathcal{L}(\mathcal{Y})$-valued kernel
	$\bm{K}_0:\mathcal{X}\times\mathcal{X}\rightarrow\mathcal{L}(\mathcal{Y})$.
	Let $\tilde{\mathcal{X}}$ and $\tilde{\mathcal{Y}}$ be two separable Hilbert spaces, and $\mathcal{F}$ be
	the set of operators from $\tilde{\mathcal{X}}$ to $\tilde{\mathcal{Y}}$ given by
	\begin{align}\label{Cin1}
		\phi(u) = \bm{M}(u)\phi_0(t(u)) \quad \forall \,\,\phi_0\in\mathcal{F}_0,
	\end{align}
	where $\bm{M}:\tilde{\mathcal{X}}\rightarrow\mathcal{L}(\mathcal{Y},\tilde{\mathcal{Y}})$ is a fixed operator
	and is assumed to be an invertible operator for all $u\in\tilde{\mathcal{X}}$, and
	$t:\tilde{\mathcal{X}}\rightarrow\mathcal{X}$ is a fixed Fr\'{e}chet differentiable one-to-one mapping.
	For all $\phi,\phi'\in\mathcal{F}$, we can identify a unique $\phi_0, \phi'_0\in\mathcal{F}_0$ such that
	$\phi(u) = \bm{M}(u)\phi_0(t(u))$ and $\phi'(u) = \bm{M}(u)\phi'_0(t(u))$.
	Define the inner product on $\mathcal{F}$ via $\langle\phi,\phi'\rangle_{\mathcal{F}} = \langle \phi_0, \phi'_0\rangle_{\mathcal{F}_0}$,
	and then $\mathcal{F}$ is also a vector-valued RKHS, whose operator-valued kernel is
	\begin{align}\label{Cin2}
		\bm{K}(u, u') = \bm{M}(u')\bm{K}_0(t(u), t(u'))\bm{M}(u)^*,
	\end{align}
	where $\bm{M}(u)^*$ denotes the Hilbert space adjoint.
\end{theorem}

\begin{proof}
	Let $\{ (u_i, g_i)_{i=1,\ldots,N} \}\subset\tilde{\mathcal{X}}\times\tilde{\mathcal{Y}}$, and we have 
	\begin{align}\label{Cth1}
		\begin{split}
			\langle \bm{K}(u_i, u_j)g_i, g_j\rangle_{\tilde{\mathcal{Y}}} &
			= \langle \bm{M}(u_j)\bm{K}_0(t(u_i), t(u_j))\bm{M}(u_i)^* g_i, g_j\rangle_{\tilde{\mathcal{Y}}} \\
			& = \langle \bm{K}_0(t(u_i), t(u_j))\bm{M}(u_i)^* g_i, \bm{M}(u_j)^* g_{j}\rangle_{\mathcal{Y}}.
		\end{split}
	\end{align}
	Then, the nonnegativity of $\bm{K}(\cdot,\cdot)$ follows from the nonnegative property of $\bm{K}_0(\cdot,\cdot)$.
	To prove the theorem, it suffices to verify the two conditions shown in Definition \ref{defRKHS}.
	For every $u,v\in\tilde{\mathcal{X}}$ and $g\in\tilde{\mathcal{Y}}$, we consider the operator
	$f(v) = \bm{K}(u,v)g = \bm{M}(v)\bm{K}_0(t(u), t(v))\bm{M}(u)^* g$.
	Because of $\bm{M}(u)^* g\in\mathcal{Y}$, we easily obtain $$\bm{K}_0(t(u), t(v))\bm{M}(u)^* g\in\mathcal{F}_0.$$
	According to (\ref{Cin1}), we conclude that $f(\cdot)\in\mathcal{F}$.
	
	Next, let us verify the reproducing property of $\bm{K}(\cdot,\cdot)$.
	For every $u\in\tilde{\mathcal{X}}, g\in\tilde{\mathcal{Y}}$, and $\phi\in\mathcal{F}$, we have
	\begin{align*}
		\langle\phi(u), g\rangle_{\tilde{\mathcal{Y}}} & =
		\langle \bm{M}(u)\phi_0(t(u)), g\rangle_{\tilde{\mathcal{Y}}} =
		\langle \phi_0(t(u)), \bm{M}(u)^* g\rangle_{\mathcal{Y}} \\
		& = \langle\phi_0(\cdot), \bm{K}_0(t(u),\cdot)\bm{M}(u)^* g\rangle_{\mathcal{F}_0} \\
		& = \langle\bm{M}(\cdot)\phi_0(t(\cdot)), \bm{M}(\cdot)\bm{K}_0(t(u), t(\cdot))\bm{M}(u)^* g\rangle_{\mathcal{F}} \\
		& = \langle\phi(\cdot), \bm{K}(u,\cdot)g\rangle_{\mathcal{F}},
	\end{align*}
	where the fourth equality follows from $$\langle\phi,\phi'\rangle_{\mathcal{F}} = \langle\phi_0, \phi'_0\rangle_{\mathcal{F}_0}$$
	with $\phi'_0(\cdot) = \bm{K}_0(t(u), \cdot)\bm{M}(u)^* g$.
\end{proof}

Now we present a key result, which characterizes the change of kernels when applying invertible transformations on the iSVGD trajectory.
\begin{theorem}\label{keyResult}
	Let $\mathcal{H}$, $\tilde{\mathcal{H}}$, $\mathcal{X}$, $\tilde{\mathcal{X}}$, $\mathcal{Y}$, and $\tilde{\mathcal{Y}}$
	be separable Hilbert spaces satisfying
	$\mathcal{X}\subset\mathcal{Y}, \, \tilde{\mathcal{X}}\subset\tilde{\mathcal{Y}}, \,
	\mathcal{X}\subset\tilde{\mathcal{Y}}, \, \tilde{\mathcal{X}}\subset\mathcal{Y}.$
	Assume that Assumption \ref{assKernel} holds for the triples $(\mathcal{X}, \mathcal{Y}, \mathcal{H})$ and $(\tilde{\mathcal{X}},\tilde{\mathcal{Y}}, \tilde{\mathcal{H}})$
	with two fixed parameters $s\in [0,1]$, respectively.
	Let $T\in\mathcal{L}(\mathcal{Y},\tilde{\mathcal{Y}})$ and assume that $T$ is a bounded operator when
	restricted to be an operator from $\mathcal{X}$ to $\tilde{\mathcal{X}}$.
	Let $\mu$, $\mu^{\bm{d}}$ be two probability measures and $\tilde{\mu}$, $\tilde{\mu}^{\bm{d}}$ be the measures of $\tilde{u}=Tu$
	when $u$ is drawn from $\mu$, $\mu^{\bm{d}}$, respectively.
	Introduce two Stein operators $\mathcal{S}$ and $\tilde{\mathcal{S}}$ as follows:
	\begin{align*}
		& \mathcal{S}\phi(u) = \langle -DV(u), \phi(u) \rangle_{\mathcal{Y}} +
		\sum_{k=1}^{\infty}D_{k}\langle\phi(u), e_k\rangle_{\mathcal{Y}},  \quad \forall\,u\in\mathcal{X},  \\
		& \tilde{\mathcal{S}}\tilde{\phi}(\tilde{u}) = \langle -D_{\tilde{u}}V(T^{-1}\tilde{u}),
		\tilde{\phi}(\tilde{u}) \rangle_{\tilde{\mathcal{Y}}} +
		\sum_{k=1}^{\infty}D_{(\tilde{u})_k}\langle\tilde{\phi}(\tilde{u}), \tilde{e}_k\rangle_{\tilde{\mathcal{Y}}},
		\quad \forall\,\tilde{u}\in\tilde{\mathcal{X}},
	\end{align*}
	where $\{e_k\}_{k=1}^{\infty}$ and $\{\tilde{e}_k\}_{k=1}^{\infty}$ are orthonormal bases in $\mathcal{Y}$
	and $\tilde{\mathcal{Y}}$, respectively. Then, we have
	\begin{align}\label{keyTheo1}
		\mathbb{E}_{u\sim\mu}[\mathcal{S}\phi(u)] = \mathbb{E}_{u\sim\tilde{\mu}}[\tilde{\mathcal{S}}\tilde{\phi}(u)] \quad
		\text{with }\phi(u):=T^{-1}\tilde{\phi}(Tu).
	\end{align}
	Therefore, in the asymptotics of infinitesimal step size ($\epsilon\rightarrow 0^+$), it is equivalent to 
	running iSVGD with kernel $\bm{K}_0$ on $\tilde{\mu}$ and running iSVGD on $\mu$ with the kernel
	$\bm{K}(u,u') = T^{-1}\bm{K}_0(Tu,Tu')(T^{-1})^*,$
	in the sense that the trajectory of these two SVGD can be mapped to each other by the map $T$ (and its inverse).
\end{theorem}

\begin{proof}
	Let us introduce a mapping defined by $u'=f(u)=u+\epsilon\phi(u)$. 
	Denote $f_{\#}\mu$ as the probability measure $\mu\circ f^{-1}$.
	Let $\tilde{u}'\sim T_{\#}(f_{\#}\tilde{\mu})$ which is obtained by
	\begin{align}
		\begin{split}
			\tilde{u}' & = Tu' = T(u + \epsilon\phi(u)) =
			T(T^{-1}\tilde{u} + \epsilon\phi(T^{-1}\tilde{u})) \\
			& = \tilde{u} + \epsilon T\phi(T^{-1}\tilde{u}) \\
			& = \tilde{u} + \epsilon \tilde{\phi}(\tilde{u}),
		\end{split}
	\end{align}
	where we use the definition $\phi(u) = T^{-1}\tilde{\phi}(Tu)$ in (\ref{keyTheo1}).
	According to \cite[Theorem 3.1 ]{Liu2016NIPS} and \cite[Theorem 3]{Wang2019NIPS},
	we have
	$\mathbb{E}_{u^N\sim P^N_{\#}\mu}[\mathcal{S}^N\phi^N(u^N)] =
	\mathbb{E}_{u^N\sim P^N_{\#}\tilde{\mu}}[\tilde{\mathcal{S}}^N\tilde{\phi}^N(u^N)],$
	where
	\begin{align*}
		& \mathcal{S}^N\phi^N(u^N) = -\langle DV(u^N), \phi^N(u^N) \rangle_{\mathcal{Y}} +
		\sum_{k=1}^{N}D_k\langle \phi^N(u^N), e_k\rangle_{\mathcal{Y}},	\\
		& \tilde{\mathcal{S}}^N\tilde{\phi}^N(\tilde{u}^N) = -\langle D_{\tilde{u}^N}V(T^{-1}\tilde{u}^N),
		\tilde{\phi}^N(\tilde{u}^N) \rangle_{\tilde{\mathcal{Y}}} +
		\sum_{k=1}^{N}D_{(\tilde{u}^N)_k}\langle\tilde{\phi}^N(\tilde{u}^N), \tilde{e}_k\rangle_{\tilde{\mathcal{Y}}}.
	\end{align*}
	It is clear to note that there is no Jacobian matrix given by the transformation in $D_{\tilde{u}^N}V(T^{-1}\tilde{u}^N)$
	since the Jacobian matrix does not depend on $\tilde{u}^N$ for linear mappings, i.e., the derivative is zero.
	Following the proof for Theorem \ref{limitDefineApp}, we take $N\rightarrow\infty$ and obtain
	$\mathbb{E}_{u\sim \mu}[\mathcal{S}\phi(u)] =
	\mathbb{E}_{u\sim \tilde{\mu}}[\tilde{\mathcal{S}}\tilde{\phi}(u)].$
	From Theorem \ref{ChangeVariableTheorem}, when $\tilde{\phi}$ is in $\tilde{\mathcal{F}}$ with kernel $\bm{K}_0(u,u')$,
	$\phi$ is in $\mathcal{F}$ with kernel $\bm{K}(u,u')$. Therefore, maximizing $\mathbb{E}_{u\sim\mu}[\mathcal{S}\phi(u)]$
	in $\mathcal{F}$ is equivalent to $\mathbb{E}_{u\sim\tilde{\mu}}[\tilde{\mathcal{S}}\tilde{\phi}(u)]$ in $\tilde{\mathcal{F}}$.
	This suggests that the trajectory of iSVGD on $\tilde{\mu}^{\bm{d}}$ with $\bm{K}_0$ and that on $\mu^{\bm{d}}$ with $\bm{K}$ are equivalent, which completes the proof. 
\end{proof}

\begin{remark}
	Similar to the matrix-valued case \cite{Wang2019NIPS}, Theorem \ref{keyResult} suggests a conceptual procedure for constructing proper operator kernels to incorporate desirable preconditioning information. Different from the finite-dimensional case, the map $T$ is only allowed to be linear at this stage. For a nonlinear map, there is a Jacobian matrix in $\tilde{\mathcal{S}}^N\tilde{\phi}^N(\tilde{u}^N)$. It is difficult to analyze the limiting behavior of the Jacobian matrix related term. Practically, linear maps seem to be enough since even in the finite-dimensional case nonlinear maps
	will yield an unnatural algorithm \cite{Wang2019NIPS}.
\end{remark}

In the last part of this subsection, we provide some examples of preconditioning operators that are frequently used in
statistical inverse problems.

\subsubsection{Fixed preconditioning operator}

In Section 5 of \cite{Dashti2017}, the Langevin equation was considered by using 
$\mathcal{C}_0$ as a preconditioner, and an analysis was carried out for the pCN algorithm.
For the Newton based iterative method, we usually take the inverse of the second-order derivative of the objective functional
as the preconditioning operator \cite{Reyes2015Book}.
Here, we consider a linear operator $T$ that has similar properties as
$\mathcal{C}_0^{-\frac{1-s}{2}}$. Specifically, we require
\begin{align}\label{fixedpo1}
	T\in\mathcal{L}(\mathcal{H}^{1-s},\mathcal{H})\cap\mathcal{L}(\mathcal{H}^{-1-s},\mathcal{H}^{-2}).
\end{align}
Then, we specify the Hilbert space appearing in Theorem \ref{ChangeVariableTheorem} as
$\mathcal{X} = \mathcal{H}^{1-s},\, \mathcal{Y} = \mathcal{H}^{-1-s}, \,
\tilde{\mathcal{X}} = \mathcal{H}, \, \tilde{\mathcal{Y}} = \mathcal{H}^{-2}$
with $s\in [0,1]$. For the kernel $\bm{K}_0(\cdot, \cdot):\tilde{\mathcal{X}}\times\tilde{\mathcal{X}}\rightarrow\tilde{\mathcal{Y}}$, we assume that
\begin{align}\label{fixedpo2}
	\sup_{\tilde{u}\in\mathcal{H}}\|\bm{K}_0(\tilde{u},\tilde{u})\|_{\mathcal{L}(\mathcal{H}^{-2})} < \infty.
\end{align}
It follows from Theorem \ref{keyResult} that we may use a kernel of the form
\begin{align}\label{fixedpo3}
	\bm{K}(u,u') := T^{-1}\bm{K}_0(Tu,Tu')(T^{-1})^*,
\end{align}
where $u,u'\in\mathcal{H}^{1-s}$. Obviously, the kernel $\bm{K}$ given above satisfies
\begin{align}\label{fixedpo4}
	\sup_{u\in\mathcal{H}^{1-s}}\|T^{-1}\bm{K}_0(Tu, Tu)(T^{-1})^*\|_{\mathcal{L}(\mathcal{H}^{-1-s})} < \infty.
\end{align}
As an example, we may take $\bm{K}_0$ to be the scalar-valued Gaussian RBF kernel composed with operator $\mathcal{C}_0^{s}$:
\begin{align}\label{fixedpo5}
	\bm{K}_0(u,u') := \exp\Big( -\frac{1}{h}\|u-u'\|_{\mathcal{H}}^{2} \Big)\mathcal{C}_0^{s}, 
\end{align}
which yields 
\begin{align}\label{fixedpo6}
	\bm{K}(u, u') = \exp\Big( -\frac{1}{h}\|T(u-u')\|_{\mathcal{H}}^{2} \Big)T^{-1}\mathcal{C}_0^{s}(T^{-1})^*,
\end{align}
where $h$ is a bandwidth parameter. Define $\bm{K}_0^{T}(u,u'):=\bm{K}_0(Tu,Tu')$.
Let $\mathcal{P} := T^{-1}\mathcal{C}_0^{s}(T^{-1})^*$.
By simple calculations, we find that the iSVGD update direction of the kernel in (\ref{fixedpo3}) is
\begin{align}\label{fixedpo7}
	\begin{split}
		\phi^*_{\bm{K}}(\cdot) \! = &  \mathcal{P} \mathbb{E}_{u\sim\mu}\Big[
		\bm{K}_0^T(u,\cdot)(-D\Phi(u)\!-\!\mathcal{C}_0^{-1}u)
		+ \! \sum_{k=1}^{\infty}D_{k}\bm{K}_0^T(u,\cdot)e_k
		\Big] = \mathcal{P}\phi^*_{\bm{K}_0^T},
	\end{split}
\end{align}
which is a linear transform of the iSVGD update direction of the kernel $\bm{K}_0^T$ with the operator $T^{-1}\mathcal{C}_0^{s}(T^{-1})^*$.

\subsubsection{The $\mathcal{C}_0$ operator}

Choosing $T := \mathcal{C}_0^{-\frac{1-s}{2}}$, we can see that 
the condition (\ref{fixedpo1}) holds. Given the Kernel $\bm{K}_0$ in (\ref{fixedpo5}), the 
kernel $\bm{K}$ defined in (\ref{fixedpo6}) can be written as 
$$\bm{K}(u, u') = \exp\Big( -\frac{1}{h}\|\mathcal{C}_0^{-\frac{1-s}{2}}(u-u')\|_{\mathcal{H}}^{2} \Big)\mathcal{C}_0.$$
The operator $\mathcal{P}$ used in (\ref{fixedpo7}) is just $\mathcal{C}_0$. If there is only one particle,
the iSVGD update direction is then reduced to
$\phi^*_{\bm{K}}(\cdot) = \mathcal{C}_0 (D\Phi(u) + \mathcal{C}_0^{-1}u).$

\subsubsection{The Hessian operator}

For statistical inverse problems, the forward operator $\mathcal{G}$ is usually nonlinear, e.g.,
the inverse medium scattering problem \cite{Jia2019IP,Jia2018JFA}.
Around each particle $u_{i}$ with $i=1,2,\ldots,m$, the forward map can be approximated by the linearized map
\begin{align}\label{approG}
	\mathcal{G}(u) \approx \mathcal{G}(u_i) + D\mathcal{G}(u_i)(u-u_i). 
\end{align}
Assume that the potential function $\Phi$ takes the form
$\Phi(u)=\frac{1}{2}\|\Sigma^{-1/2}(\mathcal{G}(u)-d)\|_{\ell^2}^2$,
where $\Sigma$ is a positive definite matrix. Using the approximate formula (\ref{approG}), we have
\begin{align*}
	V(u)\approx \tilde{V}(u):=\frac{1}{2}\|\Sigma^{-1/2}(D\mathcal{G}(u_i)u - D\mathcal{G}(u_i)u_i + \mathcal{G}(u_i)-d)\|_{\ell^2}^2 +
	\frac{1}{2}\|\mathcal{C}_0^{-1/2}u\|_{\mathcal{H}}^2.
\end{align*}
It follows from a simple calculation that 
$D^2\tilde{V}(u_i) = D\mathcal{G}(u_i)^*\Sigma^{-1}D\mathcal{G}(u_i) + \mathcal{C}_0^{-1}.$
For the Newton-type iterative method, we can take the linear transformation
$T=\mathcal{C}_0^{s/2}(\frac{1}{m}\sum_{i=1}^{m}(D\mathcal{G}(u_i)^*\Sigma^{-1}D\mathcal{G}(u_i) + \mathcal{C}_0^{-1}))^{1/2}$.
If $\mathcal{G}$ is a linear operator (e.g., the examples in \cite{Jia2020IPI}), it is easy to verify the condition (\ref{fixedpo1}).
For nonlinear problems, it is necessary to employ the regularity properties of the direct problems, which
is beyond the scope of this work. Hence we will not verify this condition in this paper and leave it as a future work.
With this choice of $T$, the kernel (\ref{fixedpo6}) and the iSVGD update direction (\ref{fixedpo7}) can be easily obtained.
If there is only one particle, the iSVGD update direction is degenerated to the usual Newton update direction when evaluating MAP estimate.

\subsubsection{Mixture preconditioning}\label{mixturePreSubSec}

Using a fixed preconditioning operator, we can not specify different preconditioning operators for different particles.
Inspired by the mixture precondition \cite{Wang2019NIPS},
we propose an approach to achieve point-wise preconditioning. The idea is to use a weighted combination of
several linear preconditioning operators. This involves leveraging a set of anchor points
$\{v_{\ell}\}_{\ell=1}^{m}$, each of which is associated with a preconditioning operator $T_{\ell}$
(e.g., $T_\ell = \mathcal{C}_0^{s/2}(D\mathcal{G}(v_\ell)^*\Sigma^{-1}D\mathcal{G}(v_\ell) + \mathcal{C}_0^{-1})^{1/2}$).
In practice, the anchor points $\{v_{\ell}\}_{\ell=1}^{m}$ can be set to be the same as
the particles $\{u_i\}_{i=1}^{m}$. We then construct a kernel by
$\bm{K}(u,u')=\sum_{\ell=1}^{m}\bm{K}_{\ell}(u,u')w_{\ell}(u)w_{\ell}(u'),$
where
\begin{align}\label{newKerMix1}
	\bm{K}_{\ell}(u,u') := T_{\ell}^{-1}\bm{K}_0(T_{\ell}u,T_{\ell}u')(T_{\ell}^{-1})^*,
\end{align}
and $w_{\ell}(u)$ is a positive scalar-valued function that determines the contribution of kernel $\bm{K}_{\ell}$
at point $u$. Here $w_{\ell}(u)$ should be viewed as a mixture probability, and hence should satisfy
$\sum_{\ell=1}^{m}w_{\ell}(u) = 1$ for all $u$. In our empirical studies, we take
\begin{align}
	w_{\ell}(u) = \frac{\exp\Big( -\frac{1}{2}\|T_{\ell}(u-v_{\ell})\|_{\mathcal{H}}^2 \Big)}{
		\sum_{\ell'=1}^{m}\exp\Big( -\frac{1}{2}\|T_{\ell'}(u-v_{\ell'})\|_{\mathcal{H}}^2 \Big)}.
\end{align}
In this way, each point $u$ is mostly influenced by the anchor point closest to it, which allows 
to apply different preconditioning for different points. In addition, the iSVGD update direction has the
form
\begin{align}\label{formulaMix}
	\begin{split}
		\phi_{\bm{K}}^*(\cdot) = \sum_{\ell=1}^m w_{\ell}(\cdot)\mathbb{E}_{u\sim\mu}\Big[ &
		-w_{\ell}(u)\bm{K}_{\ell}(u,\cdot)(D\Phi(u)+\mathcal{C}_0^{-1}u) \\
		& + \sum_{k=1}^{\infty}D_{k}(w_{\ell}(u)\bm{K}_{\ell}(u,\cdot)e_{k})
		\Big],
	\end{split}
\end{align}
which is a weighted sum of a number of iSVGD update directions with linear preconditioning operators. The implementation details of (\ref{formulaMix}) are given in the supplementary material.

\begin{remark}\label{remark_s}
	For the kernel defined above, the particles should belong to the Hilbert space $\mathcal{H}^{1-s}$.
	Based on the studies 
	the finite-dimensional problems \cite{Wang2019NIPS}, we may let the parameter $s$ be equal to $0$. However, when the parameter $s=0$,
	each particle $u_i$ belongs to $\mathcal{H}^{1}$ which is the Cameron--Martin space of the prior measure.
	By the classical Gaussian measure theory \cite{DaPrato1996Book}, we know that $\mathcal{H}^{1}$ has zero measure.
	This fact implies that all of the particles belong to a set with zero measure, which may lead to too concentrated
	particles and deviates from our purpose. Hence we should choose $s > 0$
	to ensure the effectiveness of the SVGD sampling algorithm. These observations are illustrated by our numerical experiments 
	in Section \ref{ApplicationSection}.
\end{remark}

\subsection{Some insights about iSVGD}\label{SubsectionInsight}

We have constructed the well-defined iSVGD algorithms with or without preconditioning operators,
which is the first step to extend the finite-dimensional SVGD to the infinite-dimensional space.
Some mathematical studies have been carried out for the finite-dimensional SVGD, e.g., gradient flow on probability space \cite{Liu2017NIPS} and mean field limit theory related to the macroscopic behavior \cite{Lu2019SIAM}.
These results provide in-depth understandings of the SVGD algorithm and motivate many new algorithms \cite{pmlr-v97-liu19i}.
In this subsection, we intend to provide a preliminary mathematical study on the iSVGD under a simpler setting. 

We consider the kernel operator $\bm{K}(u,v):=K(\|u-v\|_{\mathcal{H}})\text{Id}$ with $u,v\in\mathcal{H}$ and $K(\cdot)$ being a scalar function. 
Let $m$ be the sample number and $V(u)$ be defined in (\ref{defV}). Similar to the finite-dimensional case,
the iterative procedure in Algorithm \ref{infiniteSVGD} can be viewed as a particle system:
\begin{align}\label{particlesystem}
	\begin{split}
		& \frac{d}{dt}u_i(t) = -(\tilde{D}\bm{K}*\mu_{m}(t))(u_i(t)) - (\bm{K}*DV \mu_{m}(t))(u_i(t)), \\
		& \mu_{m}(t) = \frac{1}{m}\sum_{j=1}^{m}\delta_{u_j(t)}, \\
		& u_i(0) = u_i^0,  \quad i=1,2,\ldots,m,
	\end{split}
\end{align}
where $\{u_i^0\}_{i=1}^{m}$ are the initial particles, $\delta_{u_i(t)}$ denotes the Dirac measure concentrated on $u_i(t)$ with $i=1,2,\ldots,m$,
``$*$'' denotes the usual convolution operator, and $\tilde{D}\bm{K}(u-v) = \sum_{k=1}^{\infty}D_{u_k}\bm{K}(u-v)e_k$.
For convenience, we write the two convolution terms in the following forms:
\begin{align*}
	(\tilde{D}\bm{K}*\mu_{m}(t))(u_i(t)) & = \frac{1}{m}\sum_{j=1}^{m}\tilde{D}\bm{K}(u_i(t)-u_j(t)), \\
	(\bm{K}*DV \mu_{m}(t))(u_i(t)) & = \frac{1}{m}\sum_{j=1}^{m}\bm{K}(u_i(t)-u_j(t))DV(u_j(t)).
\end{align*}
Similarly, we consider the weak form equation about the measure-valued function:
\begin{align}\label{PDEonHilbert}
	\begin{split}
		& \frac{d}{dt}\langle\mu(t),\varphi\rangle = \langle\mu(t), L(\mu(t))\varphi\rangle, \\
		& \mu(0) = \nu,
	\end{split}
\end{align}
where $\nu$ is the probability measure employed to generate initial particles, $\varphi$ is the test function, and
\begin{align}
	L(\mu(t))\varphi = \langle \tilde{D}\bm{K}*\mu(t), D\varphi\rangle_{\mathcal{H}} + \langle\bm{K}*DV\mu(t), D\varphi\rangle_{\mathcal{H}}.
\end{align}
Let $W^{1,2}(\mathcal{H},\mu)$ be the usual Sobolev space defined for a Gaussian measure $\mu$ \cite{Prato2004Book}.

\begin{theorem}\label{invarianceSolution}
	Let $\mu_0$ and $\Phi$ be the prior measure and potential function defined in (\ref{BasicBayesFormula}), respectively.
	Assume $K(\cdot)\in W^{1,2}(\mathcal{H},\mu_0)$ and $e^{-\Phi(\cdot;\bm{d})}\in L^2(\mathcal{H},\mu_0)$.
	Then, the posterior measure $\mu^{\bm{d}}$ defined in (\ref{BasicBayesFormula}) is an invariant solution to Eq. (\ref{PDEonHilbert}),
	i.e., when $\nu:=\mu^{\bm{d}}$, the solution $\mu(t)$ of (\ref{PDEonHilbert}) is equal to $\mu^{\bm{d}}$.
\end{theorem}

The proof is given in the supplementary material. Clearly, this theorem holds in the finite-dimensional setting.
We point out that the integration by parts may not hold for the infinite-dimensional case. 
In the finite-dimensional setting, the analysis of the corresponding particle system (\ref{particlesystem}) and Eq. (\ref{PDEonHilbert}) have been given recently in \cite{Lu2019SIAM}. It is sophisticated to define meaningful solutions for the above interacting particle system (\ref{particlesystem}) and the measure-valued function equation (\ref{PDEonHilbert}), which  are beyond the scope of this study and are left for future work.
One of the major difficulties for the infinite-dimensional case is that $\mathcal{C}_0^{-1}$
(the precision operator of the prior measure) is usually an unbounded operator \cite{Dashti2017}. Nearly all of the estimates presented in \cite{Lu2019SIAM}
for the finite-dimensional case cannot be adopted for the infinite-dimensional setting.

Numerical experiments indicate that the SVGD without preconditioning operators
can hardly provide accurate estimates for some inverse problems. The SVGD with preconditioning operators can accelerate the convergence and give reliable estimates efficiently. In addition, the unboundedness issue induced by the precision operator
$\mathcal{C}_0^{-1}$ may be overcome by introducing preconditioning operators.
A detailed analysis of the iSVGD with preconditioning operators may be a good starting point for future theoretical studies.

At the end of this subsection, we mention a critical difference between finite- and infinite-dimensional theories. 
It follows from Theorem 2.7 in \cite{Lu2019SIAM} and Theorem 1.1 in \cite{Trillos2015CJM} that the empirical 
measure constructed by particles in finite-dimensional SVGD can approximate the continuous counterpart with accuracy $\epsilon$ when the number of particles are of order $O(\epsilon^{d})$, where $d$ is the discrete dimension. 
Obviously, an infinite number of particles is needed if the dimension $d$ goes to infinity, which indicates that 
the infinite-dimensional theory may be meaningless. 

The above statement explains that not every finite-dimensional setting can be meaningfully generalized to the 
infinite-dimensional space. The assumption on prior measure is important for the infinite-dimensional theory (the current assumption may be slightly relaxed, e.g., the Besov type measure). According to the general analysis for the convergence and concentration of empirical measures given in \cite{Lei2020Bernoulli}, we believe that the prior measures used here can be approximated by the empirical measures under the Wasserstein distance on infinite-dimensional Hilbert space. Specifically, the estimate of the convergence speed is not relevent to the dimension when considering some finite-dimensional spaces as the projected infinite-dimensional space. If a theorem similar as Theorem 2.7 in \cite{Lu2019SIAM} for
the system (\ref{particlesystem})--(\ref{PDEonHilbert}) can be proved, we are able to confirm that the particles obtained by iSVGD can approximate
the posterior measure for certain accuracy with particle numbers independent of the discrete dimension. 
However, it is higly non-trivial to carry out an in-depth study of the system (\ref{particlesystem})--(\ref{PDEonHilbert})  and is beyond the scope of the current work. In Subsection 6.3 of the supplementary material, we give a numerical illustration to address this issue. 

\section{Applications}\label{ApplicationSection}

The proposed framework is valid for Bayesian inverse problems governed by any systems of PDEs.
Due to the page limitation, we present one example of an inverse problem governed
by the steady state Darcy flow equation. The second example concerns an inverse problem of the Helmholtz equation and is given in the supplementary material. 

Consider the following PDE model:
\begin{align}\label{DarcyEq}
	\begin{split}
		-\nabla\cdot(e^u\nabla w) & = f \quad\text{in }\Omega, \\
		w & = 0 \quad\text{on }\partial\Omega,
	\end{split}
\end{align}
where $\Omega\subset\mathbb{R}^2$ is a bounded Lipschitz domain,
$f(x)$ denotes the sources, and $e^{u(x)}$ describes the permeability of the porous medium.
This model is used as a benchmark problem in many works, e.g., the preconditioned Crank--Nicolson (pCN) algorithm \cite{Cotter2013SS} and the sequential Monte Carlo method \cite{Beskos2015SC}. We will compare the performance of the proposed iSVGD approach
with the pCN \cite{Cotter2013SS} and the randomized maximum a posterior (rMAP) method \cite{Wang2018SISC}.

\subsection{Basic settings and finite-element discretization}

For numerical implementations, it is essential to compute all of the related gradients and Hessian operators before discretization (i.e., pushing the discretization to the last step).
A direct calculation yields the gradient and Hessian operators of the operator-valued kernel, but the adjoint method \cite{Reyes2015Book} needs to be employed for the potential $\Phi$ involving PDEs.
More discussions on finite- and infinite-dimensional approaches can be found in the supplementary material, 
which might be helpful for readers who are not familiar with infinite-dimensional approach.
Let $\mathcal{F}$ be the solution operator that maps the parameter $u$ to the solution of (\ref{DarcyEq}),
and $\mathcal{M}$ be the measurement operator defined as
$\bm{d} = \mathcal{M}(w) = (\ell_{x_1}(w), \ell_{x_2}(w),\ldots, \ell_{x_{N_d}}(w))^T,$
where
\begin{align}
	\ell_{x_j}(w) = \int_{\Omega}\frac{1}{2\pi\delta^2}e^{-\frac{1}{2\delta^2}\|x-x_j\|^2} w(x)dx
\end{align}
with $\delta>0$ being a sufficiently small number and $x_i\in\Omega$ for $i=1,\ldots,N_d$.
The forward map can be defined as
$\mathcal{G} := \mathcal{M}\circ\mathcal{F}$, and the problem can be written in the abstract form
$\bm{d} = \mathcal{G}(u) + \bm{\epsilon}$ with $\bm{\epsilon}\sim\mathcal{N}(0,\sigma^2\text{Id})$.
Then we have
$\Phi(u) = \frac{1}{2\sigma^2}\|\mathcal{M}(w) - \bm{d}\|^2$.
The gradient $D\Phi(u)$ acting in any direction $\tilde{u}$ is given by
\begin{align}\label{DradientForm}
	\langle D\Phi(u), \tilde{u}\rangle = \int_{\Omega}\tilde{u}e^u\nabla w\cdot \nabla p dx,
\end{align}
where the adjoint state $p$ satisfies \emph{the adjoint equation}
\begin{align}\label{DarcyAdjoint}
	\begin{split}
		-\nabla\cdot(e^u\nabla p) & = -\frac{1}{\sigma^2}\sum_{j=1}^{N_d}\frac{1}{2\pi\delta^2}e^{-\frac{1}{2\delta^2}\|x-x_j\|^2}(\ell_{x_j}(w)-d_j)\quad\text{in }\Omega, \\
		p & = 0\quad\text{on }\partial\Omega.
	\end{split}
\end{align}
The Hessian acting in direction $\tilde{u}$ and $\hat{u}$ reads
\begin{align}\label{HessianForm}
	\begin{split}
		\langle\langle D^2\Phi(u), \hat{u}\rangle, \tilde{u}\rangle = & \int_{\Omega} \hat{u}\tilde{u}e^u\nabla w\cdot\nabla p dx
		+ \int_{\Omega}\tilde{u}e^u\nabla w\cdot\nabla\hat{p} dx \\
		& + \int_{\Omega}\tilde{u}e^u\nabla p\cdot\nabla\hat{w}dx,
	\end{split}
\end{align}
where the state $\hat{w}$ satisfies \emph{the incremental forward equation}
\begin{align}
	\begin{split}
		-\nabla\cdot(e^u\nabla\hat{w}) & = \nabla\cdot(\hat{u}e^u\nabla w) \quad \text{in }\Omega, \\
		\hat{w} & = 0 \quad\text{on }\partial\Omega,
	\end{split}
\end{align}
and the state $\hat{p}$ satisfies \emph{the incremental adjoint equation}
\begin{align}
	\begin{split}
		-\nabla\cdot(e^u\nabla\hat{p}) & = \nabla\cdot(\hat{u}e^u\nabla p) - \frac{1}{2\pi\delta^2\sigma^2}\sum_{j=1}^{N_d}\hat{w}e^{-\frac{1}{2\delta^2}\|x-x_j\|^2} \quad\text{in }\Omega, \\
		\hat{p} & = 0 \quad\text{on }\partial\Omega.
	\end{split}
\end{align}

In experiments, we choose $\Omega$ to be a rectangular domain $\Omega = [0,1]^2 \subset\mathbb{R}^2$,
set $\mathcal{H}=L^{2}(\Omega)$, and
consider the prior measure $\mu_0 = \mathcal{N}(u_0, \mathcal{C}_0)$ with the mean function $u_0$ and
the covariance operator $\mathcal{C}_0:=A^{-2}$, where
$A = \alpha (I-\Delta) \, (\alpha > 0)$
with the domain of $A$ given by
$D(A):=\Big\{ u\in H^2(\Omega) \,:\, \frac{\partial u}{\partial\bm{n}}=0\text{ on }\partial\Omega \Big\}.$
Here, $H^2(\Omega)$ is the usual Sobolev space.
Assume that the mean function $u_0$ resides in the Cameron--Martin space of $\mu_0$.

Based on (\ref{DradientForm}) and (\ref{HessianForm}), we can prove the following results, which satisfy Assumptions \ref{assumPhi}. The proof is given in the supplementary material.

\begin{theorem}\label{ConditionVerifyThm}
	Let $H^{-1}(\Omega)$ be the usual Sobolev space with the regularity index $-1$. Assume $\mathcal{X}=\mathcal{H}^{1-s}$ with the parameter $s<0.5$, and then we have
	\begin{align*}
		0 \leq \Phi(u) &\leq C(1+\|f\|_{H^{-1}})^2 e^{2\|u\|_{\mathcal{X}}}, \\
		\|D\Phi(u)\|_{\mathcal{X}^*} &\leq C (1+\|f\|_{H^{-1}})^2 e^{4\|u\|_{\mathcal{X}}},  \\
		\|D^2\Phi(u)\|_{\mathcal{L}(\mathcal{X},\mathcal{X}^*)} &\leq C(1+\|f\|_{H^{-1}})^2 e^{6\|u\|_{\mathcal{X}}}.
	\end{align*}
\end{theorem}

In the following, we use the Gaussian kernel, i.e., $\bm{K}(u,u') = \exp\left(-\frac{1}{h}\|u-u'\|_{\mathcal{H}}^2\right)$,
for the iSVGD without preconditioning operators. For numerical examples with preconditioning operators, we employed the kernel
given in Subsection \ref{mixturePreSubSec}.

For finite-dimensional approximations, we consider a finite-dimensional subspace $V_h$ of $L^2(\Omega)$ originating
from the finite element discretization with the continuous Lagrange basis functions $\{ \phi_j \}_{j=1}^{n}$, which correspond to the nodal
points $\{x_j\}_{j=1}^n$, such that $\phi_j(x_i)=\delta_{ij}$ for $i,j\in\{1,\ldots,n\}$.
Instead of statistically inferring parameter functions $u\in L^2(\Omega)$, we consider the
approximation $u_h=\sum_{j=1}^n u_j\phi_j \in V_h$. Under this finite-dimensional approximation,
we can employ the numerical method provided in \cite{Tan2013SISC} to discretize the prior, and construct
finite-dimensional approximations of the Gaussian approximation of the posterior measure.
Based on our analysis in Subsection \ref{SVGDPre}, we need to calculate the fractional powers of the operator $\mathcal{C}_0$.
Here, we employ the matrix transfer technique (MTT) \cite{Thanh2016IPI}. The main idea of MTT is to indirectly discretize a fractional Laplacian using a discretization of the standard Laplacian.
As discussed in \cite{Tan2013SISC}, the operator $M$ is taken as 
\begin{align}
	M=(M_{ij})_{i,j=1}^n \quad\text{and}\quad M_{ij} = \int_{\Omega}\phi_i(x)\phi_j(x)dx, \quad i,j\in\{1,\ldots,n\}.
\end{align}
The matrix $M^{1/2}$ is approximated by
the diagonal matrix $\text{diag}(M_{11}^{1/2},\ldots,M_{nn}^{1/2})$.

Finally, we mention that the finite element discretization is implemented by employing the open
software FEniCS (Version 2019.1.0) \cite{Logg2012Book}.
All programs were run on a personal computer with Intel(R) Core(TM) i7-7700 at 3.60 GHz (CPU), 32 GB (memory),
and Ubuntu 18.04.2 LTS (OS).

\subsection{Numerical results}

In the experiments, the noise level is fixed to be $1\%$ since the goal is to test algorithms rather than demonstrating the Bayesian modeling.  We compare the iSVGD with the mixture preconditioning operator (iSVGDMPO) with
the preconditioned Crank--Nicolson (pCN) sampling algorithm \cite{Dashti2017} and the 
randomized maximum a posteriori (rMAP) algorithm \cite{Wang2018SISC}. Since the rMAP sampling algorithm is not accurate for nonlinear problems, we choose $\alpha=0.5$ in the prior probability measure.
It should be mentioned that we choose the anchor points in the iSVGDMPO just to be the same as the particles 
and the anchor points will be updated during the iterations. 
The initial particles of the iSVGD are generated from a probability measure by using the method proposed in \cite{Tan2013SISC}. 

For the current settings, the gradient descent based method seems hardly to
find appropriate solutions in reasonable iterative steps. Hence, the optimization method with preconditioning operators, e.g., the Newton-conjugate gradient method, is employed.
The term $\mathbb{E}_{u'\sim\mu_{\ell}}[\bm{K}(u',u)DV(u')]$ in (\ref{algformul1}) is
an averaged gradient descent component in the whole iterative term, which drives all of the particles to be concentrated.
We anticipate that Algorithm \ref{infiniteSVGD} cannot work well due to the inefficiency of the gradient descent algorithm.
Due to the page limitation, numerical results are given in the supplementary material,
which show that Algorithm \ref{infiniteSVGD} does not perform well in some cases.
This is one of the main motivations for us to study the iSVGD with preconditioning operators.

We compare the iSVGD with the mixture preconditioning operator (iSVGDMPO) with
those obtained by the pCN and rMAP sampling algorithms. 
As illustrated in Remark \ref{remark_s}, the parameter $s$ should not be zero. Intuitively, the particles should belong to a space with probability approximately equal to one under the prior measure $\mu_0$. By the Gaussian measure theory \cite{DaPrato1996Book}, we may take $s > 0.5$ since $\mu_0(\mathcal{H}^{1-s}) = 1$
for any $s>0.5$. Since the posterior measure is usually concentrated on
a small support set of the prior measure, the parameter $s$ should be slightly smaller than $0.5$.
Thus, we set $s=0.3$ or $0.4$ in our examples.
Usually, the initial particles are scattered, and the variances of the initial particles are larger than the
final particles obtained by the iSVGDMPO. We design the following adaptive empirical strategy for $s$:
\begin{align}\label{adaptive_s}
	s = -0.5\frac{\|\text{var}\|_{\ell^2}}{\|\text{var}_{0}\|_{\ell^2}} + 0.5,
\end{align}
where $\text{var}$ is the current estimated variance, $\text{var}_{0}$ is the estimated variance of the initial particles,
and $\|\cdot\|_{\ell^2}$ is the usual $\ell^2$-norm. Obviously, for the initial particles, we have $s=0$.
The particles are forced to be concentrated. When the variance is reduced, the parameter $s$ approaches $0.5$ to avoid
that the particles are concentrated on a set with zero measure.
Since the pCN is a dimension independent MCMC type sampling algorithm, we take the results obtained
by the pCN as the baseline (accurate estimate). To make sure that the pCN algorithm yields an accurate estimate, we
iterate $10^6$ steps and withdraw the first $10^5$ samples.
Several different step-sizes are tried and the traces of some parameters are plotted,
and then the most reliable one is picked as the baseline.

\begin{figure}[t]
	\centering
	\includegraphics[width=0.95\textwidth]{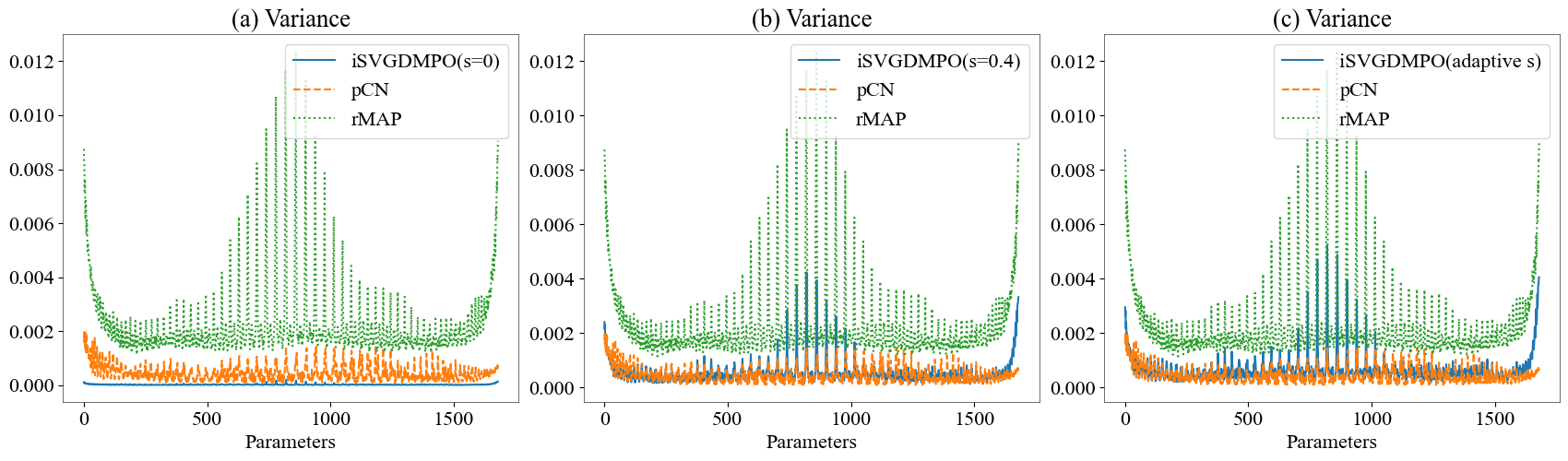}\\
	\vspace*{0.0em}
	\caption{The comparison of the variances estimated by the pCN, rMAP, iSVGDMPO with different $s$.
		(a): $s=0$;
		(b): $s=0.4$;
		(c): adpatively chosen $s$.
	}\label{Fig_compare_s}
\end{figure}

In Figure \ref{Fig_compare_s}, we show the estimated variances obtained by the iSVGDMPO (blue solid line),
rMAP (green dotted line), and the pCN (orange dashed line) sampling algorithms.
The estimated variances of the iSVGDMPO are shown for $s=0$ and $s=0.4$ on the left and in the middle, respectively. On the right, we exhibit the estimated variances when the empirical adaptive strategy (\ref{adaptive_s}) is employed.
As expected, the estimated variances are too small when $s=0$, which indicates that the particles are concentrated on
a small set. Choosing $s=0.4$ or using the empirical strategy, we obtain similar estimates, which is more similar to the baseline obtained by the pCN compared with the estimates obtained by the rMAP.

One important question arises: how does $s$ influence the convergence of the iSVGDMPO? The detailed numerical comparisons are given in the supplementary material. Here we state the conclusions:
The convergence speeds are similar for $s=0.4$ and the adaptively chosen $s$. When specifying $s=0.5$, the variances will gradually approach the background truth,
but the convergence speed seems much slower than $s=0.4$ or the adaptively chosen $s$.
In the following numerical experiments, we use the empirical adaptive strategy to specify the parameter $s$.

In addition, we provide three videos to exhibit the dynamic changing procedure of the
estimated variances in the supplementary material. The update perturbation with and without repulsive force term are exhibited.
These videos can further illustrate our theoretical findings.
We can see that the repulsive force terms indeed prevent the particles from being over concentrated.

\begin{figure}[t]
	\centering
	\includegraphics[width=0.8\textwidth]{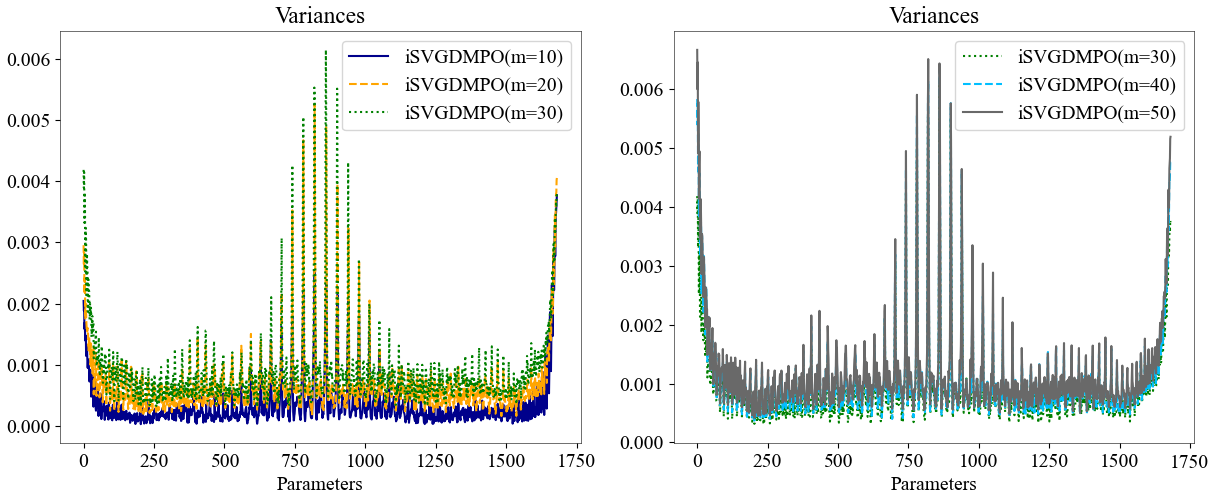}\\
	\vspace*{0.0em}
	\caption{{\small The comparison of the variances estimated by the iSVGDMPO with $s=10,20,30,40,50$.}
	}\label{Fig_compare_sample_number}
\end{figure}

Apart from the parameter $s$, how many samples should be taken to guarantee a stable statistical quantity estimate is important for using the iSVGDMPO.
When the particle number is too small, we cannot obtain reliable estimates. However, the computational complexity
increases when the particle number increases. In Figure \ref{Fig_compare_sample_number},
we show the estimated variances when particle number equals to $10$, $20$, $30$, $40$, and $50$.
Denote by $m$ the number of samples.
On the left in Figure \ref{Fig_compare_sample_number}, we show the results obtained when $m=10,20,30$.
Obviously, when $m=10$, the estimated variances are significantly smaller than those obtained when $m=20,30$.
On the right in Figure \ref{Fig_compare_sample_number}, we find that the estimated variances are similar
when $m=30,40,50$. Hence, it is enough for our numerical examples to take $m=20$ or $30$, which attains a balance
between efficiency and accuracy. So far, we have only compared the variances with different parameters in the iSVGDMPO.
In the following, qualitative and quantitative comparisons of other statistical quantities are provided to
illustrate the effectiveness of the iSVGDMPO.

\begin{figure}[htbp]
	\centering
	\includegraphics[width=1.0\textwidth]{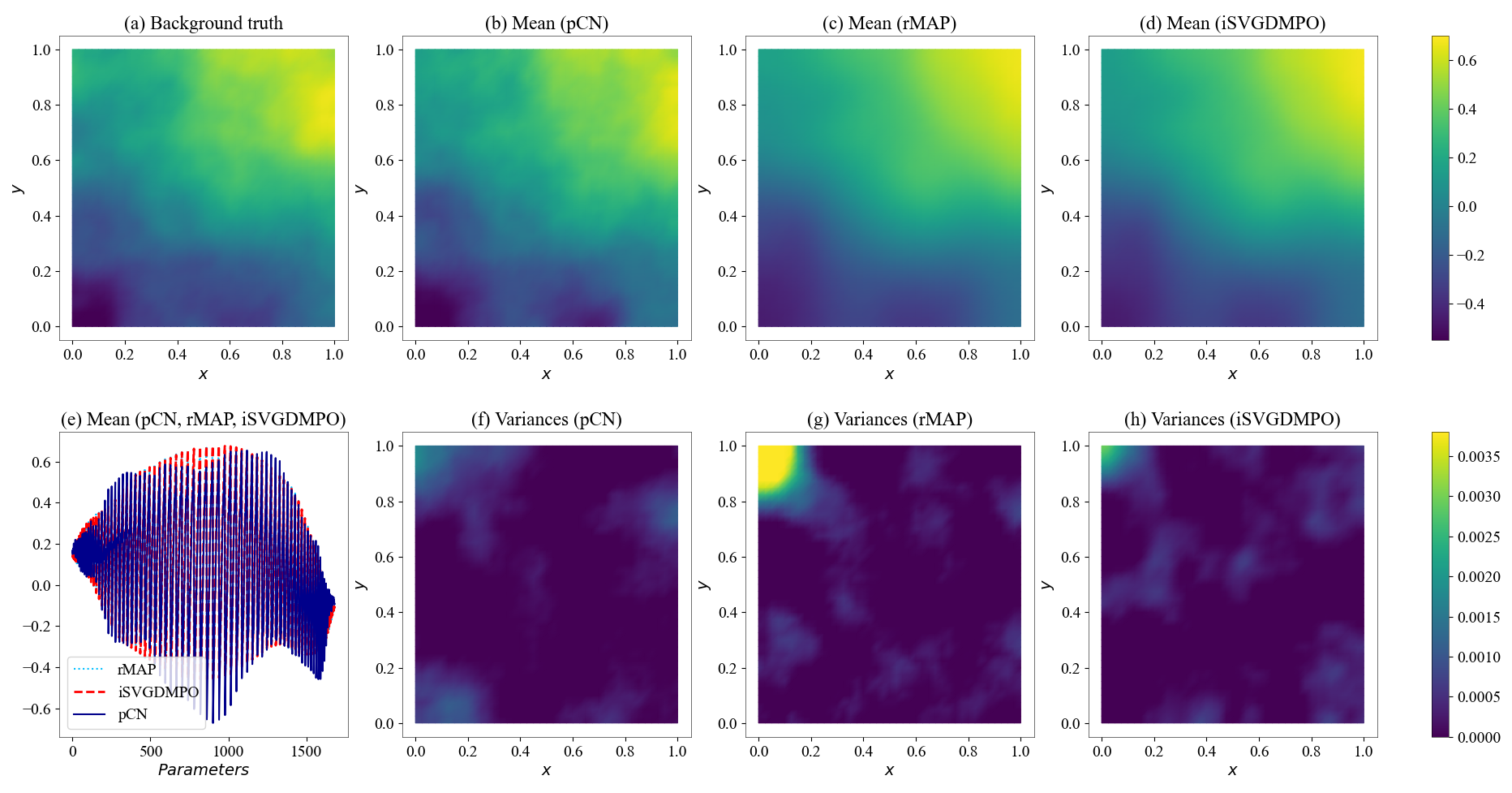}\\
	\vspace*{0.0em}
	\caption{The background truth and the estimated mean and variance functions by the pCN, rMAP, and iSVGDMPO.
		(a): The background truth; (b): The estimated mean function by the pCN;
		(c): The estimated mean function by the rMAP; (d): The estimated mean function by the iSVGDMPO;
		(e): The estimated mean function on mesh points by the pCN (blue solid line), rMAP (light blue dotted line), and iSVGDMPO (red dashed line);
		(f): The estimated variances by the pCN; (g): The estimated variances by the rMAP;
		(h): The estimated variances by the iSVGDMPO.
	}\label{Fig_truth_mean_variances}
\end{figure}

\begin{figure}[htbp]
	\centering
	\includegraphics[width=0.95\textwidth]{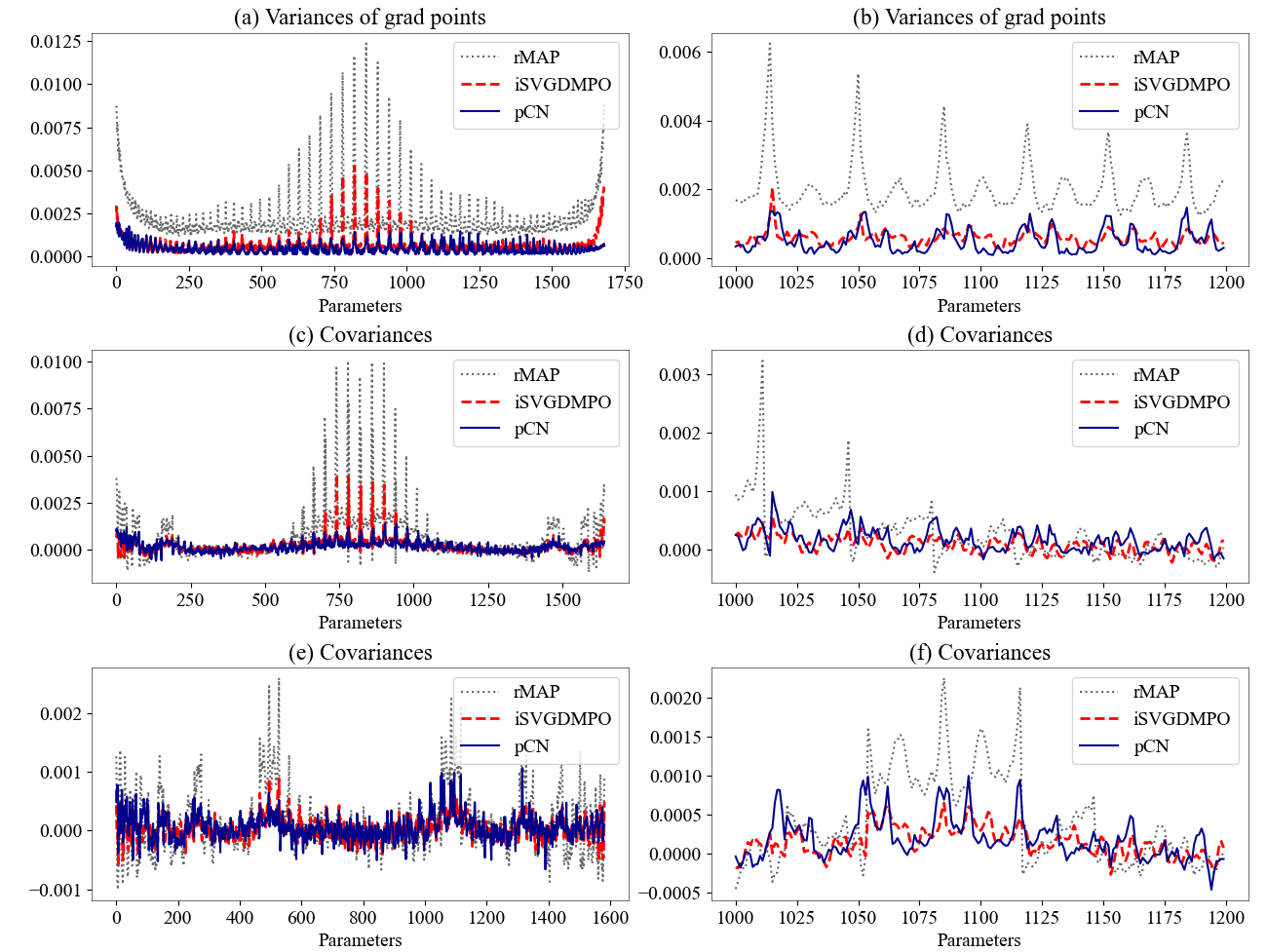}\\
	\vspace*{0.0em}
	\caption{The estimated variances and covariances by the pCN (blue solid line), rMAP (gray dotted line), and iSVGDMPO (red dashed line).
		(a): The estimated variances $\{\text{var}_u(x_i)\}_{i=1}^{N_g}$ on all mesh points;
		(b): The estimated variances for mesh points with indexes from $1000$ to $1200$ (show details);
		(c): The estimated covariances $\{\text{cov}_u(x_i, x_{i+50})\}_{i=1}^{N_g-50}$ on mesh point pairs $\{(x_i,x_{i+50})\}_{i=1}^{N_g-50}$;
		(d): The estimated covariances shown in (c) with indexes from $1000$ to $1200$ (show details);
		(e): The estimated covariances $\{\text{cov}_u(x_i, x_{i+100})\}_{i=1}^{N_g-100}$ on mesh point pairs $\{(x_i,x_{i+100})\}_{i=1}^{N_g-100}$;
		(f): The estimated covariances shown in (e) with indexes from $1000$ to $1200$ (show details).
	}\label{Fig_variance_covariance}
\end{figure}

Now, we specify the sampling number $m=30$ and set the parameter $s$ by the proposed empirical strategy (\ref{adaptive_s}).
In Figure \ref{Fig_truth_mean_variances}, we show the background truth and the estimated mean and variance functions
obtained by the pCN, rMAP, and iSVGDMPO, respectively. The iterative number of the iSVGDMPO is set to be $30$.
From the first line, we observe that the mean functions obtained by the rMAP and iSVGDMPO are similar,
which are slightly smoother than the one obtained by the pCN algorithm. This may be caused by the inexact matrix-free Newton-conjugate
gradient algorithm \cite{Tan2013SISC}. As investigated in \cite{Wang2018SISC},
many more powerful Newton-type algorithms can be employed to improve the performance both of the rMAP and iSVGDMPO.
For the variances, the iSVGDMPO gives more reliable estimates compared with the rMAP, as can be seen from Figure \ref{Fig_truth_mean_variances} (f), (g), and (h).

Next, we provide some more comparisons of statistical quantities between the results obtained by
the pCN, rMAP, and iSVGDMPO. The samples are discretization of functions.
As introduced in \cite{Ramsay2005Book}, the mean, variance and covariance functions are the main statistics for functional data.
The variance function denoted by $\text{var}_u(x)$ can be defined as $\text{var}_u(x) = \frac{1}{m}\sum_{i=1}^{m}(u_i(x)-\bar{u}(x))^2$,
where $x\in\Omega$ is a point residing in the domain $\Omega$, $\bar{u}$ is the mean function, and $m$ is the sample number.
The covariance function can be defined as $\text{cov}_u(x_1, x_2) = \frac{1}{m-1}\sum_{i=1}^{m}(u_i(x_1)-\bar{u}(x_1))(u_i(x_2)-\bar{u}(x_2))$,
where $x_1, x_2\in\Omega$ and $m, \bar{u}$ are defined as in $\text{var}_u(x)$.
For simplicity, we compute these quantities on the mesh points and exhibit the results in Figure \ref{Fig_variance_covariance}.
In all of the subfigures in Figure \ref{Fig_variance_covariance}, the 
estimates obtained by the pCN, rMAP, and iSVGDMPO are drawn in blue solid line, gray dotted line, and red dashed line, respectively.
In Figure \ref{Fig_variance_covariance} (a), we show the variance function calculated on all of the mesh points, i.e.,
$\{\text{var}_u(x_i)\}_{i=1}^{N_g}$ ($N_g$ is the number of mesh points).
In Figure \ref{Fig_variance_covariance} (c) and (e), we show the covariance function calculated on the pairs of points
$\{(x_i,x_{i+50})\}_{i=1}^{N_g-50}$ and $\{(x_i,x_{i+100})\}_{i=1}^{N_g-100}$, respectively.
Compared with the estimates given by the rMAP, we can find that the estimates obtained by the iSVGDMPO are visually more similar to the estimates provided by the pCN. In Figure \ref{Fig_variance_covariance} (b), (d), and (f), we provide the same estimates
shown in (a), (c), and (e) with points indexing from $1000$ to $1200$, which give more detailed comparisons.
The results also confirm that the iSVGDMPO provides more similar estimates to the pCN.

\begin{table}[htbp]\label{tab1}
	\caption{The $\ell^2$-norm error of the variance and covariance functions on mesh points for the rMAP and iSVGDMPO (the estimates of the pCN are seen as the background truth).}
	{\small
		\begin{tabular}{ccccc}
			\Xhline{1.1pt}
			& $\text{var}_u(x_i)$ & $\text{cov}_u(x_i,x_{i+10})$ & $\text{cov}_u(x_i,x_{i+20})$ & $\text{cov}_u(x_i,x_{i+30})$  \\
			\hline
			{\small rMAP} & $0.00759$ & $0.00100$ & $0.00075$ & $0.00092$  \\
			\hline
			{\small iSVGDMPO} & $0.00038$ & $0.00012$ & $0.00009$ & $0.00010$  \\
			\Xhline{1.1pt}
			& $\text{cov}_u(x_i,x_{i+40})$ & $\text{cov}_u(x_i,x_{i+50})$ & $\text{cov}_u(x_i,x_{i+60})$ & $\text{cov}_u(x_i,x_{i+70})$ \\
			\hline
			{\small rMAP} & $0.00227$ & $0.00038$ & $0.00043$ & $0.00056$  \\
			\hline
			{\small iSVGDMPO} & $0.00015$ & $0.00007$ & $0.00006$ & $0.00007$  \\
			\Xhline{1pt}	
			& $\text{cov}_u(x_i,x_{i+80})$ & $\text{cov}_u(x_i,x_{i+90})$ & $\text{cov}_u(x_i,x_{i+100})$ & $\text{cov}_u(x_i,x_{i+110})$ \\
			\hline
			{\small rMAP} & $0.00142$ & $0.00029$ & $0.00031$ & $0.00047$ \\
			\hline
			{\small iSVGDMPO} & $0.00012$ & $0.00006$ & $0.00006$ & $0.00007$ \\
			\Xhline{1.1pt}
	\end{tabular}}
\end{table}

In addition, a quantitative comparison among the pCN, rMAP, and iSVGDMPO are given in Table \ref{tab1}. We compute the $\ell^2$-norm differences of the variance and covariance functions on the mesh points obtained by the pCN, rMAP, and iSVGDMPO.
In the table, the notation $\text{cov}_u(x_i,x_{i+k})$ ($k=10, 20, \ldots, 110$) means the covariance function values on the pair of mesh points
$\{(x_i,x_{i+k})\}_{i=1}^{N_g}$. The numbers below this notation are the $\ell^2$ differences between the vectors obtained by the  rMAP and iSVGDMPO with the 
pCN, respectively. All of the $\ell^2$ differences of the iSVGDMPO with the pCN are much smaller than the corresponding values of rMAP, which show the superiority of the iSVGDMPO.

\section{Conclusion}

In this paper, the approximate sampling algorithm is proposed for the infinite-dimensional Bayesian approach. We introduce the Stein operator on Hilbert spaces and show that it is the limit of a particular finite-dimensional version. Besides, we construct the update perturbation of the SVGD on infinite-dimensional space (called iSVGD) by using the properties of operator-valued RKHS. To accelerate the convergence speed of iSVGD, we investigate the change of variables formula and introduced preconditioning operators. As examples, we present the fixed preconditioning operators and mixture preconditioning operators.
Then, we calculate the explicit form of the update directions for the iSVGD with mixture preconditioning operators (iSVGDMPO).
Finally, we apply the constructed algorithms to an inverse problem of the steady state Darcy flow equation.
Comparing with the pCN and rMAP sampling algorithms, we demonstrate by numerical experiments that the proposed algorithms
can generate accurate estimates efficiently.

The iSVGD is analyzed by studying the limiting behavior of the finite-dimensional objects. This work presents an infinite-dimensional version of the approach given in \cite{Wang2019NIPS}.
It is worth mentioning that our results not only provide an infinite-dimensional version
but also indicate that an intuitive trivial generalization of algorithms given in \cite{Wang2019NIPS} may not be suitable since
particles will belong to a set with zero measure. Our results also show that it is necessary
to introduce the parameter $s$, which has not been considered in the existing work.

The current work may be extended to combine the generalizations of the kernel using Hessian operators
in the Wasserstein space \cite{Li2021JMP}. The proposed approach may be combined with other algorithms, such as the accelerated information gradient flows \cite{Wang2020arXiv} and the mean-field type MCMC algorithms \cite{alfredo2020SIADDS}, to generate new and more efficient algorithms. It is also interesting and important to do more theoretical studies, e.g., 
introduce infinite-dimensional Stein geometry \cite{Korba2020NIPS} and 
develop systematic theories of the interacting particle system and the mean field limit equation \cite{Lu2019SIAM}.
We will report the progress on these aspects elsewhere in the future.  



\bibliographystyle{siamplain}
\bibliography{references}
\end{document}


\title[Infinite-dimensional SVGD]
 {Supplementary Material: Stein variational gradient descent on infinite-dimensional space and applications to statistical inverse problems}

\author[J. Jia]{Junxiong Jia}
\address{School of Mathematics and Statistics,
Xi'an Jiaotong University,
Xi'an,
710049, China}
\email{jjx323@xjtu.edu.cn}

\author[P. Li]{Peijun Li}
\address{Department of Mathematics,
Purdue University,
West Lafayette, Indiana, 47907, USA}
\email{lipeijun@math.purdue.edu}

\author[D. Meng]{Deyu Meng}
\address{School of Mathematics and Statistics,
	Xi'an Jiaotong University,
	Xi'an,
	710049, China}
\email{dymeng@mail.xjtu.edu.cn}


\subjclass[2010]{}



\keywords{}

\begin{abstract}
In this supplementary material, we present the details for some of the results and examples given in the main text. 
\end{abstract}

\maketitle

\tableofcontents

\section{Example of kernel satisfying assumption (27) in Theorem 11}\label{ExampleOfThm12}

In this section, we present an example of the kernel that satisfies the assumption (27) given in Theorem 11. Let us recall the assumption (27)
\begin{align}\label{optimalAssum}
\mathbb{E}_{u\sim\mu}\Big[
D_{u'}\bm{K}(u,u')\mathcal{C}_0^{-1/2}g+\sum_{k=1}^{\infty}D_{k}D_{u'}\bm{K}(u,u')e_k
\Big],
\end{align}
which belongs to $\mathcal{L}_1(\mathcal{X},\mathcal{Y})$ for each $u'\in\mathcal{X}$ and $g\in\mathcal{H}^{-s}$.

Taking $\bm{K}(u,u')=K(u,u')\text{Id}$ with
\begin{align*}
K(u,u')=\exp\left( -\frac{1}{h}\|u-u'\|_{\mathcal{X}}^2 \right)
\end{align*}
being a scalar-valued kernel, we have
\begin{align*}
D_{u'}\bm{K}(u,u') & = -\frac{1}{h}\langle u-u', \cdot\rangle_{\mathcal{X}} K(u,u'),  \\
D_{k}D_{u'}\bm{K}(u,u')e_k & =
\frac{1}{h^2}(u_k-u_k')\langle u-u',\cdot\rangle_{\mathcal{X}} K(u,u')e_k.
\end{align*}
Let $\{\varphi_j\}_{j=1}^{\infty}$ be an orthonormal basis of $\mathcal{X}$, and
recall that $\{e_j\}_{j=1}^{\infty}$ represents an orthonormal basis of $\mathcal{Y}$.
Plugging the above formula into (\ref{optimalAssum}), we find that
\begin{align*}
\sum_{j=1}^{\infty}\langle\mathbb{E}_{u\sim\mu}\Big[
D_{u'}\bm{K}(u,u')\mathcal{C}_0^{-1/2}g+\sum_{k=1}^{\infty}D_{k}D_{u'}\bm{K}(u,u')e_k
\Big]\varphi_j, e_j\rangle =
\mathbb{E}_{u\sim\mu}\Big\{ \text{I} + \text{II} \Big\},
\end{align*}
where
\begin{align*}
\text{I} & = -\frac{1}{h}\sum_{j=1}^{\infty}\langle \langle u-u',\varphi_j\rangle_{\mathcal{X}}K(u,u')
\mathcal{C}_0^{-1/2}g,e_{j}\rangle_{\mathcal{Y}}, \\
\text{II} & = \frac{1}{h^2}\sum_{j=1}^{\infty}\langle \sum_{k=1}^{\infty}(u_k-u_k')e_k
\langle u-u',\varphi_j\rangle_{\mathcal{X}},e_j\rangle_{\mathcal{Y}}K(u,u').
\end{align*}
For term I, we have
\begin{align}\label{speciTerm1}
\begin{split}
\text{I} \leq & \frac{1}{h}K(u,u')\Big( \sum_{j=1}^{\infty}\langle u-u',\varphi_j\rangle_{\mathcal{X}}^2 \Big)^{1/2}
\Big( \sum_{j=1}^{\infty}\langle \mathcal{C}_0^{-1/2}g, e_j\rangle \Big)^{1/2} \\
\leq & \frac{C}{h}K(u,u')\|u-u'\|_{\mathcal{X}}\|\mathcal{C}_0^{-1/2}g\|_{\mathcal{Y}} \\
\leq & \frac{C}{h}K(u,u')\|u-u'\|_{\mathcal{X}}\|g\|_{\mathcal{H}^{-s}} < \infty.
\end{split}
\end{align}
For term II, we have
\begin{align}\label{speciTerm2}
\begin{split}
\text{II} \leq & \frac{1}{h^2}K(u,u')\sum_{j=1}^{\infty}\langle u-u',\varphi_j\rangle_{\mathcal{X}}
\langle (u_j-u_j')e_j, e_j\rangle_{\mathcal{Y}} \\
\leq & \frac{1}{h^2}K(u,u')\Big( \sum_{j=1}^{\infty}\langle u-u', \varphi_j\rangle_{\mathcal{X}}^2 \Big)^{1/2}
\Big( \sum_{j=1}^{\infty}(u_j-u_j')^2 \Big)^{1/2}	\\
\leq & \frac{1}{h^2}K(u,u')\|u-u'\|_{\mathcal{X}}^2 < \infty.
\end{split}
\end{align}
Combining estimates (\ref{speciTerm1}) and (\ref{speciTerm2}) yields 
\begin{align}
\sum_{j=1}^{\infty}\langle\mathbb{E}_{u\sim\mu}\Big[
D_{u'}\bm{K}(u,u')\mathcal{C}_0^{-1/2}g+\sum_{k=1}^{\infty}D_{k}D_{u'}\bm{K}(u,u')e_k
\Big]\varphi_j, e_j\rangle < \infty,
\end{align}
which implies that (\ref{optimalAssum}) belongs to $\mathcal{L}_1(\mathcal{X},\mathcal{Y})$.
Taking $\mathcal{X}=\mathcal{H}^{1}, \mathcal{Y}=\mathcal{H}^{-1}$, $s=0$, and projecting all of the quantities to $\mathcal{X}^{N}$, we then obtain the finite-dimensional SVGD as reviewed in Section 2 of the main text.

\section{Implementation details for the mixture preconditioning}\label{ImplementDetailOfMixturePrecond}

In Subsection 3.3, we present the mixture preconditioning operators, which can specify different preconditioning operators for different particles. Here, we provide some more implementation details.

In practice, we approximate the expectation $\mathbb{E}_{u\sim \mu}$
by empirical mean of particles $\{u_i\}_{i=1}^m$. Hence, the formula (56) reduces to
\begin{align}
\begin{split}
\phi^*_{\bm{K}}(\cdot) = \sum_{\ell=1}^{m}w_{\ell}(\cdot)\sum_{j=1}^{m}\Bigg[
& -w_{\ell}(u_j)\bm{K}_{\ell}(u_j,\cdot)(D_{u_j}\Phi(u_j) + \mathcal{C}_0^{-1}u_j)  \\
& + \sum_{k=1}^{\infty}D_k(w_{\ell}(u_j)\bm{K}_{\ell}(u_j,\cdot)e_k) \Bigg].
\end{split}
\end{align}
Taking $\bm{K}_{\ell}$ in (54) with
$$T_\ell = \mathcal{C}_0^{s/2}(D\mathcal{G}(u_\ell)^*\Sigma^{-1}D\mathcal{G}(u_\ell) + \mathcal{C}_0^{-1})^{1/2},$$ we get
\begin{align*}
& \bm{K}_{\ell}(u_j,\cdot)(D_{u_j}\Phi(u_j) + \mathcal{C}_0^{-1}u_j) \\
= &
(D\mathcal{G}(u_\ell)^*\Sigma^{-1}D\mathcal{G}(u_\ell) + \mathcal{C}_0^{-1})^{-1}(D_{u_j}\Phi(u_j) + \mathcal{C}_0^{-1}u_j)
\exp\Big( -\frac{1}{h}\|T_{\ell}(u_j - \cdot)\|_{\mathcal{H}}^2 \Big).
\end{align*}
For the term $D_{k}(w_{\ell}(u_j)\bm{K}_{\ell}(u_j,\cdot)e_k)$, it is clear to note 
\begin{align}
D_{k}(w_{\ell}(u_j)\bm{K}_{\ell}(u_j,\cdot)e_k) =
D_{k}w_{\ell}(u_j)\bm{K}_{\ell}(u_j, \cdot)e_k + w_{\ell}(u_j)D_k\bm{K}_{\ell}(u_j,\cdot)e_k.
\end{align}
For the first term, we have
\begin{align}\label{jisuande1}
\begin{split}
D_{k}w_{\ell}(u_j)\bm{K}_{\ell}(u_j, \cdot)e_k =
& -\langle T_{\ell}(u_j - u_{\ell}), T_{\ell}\varphi_k \rangle_{\mathcal{H}} w_{\ell}(u_j)\bm{K}_{\ell}(u_j, \cdot)e_k \\
& - J_k w_{\ell}(u_j)\bm{K}_{\ell}(u_j, \cdot)e_k,
\end{split}
\end{align}
where
\begin{align}\label{jisuande2}
J_k = \frac{\sum_{\ell'=1}^{m}\langle T_{\ell}(u_j-u_{\ell'}), T_{\ell}\varphi_k\rangle_{\mathcal{H}}
	\exp\Big( -\frac{1}{2}\|T_{\ell'}(u_j - u_{\ell'})\|_{\mathcal{H}}^2 \Big)}
{\sum_{\ell'=1}^{m}\exp\Big( -\frac{1}{2}\|T_{\ell'}(u_j - u_{\ell'})\|_{\mathcal{H}}^2 \Big)}.
\end{align}
For the second term, we have
\begin{align}\label{jisuande3}
w_{\ell}(u_j)D_k\bm{K}_{\ell}(u_j,\cdot)e_k =
-\frac{2}{h}w_{\ell}(u_j) \langle T_{\ell}(u_j-\cdot), T_{\ell}\varphi_k\rangle_{\mathcal{H}} \bm{K}_{\ell}(u_j,\cdot)e_k.
\end{align}
Combining (\ref{jisuande1}), (\ref{jisuande2}) and (\ref{jisuande3}), we obtain
\begin{align}
\sum_{k=1}^{\infty}D_k(w_{\ell}(u_j)\bm{K}_{\ell}(u_j,\cdot)e_k) & =  -\frac{2}{h}
w_{\ell}(u_j)\sum_{k=1}^{\infty}\langle T_{\ell}(u_j-\cdot), T_{\ell}\varphi_k\rangle_{\mathcal{H}} \bm{K}_{\ell}(u_j,\cdot)e_k  \nonumber \\
&
- w_{\ell}(u_j)\sum_{k=1}^{\infty}\langle T_{\ell}(u_j-u_{\ell}), T_{\ell}\varphi_k\rangle_{\mathcal{H}} \bm{K}_{\ell}(u_j,\cdot)e_k \\
& - w_{\ell}(u_j)\sum_{\ell'=1}^m \sum_{k=1}^{\infty}\langle T_{\ell'}(u_j-\cdot),T_{\ell'}\varphi_k\rangle_{\mathcal{H}}\bm{K}_{\ell}(u_j,\cdot)e_k M_{\ell'},
\nonumber
\end{align}
where
\begin{align}
M_{\ell'} = \frac{\exp\Big( -\frac{1}{2}\|T_{\ell'}(u_j - u_{\ell'})\|_{\mathcal{H}}^2 \Big)}
{\sum_{\ell''=1}^{m}\exp\Big( -\frac{1}{2}\|T_{\ell''}(u_j - u_{\ell''})\|_{\mathcal{H}}^2 \Big)}.
\end{align}
For specific examples, we have the explicit form 
\begin{align}
\sum_{k=1}^{\infty}\langle T_{\ell}(u_j-\cdot), T_{\ell}\varphi_k\rangle_{\mathcal{H}} \bm{K}_{\ell}(u_j,\cdot)e_k.
\end{align}
For example, we take $\mathcal{X}$, $\mathcal{Y}$, $\tilde{\mathcal{X}}$, and $\tilde{\mathcal{Y}}$ as in the fixed precondition case
and specify $\bm{K}_{\ell}$ as in (51) with $T$ replaced by $T_{\ell}$. Then, we have
\begin{align}
\begin{split}
& \sum_{k=1}^{\infty}\langle T_{\ell}(u_j-\cdot), T_{\ell}\varphi_k\rangle_{\mathcal{H}} \bm{K}_{\ell}(u_j,\cdot)e_k  \\
& \qquad\qquad\quad
= \exp\Big( -\frac{1}{h}\|T_{\ell}(u_j-\cdot)\|_{\mathcal{H}}^2 \Big) T_{\ell}^{-1}\mathcal{C}_0^s(T_{\ell}^{-1})^*
\mathcal{C}_0^{-s}T_{\ell}^* T_{\ell}(u_j-\cdot).
\end{split}
\end{align}
Hence, it is not required to calculate the orthonormal basis $\{e_k\}_{k=1}^{\infty}$ and $\{\varphi_i\}_{i=1}^{\infty}$
in spaces $\mathcal{Y}$ and $\mathcal{X}$ explicitly in the implementations.

\section{Proof of Theorem 17}\label{InvariantSol}

Blow is the proof of Theorem 17.

\begin{proof}
Denote by $\mathcal{E}(\mathcal{H})$ the set of all the exponential functions and let 
\begin{align}
\varphi_h(x):=e^{i\langle x, h\rangle_{\mathcal{H}}}, \quad x,h\in\mathcal{H}.
\end{align}
By \cite{Prato2004Book}, the function space $\mathcal{E}(\mathcal{H})$ is dense in $L^2(\mathcal{H}, \mu_0)$, 
where $\mu_0$ is the prior measure.  Let $K_n, \varphi_n \in \mathcal{E}(\mathcal{H})$ satisfy
\begin{align*}
\lim_{n\rightarrow\infty}\|K_n - K\|_{W^{1,2}(\mathcal{H},\mu_0)} = 0, \quad
\lim_{n\rightarrow\infty}\|\psi_n - \exp(-\Phi)\|_{L^2(\mathcal{H},\mu_0)} = 0.
\end{align*}
For the prior probability measure, we have $\mathcal{C}_0\varepsilon_k = \lambda_k^2\varepsilon_k$ with $k=1,2,\ldots$,
i.e., $\{\lambda_k^2,\varepsilon_k\}_{k=1}^{\infty}$ is the eigensystem of $\mathcal{C}_0$.
It follows from \cite[Lemma 1.5]{Prato2004Book} that we have 
\begin{align}\label{BasicIntegralbyPartsIdentity}
\begin{split}
\int_{\mathcal{H}}D_k K_n(u-\tilde{u})\psi_{n'}(\tilde{u})\mu_0(d\tilde{u}) = &
-\int_{\mathcal{H}}K_n(u-\tilde{u})D_k\psi_{n'}(\tilde{u})\mu_0(d\tilde{u}) \\
& + \frac{1}{\lambda_k^2}\int_{\mathcal{H}}\tilde{u}_k K_n(u-\tilde{u})\psi_{n'}(\tilde{u})\mu_0(d\tilde{u}),
\end{split}
\end{align}
where $\tilde{u}_k = \langle u, \varepsilon_k \rangle_{\mathcal{H}}, k=1,2,\ldots$. By a simple calculation, we have
\begin{align*}
-\int_{\mathcal{H}}D_k K_n(u-\tilde{u})\psi_{n'}(\tilde{u})\mu_0(d\tilde{u}) =
\int_{\mathcal{H}}K_n(u-\tilde{u})\left( D_k\psi_{n'}(\tilde{u}) + \frac{\tilde{u}_k}{\lambda_k^2} \right)\psi_{n'}(\tilde{u})\mu_0(d\tilde{u}).
\end{align*}
Taking $n'\rightarrow\infty$ in the above equality leads to 
\begin{align}
\int_{\mathcal{H}}D_kK_n(u-\tilde{u})e^{-\Phi(\tilde{u};\bm{d})} -
K_n(u-\tilde{u})D_kV(\tilde{u})e^{-\Phi(\tilde{u};\bm{d})}\mu_0(d\tilde{u}) = 0,
\end{align}
where $V(\cdot)$ is defined in (9). Taking $n\rightarrow\infty$, we arrive at
\begin{align}
\int_{\mathcal{H}} D_kK(u-\tilde{u})e^{-\Phi(\tilde{u};\bm{d})} - K(u-\tilde{u})D_kV(\tilde{u})e^{-\Phi(\tilde{u};\bm{d})}\mu_0(d\tilde{u}) = 0,
\end{align}
which implies
\begin{align}\label{invarformula}
\int_{\mathcal{H}}\langle DK(u-\tilde{u}), D\varphi(u) \rangle_{\mathcal{H}} +
\langle K(u-\tilde{u})DV(\tilde{u}), D\varphi(u) \rangle_{\mathcal{H}} \mu^{\bm{d}}(d\tilde{u}) = 0,
\end{align}
where $\varphi\in\mathcal{E}(\mathcal{H})$ is a test function. Through simple calculations based on (\ref{invarformula}), we further obtain
\begin{align}
\int_{\mathcal{H}} \langle DK*\mu^{\bm{d}}, D\varphi \rangle_{\mathcal{H}} + \langle K*DV\mu^{\bm{d}}, D\varphi \rangle_{\mathcal{H}} \mu^{\bm{d}}(du) = 0,
\end{align}
which implies
\begin{align}
\langle \mu^{\bm{d}}, L(\mu^{d})\varphi \rangle = 0
\end{align}
with $L$ being defined in (59). Recalling the weak form of the equation (58), we complete the proof.
\end{proof}

\section{Proof of Theorem 18}\label{ProofOfThm17}

Let $H^{1}_0(\Omega)$ and $H^{-1}(\Omega)$ be the usual Sobolev spaces.  Consider the boundary value problem 
\begin{align}\label{DarcyEq}
\begin{split}
-\nabla\cdot(e^u\nabla w) & = f \quad\text{in }\Omega, \\
w & = 0 \quad\text{on }\partial\Omega. 
\end{split}
\end{align}
The following estimate is crucial to our proofs.
\begin{theorem}\label{BasicEstThm}
Let $u\in L^{\infty}(\Omega)$ and $f\in H^{-1}(\Omega)$, then Eq. (\ref{DarcyEq}) has a unique solution $w\in H_0^1(\Omega)$ satisfies
\begin{align}\label{BasicEst}
\|w\|_{H_0^1(\Omega)} \leq C e^{\|u\|_{L^{\infty}(\Omega)}} \|f\|_{H^{-1}(\Omega)},
\end{align}
where $C$ is a positive constant independent of $u$.
\end{theorem}

Using Theorem \ref{BasicEstThm}, we can derive the estimates for the adjoint, incremental forward, and incremental adjoint equations. For the adjoint equation, we have
\begin{align}\label{EstAdjoint}
\|p\|_{H_0^1(\Omega)} \leq C e^{\|u\|_{L^{\infty}}} \left\| \sum_{j=1}^{N_d}e^{\frac{1}{2\delta^2}\|x-x_j\|^2}(\ell_{x_j}(w) - d_j) \right\|_{L^2}.
\end{align}
Let $|\Omega|$ be the volume of domain $\Omega$.
Since
\begin{align}
\left\| e^{-\frac{1}{2\delta^2}\|x-x_j\|^2} \right\|_{L^2} \leq |\Omega|^{1/2}
\end{align}
and
\begin{align}
\ell_{x_j}(w) \leq |\Omega|\|w\|_{L^2} \quad \text{for }j=1,\ldots,N_d,
\end{align}
we deduce
\begin{align}
\begin{split}
\|p\|_{H_0^1(\Omega)} & \leq \frac{|\Omega|^{1/2}}{2\pi\delta^2}\sum_{j=1}^{N_d}(\|\bm{d}\| + |\Omega|\|w\|_{L^2}) \\
& \leq \frac{N_d |\Omega|^{1/2}}{2\pi\delta^2}\left(\|\bm{d}\| + C|\Omega|e^{\|u\|_{L^{\infty}}}\|f\|_{H^{-1}}\right) \\
& \leq C \left( 1 + e^{\|u\|_{L^{\infty}}}\|f\|_{H^{-1}} \right),
\end{split}
\end{align}
which implies
\begin{align}\label{EstAdjointFinal}
\|p\|_{H_0^1(\Omega)} \leq C(1+\|f\|_{H^{-1}})e^{2\|u\|_{L^{\infty}}}.
\end{align}
For the incremental forward equation, we have
\begin{align}\label{EstIncrementalForwardFinal}
\begin{split}
\|\hat{w}\|_{H_0^1} & \leq C e^{\|u\|_{L^{\infty}}} \|\nabla\cdot(\hat{u}e^u\nabla w)\|_{H^{-1}} \\
& \leq C e^{\|u\|_{L^{\infty}}} \|\hat{u}e^u\nabla w\|_{L^2} \\
& \leq C e^{2\|u\|_{L^{\infty}}} \|\hat{u}\|_{L^{\infty}} \|\nabla w\|_{L^2} \\
& \leq C e^{3\|u\|_{L^{\infty}}} \|f\|_{H^{-1}} \|\hat{u}\|_{L^{\infty}}.
\end{split}
\end{align}
Similarly, based on Theorem \ref{BasicEstThm}, we have
\begin{align}\label{estIAE}
\begin{split}
\|\hat{p}\|_{H_0^1} & \leq Ce^{\|u\|_{L^{\infty}}}\left[I_1 + I_2\right],
\end{split}
\end{align}
where
\begin{align}
I_1 & = \|\nabla\cdot(\hat{u}e^u\nabla p)\|_{H^{-1}}, \\
I_2 & = \frac{1}{2\pi\delta^2\sigma^2}\sum_{j=1}^{N_d}\|\hat{w}e^{-\frac{1}{2\delta^2}\|x-x_j\|^2} \|_{L^2}.
\end{align}
For $I_1$, we have
\begin{align}\label{estI1}
\begin{split}
I_1 & \leq \|\hat{u}e^u\nabla p\|_{L^2} \leq e^{\|u\|_{L^{\infty}}}\|\hat{u}\|_{L^{\infty}}\|\nabla p\|_{L^2} \\
& \leq C (1+\|f\|_{H^{-1}}) e^{3\|u\|_{L^{\infty}}} \|\hat{u}\|_{L^{\infty}}.
\end{split}
\end{align}
For $I_2$, we get
\begin{align}\label{estI2}
\begin{split}
I_2 & \leq C \|\hat{w}\|_{L^2} \leq C e^{3\|u\|_{L^{\infty}}} \|f\|_{H^{-1}}\|\hat{u}\|_{L^{\infty}}.
\end{split}
\end{align}
Combining (\ref{estIAE}) with estimates of $I_1$ and $I_2$, we obtain the estimate of the
adjoint equation
\begin{align}\label{EstIncrementalAdjointFinal}
\begin{split}
\|\hat{p}\|_{H_0^1} \leq C (1+\|f\|_{H^{-1}})e^{4\|u\|_{L^{\infty}}}\|\hat{u}\|_{L^{\infty}}.
\end{split}
\end{align}

It is clear to note that 
\begin{align}\label{LinftyControl}
\|u\|_{L^{\infty}} \leq C\|u\|_{\mathcal{H}^{1-s}} = C\|u\|_{\mathcal{X}}
\end{align}
holds for $s < 0.5$, which can be deduced based on similar arguments given in \cite[Lemma 16 or Theorem 28]{Dashti2017}.
Since the Hilbert scale is based on the covariance operator $\mathcal{C}_0$ \cite{Agapiou2013SPA,Jia2020IPI},
the space $\mathcal{H}^{1-s}$ is different from the one introduced in \cite{Dashti2017}.
The space $\mathcal{H}^{1-s}$ in our paper is approximately
equal to the space $\mathcal{H}^{2(1-s)}$ defined in \cite{Dashti2017}.
Next, we give the three estimates shown in Theorem 17.

First is to estimate $\Phi(u)$. A simple calculation gives 
\begin{align}
\begin{split}
\Phi(u) & = \frac{1}{2}\|\mathcal{M}(w) - \bm{d}\|^2 \leq C(1 + \|w\|_{L^2})^2  \\
& \leq C(1+\|f\|_{H^{-1}})^2 e^{2\|u\|_{\mathcal{X}}},
\end{split}
\end{align}
where the last inequality used estimates (\ref{BasicEst}) and (\ref{LinftyControl}).

Next is to estimate $D\Phi(u)$.  For any $\tilde{u}\in\mathcal{X}$, we get
\begin{align}\label{estDPhi}
\begin{split}
\langle D\Phi(u), \tilde{u}\rangle & = \int_{\Omega}\tilde{u}e^u\nabla w\cdot\nabla p dx
\leq \|\tilde{u}\|_{L^{\infty}}e^{\|u\|_{L^{\infty}}}\|\nabla w\|_{L^2} \|\nabla p\|_{L^2} \\
& \leq C (1+\|f\|_{H^{-1}})^2 e^{4\|u\|_{L^{\infty}}} \|\tilde{u}\|_{L^{\infty}}	\\
& \leq C (1+\|f\|_{H^{-1}})^2 e^{4\|u\|_{\mathcal{X}}} \|\tilde{u}\|_{\mathcal{X}},
\end{split}
\end{align}
where estimates (\ref{BasicEst}) and (\ref{EstAdjointFinal}) are used to derive the second inequality and
estimate (\ref{LinftyControl}) is used for obtaining the third inequality.
Clearly, it follows from (\ref{estDPhi}) that 
\begin{align}
\|D\Phi(u)\|_{\mathcal{X}^*} \leq C (1+\|f\|_{H^{-1}})^2 e^{4\|u\|_{\mathcal{X}}}.
\end{align}

It is also required to estimate $D^2\Phi(u)$. For any $\tilde{u}, \hat{u}\in \mathcal{X}$, we obtain
\begin{align}\label{estD2PhiZong}
\begin{split}
\langle \langle D^2\Phi(u),\hat{u}\rangle, \tilde{u}\rangle & = I_1 + I_2 + I_3,
\end{split}
\end{align}
where
\begin{align}
I_1 & = \int_{\Omega}\hat{u}\tilde{u}e^u\nabla w\cdot\nabla p dx, \\
I_2 & = \int_{\Omega}\tilde{u}e^u\nabla w\cdot\nabla \hat{p} dx, \\
I_3 & = \int_{\Omega}\tilde{u}e^u\nabla p\cdot\nabla \hat{w} dx.
\end{align}
For $I_1$, we have
\begin{align}\label{estD2PhiI1}
\begin{split}
I_1 & \leq \|\hat{u}\|_{L^{\infty}}\|\tilde{u}\|_{L^{\infty}} e^{\|u\|_{L^{\infty}}} \|\nabla w\|_{L^2} \|\nabla p\|_{L^2} \\
& \leq C \|\hat{u}\|_{L^{\infty}}\|\tilde{u}\|_{L^{\infty}} e^{\|u\|_{L^{\infty}}} (1+\|f\|_{H^{-1}})^2 e^{3\|u\|_{L^{\infty}}}  \\
& \leq C (1+\|f\|_{H^{-1}})^2 e^{4\|u\|_{\mathcal{X}}} \|\hat{u}\|_{\mathcal{X}}\|\tilde{u}\|_{\mathcal{X}},
\end{split}
\end{align}
where (\ref{BasicEst}) and (\ref{EstAdjointFinal}) are used for deriving the second inequality and
(\ref{LinftyControl}) is employed to derive the third inequality. By similar calculations, we obtain from (\ref{EstAdjoint}), (\ref{EstIncrementalForwardFinal}), (\ref{EstIncrementalAdjointFinal}), and (\ref{LinftyControl}) that
\begin{align}\label{estD2PhiI2}
I_2 \leq C(1+\|f\|_{H^{-1}})^2 e^{6\|u\|_{\mathcal{X}}} \|\hat{u}\|_{\mathcal{X}}\|\tilde{u}\|_{\mathcal{X}}
\end{align}
and
\begin{align}\label{estD2PhiI3}
I_3 \leq C (1+\|f\|_{H^{-1}})^2 e^{6\|u\|_{\mathcal{X}}} \|\hat{u}\|_{\mathcal{X}}\|\tilde{u}\|_{\mathcal{X}}.
\end{align}
substituting (\ref{estD2PhiI1}), (\ref{estD2PhiI2}), and (\ref{estD2PhiI3}) into (\ref{estD2PhiZong}), we obtain
\begin{align}\label{estD2PhiIm}
\langle \langle D^2\Phi(u),\hat{u}\rangle, \tilde{u}\rangle \leq C (1+\|f\|_{H^{-1}})^2 e^{6\|u\|_{\mathcal{X}}} \|\hat{u}\|_{\mathcal{X}}\|\tilde{u}\|_{\mathcal{X}}. 
\end{align}
Hence
\begin{align}\label{estD2PhiFinal}
\|D^2\Phi(u)\|_{\mathcal{L}(\mathcal{X},\mathcal{X}^*)} \leq C (1+\|f\|_{H^{-1}})^2 e^{6\|u\|_{\mathcal{X}}},
\end{align}
which completes the proof.

\section{More numerical results for the Darcy flow model}\label{MoreNumResults}

In this section, we provide more numerical results for the Darcy flow model given in Section 4 of the main text.
We intend to answer the following two questions: how do different optimization methods affect the estimates;
how does $s$ influence the convergence speed of the algorithm.

\subsection{Comparison of optimization methods}\label{CompareGDandNewton}

We compare different optimization methods for solving the inverse problem of the Darcy flow equation.
Specifically, we present the maximum a posteriori (MAP) estimate obtained by the gradient descent (GD) algorithm and an inexact matrix-free Newton-conjugate gradient (IMFNCG) algorithm. The latter is suitable for computing large-scale inverse problems.
For more details about the IMFNCG algorithm, we refer to \cite{Tan2013SISC,Wang2018SISC} and references therein.
The step length of GD and IMFNCG are determined by the Armijo line search, and the initial
guess is set to be a zero function.

\begin{figure}[t]
	\centering
	\includegraphics[width=0.6\textwidth]{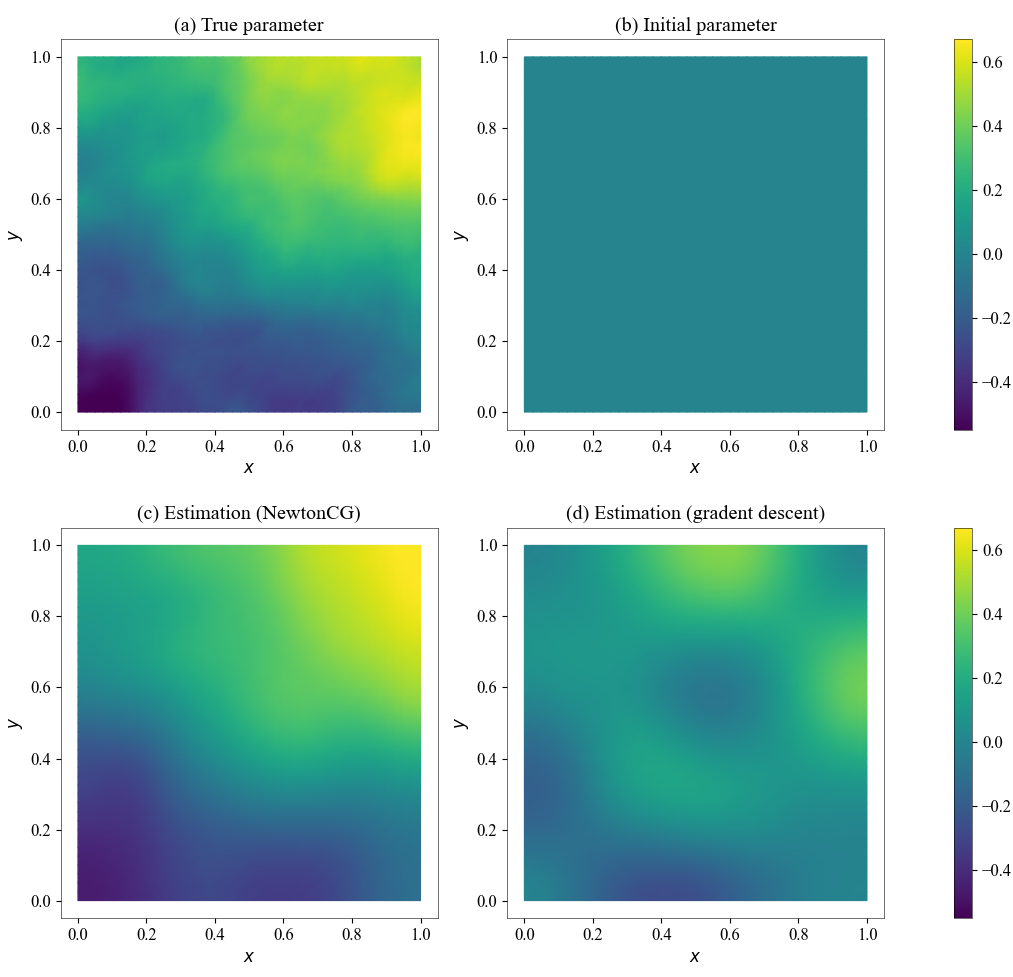}\\
	\vspace*{0.0em}
	\caption{{\small Comparison between the results obtained by the gradient descent (GD) and IMFNCG algorithms
		for the Darcy flow model.
		(a): The background truth; (b): The initial guess of the parameter;
		(c): The MAP estimate obtained by the IMFNCG algorithm with $10$ steps;
		(d): The MAP estimate obtained by the GD algorithm with $1000$ steps.}
	}\label{Fig_compare1}
\end{figure}

\begin{figure}[t]
	\centering
	\includegraphics[width=0.7\textwidth]{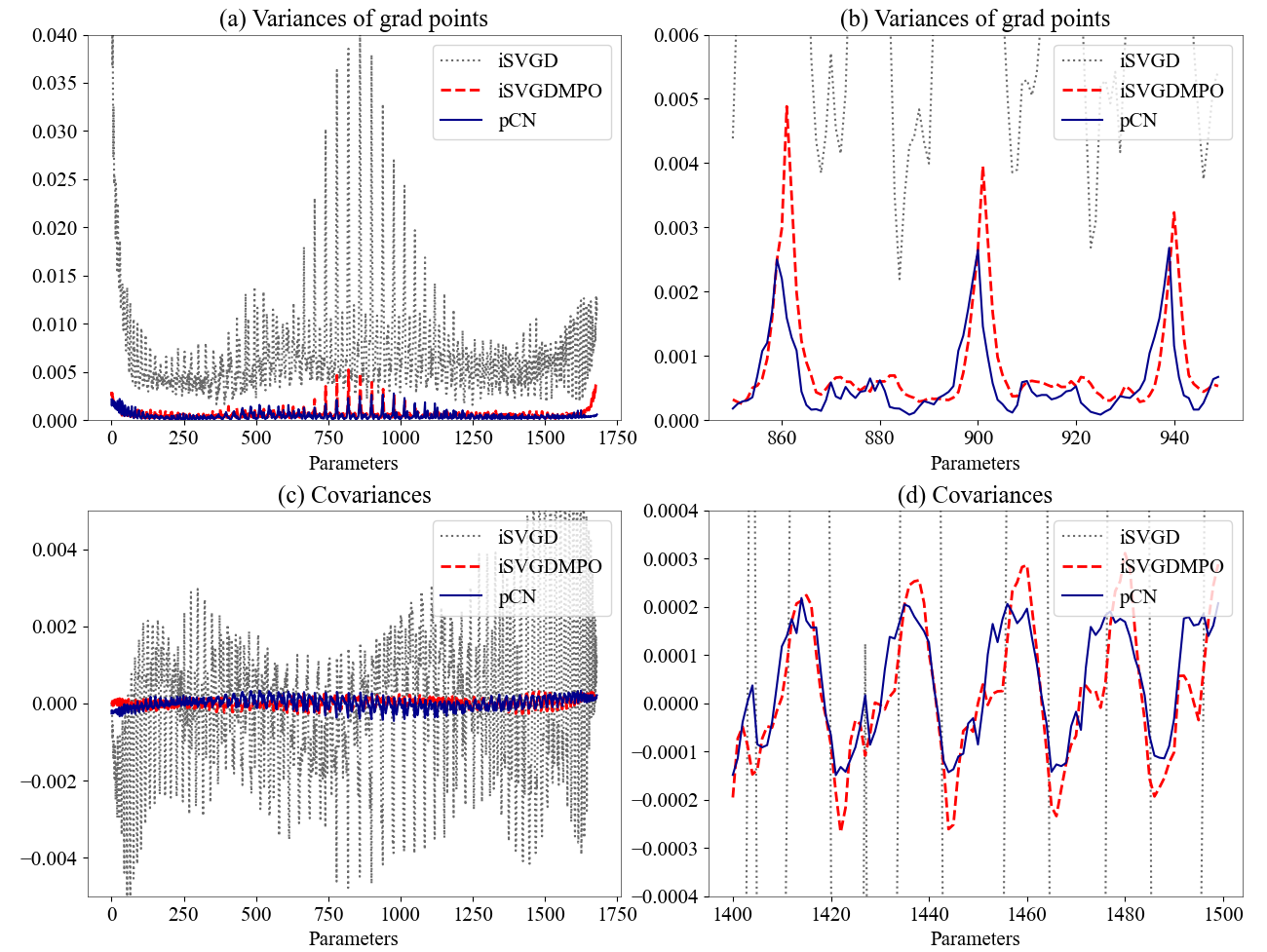}\\
	\vspace*{0.0em}
	\caption{{\small Comparison of the variances and covariance estimated by the pCN, iSVGD ($1000$ iterative steps) and iSVGDMPO ($25$ iterative steps)
		for the Darcy flow model.
		(a): Variances of grad points estimated by pCN, iSVGD and iSVGDMPO (adpative s);
		(b): Local enlarged draw of variances in (a);
		(c): Covariances of point with coordinate $(0.465, 0.035)$ with all other points on the grid estimated by pCN, iSVGD and iSVGDMPO (adpative s);
		(d): Local enlarged draw of covariances in (c).}
	}\label{Fig_compare2}
\end{figure}

Figure \ref{Fig_compare1} shows the estimates obtained by the GD and IMFNCG. On the top left, we show the background truth
function $u$. On the top right, we show the initial zero function. In the second row, we show the MAP estimates obtained
by the IMFNCG and GD algorithms, respectively. It can be seen that the IMFNCG algorithm with only $10$ steps of iteration gives a
reasonable estimate. However, the GD algorithm with Armijo line search cannot provide an accurate estimate even after $1000$ iterative steps. The iSVGD sampling algorithm with no precondition is reduced to the GD algorithm when only one particle is considered. Hence, it is expected that the iSVGD sampling algorithm cannot work well since
particles can hardly concentrate due to the inefficiency of the optimization procedure.
Figure \ref{Fig_compare2} exhibits the estimates of the variance and covariance functions calculated on mesh points by iSVGD and iSVGDMPO when the initial particles are generated by Gaussian approximation of the posterior measure \cite{Tan2013SISC}.
The results shown in Figure \ref{Fig_compare2} confirm our intuition.

In addition, these numerical results verify that it is necessary to introduce the iSVGD with preconditioning operators to
enhance the optimization procedure. Only with an efficient optimization procedure, the concentrate force
(i.e., the first term in the bracket of (40)) and the repulsive force
(i.e., the second term in the bracket of (40)) can sufficiently play their roles to provide accurate samplings.

\subsection{Convergence speed comparison for different values of $s$}\label{compareSpeed}

When choosing a kernel and the prior measure as in Section 4 of the main text, the particles should belong to the Hilbert space $\mathcal{H}^{1-s}$.
From the analysis, we know that $\mu_0(\mathcal{H}^{1-s}) = 0$ or $1$,
when $s=0$ or $s > 0.5$, respectively.
The intuitive idea for specifying the parameter $s$ can be explained as follows:
\begin{enumerate}
	\item The particles should not belong to a set with zero measure, which may lead to inaccurate estimates;
	\item The particles should reside in a small support region of the prior probability measure.
\end{enumerate}
Based on the above two criteria, we may choose $s$ around $0.5$. Here, we provide some numerical results to answer the important question: how does $s$ influence the convergence speed of the iSVGDMPO algorithm.

\begin{figure}[t]
	\centering
	\includegraphics[width=0.9\textwidth]{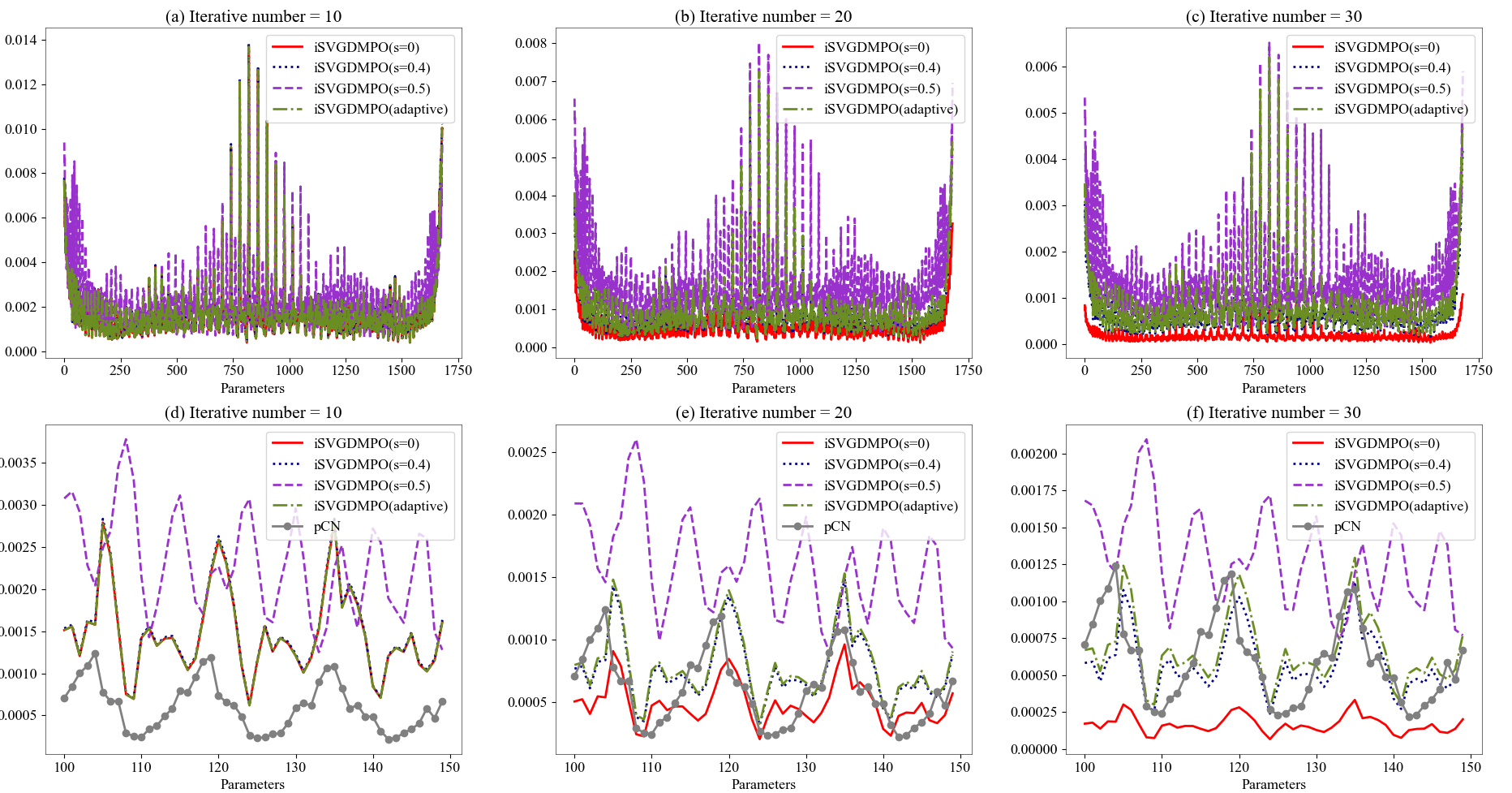}\\
	\vspace*{0.0em}
	\caption{{\small
			For $s=0, 0.4, 0.5$ or adaptively chosen $s$, comparison for the estimated variances of iSVGDMPO when iterative numbers are $10, 20, 30$, respectively.
			(a): The estimated variances when iterative number equal to $10$;
			(b): The estimated variances when iterative number equal to $20$;
			(c): The estimated variances when iterative number equal to $30$;
			(e): The estimated variances (part of the parameters) for the pCN and iSVGDMPO (iterative number equal to $10$);
			(f): The estimated variances (part of the parameters) for the pCN and iSVGDMPO (iterative number equal to $20$);
			(g): The estimated variances (part of the parameters) for the pCN and iSVGDMPO (iterative number equal to $30$).}}
	\label{Fig_compare_ConvSpeed}
\end{figure}

Figure \ref{Fig_compare_ConvSpeed} show the detailed comparisons for the Darcy flow model.
We present the estimated variances when the iterative numbers equal to $10, 20$ and $30$ in (a), (b) and (c) of Figure \ref{Fig_compare_ConvSpeed}, respectively.
In (d), (e) and (f) of Figure \ref{Fig_compare_ConvSpeed}, we depict the estimated variances only for some parameters, which provide
more detailed illustrations.
In these figures, estimated variances for $s=0, 0.4, 0.5$, and the adaptively chosen one are shown, which indicate the following results:
\begin{enumerate}
	\item When the iterative number is smaller than $10$, the convergence speeds for $s=0, 0.4$, and that adaptively chosen are almost the same. The convergence speed for $s=0.5$ is obviously slower than other cases;
	\item When the iterative number approximates $30$, the estimated variances for $s=0$ is much smaller than the estimations given by the pCN and
	iSVGDMPO algorithm with $s=0.4, 0.5$, and the adaptively chosen $s$.
\end{enumerate}
In summary, the convergence speeds are similar for $s=0.4$ or that adaptively chosen.
The obtained estimates, at least for the variance function, are more accurate
when the results of pCN are chosen as the background truth. In the main text, the comparisons for other statistical quantities are given when the parameter $s$ is specified adaptively. When specifying $s=0.5$, the variances will gradually approach the background truth, but the convergence speed seems much slower than $s=0.4$ or the adaptively chosen $s$.

\section{Discussions on the finite- and infinite-dimensional approaches}\label{compareFiniteInfinite}

Since SVGD is constructed usually for the finite-dimensional problems in the field of machine learning, it would be better
for us to provide some detailed explanations about finite- and infinite-dimensional approaches, which should be useful for readers
who are not familiar with the infinite-dimensional approach.  

\subsection{General illustration}\label{compareFiniteInfinite1}

The SVGD algorithm is related to optimization problems since it reduces to an optimization problem for computing maximum a posterior estimate when only one particle is considered. In the following, we firstly recall some discussions from the perspective of PDE-constrained optimization problems. 
For PDE-constrained optimization problems, there are two typical approaches: 
\begin{itemize}
	\item \emph{Discretize-then-optimize}: Discretize the PDEs to formulate a finite dimensional optimization problem, 
	then all of the optimization techniques developed on finite-dimensional space can be applied. 
	\item \emph{Optimize-then-discretize}: Formulate infinite-dimensional optimization problems and construct the optimization
	schemes on some appropriate infinite-dimensional spaces. The discretizations are pushed to the last step to generate 
	practical numerical schemes. 
\end{itemize}
\emph{Discretize-then-optimize} and \emph{optimize-then-discretize} are the finite- and infinite-dimensional approaches
mentioned in the main context, respectively. 
For the advantages of the approach of \emph{optimize-then-discretize}, we refer to page 43--44 of \cite{Reyes2015Book}
and Chapters 2 and 3 of \cite{Hinze2009Book}.
More specifically, the advantages of infinite-dimensional approach are mainly two-folds: 
\begin{itemize}
	\item It is important to have a better understanding of the function space structure of the numerical algorithms in order to design optimal numerical schemes for related PDEs 
	(e.g., when forward PDEs are not self-adjoint, we may need to design certain numerical schemes to calculate 
	forward PDEs and adjoint PDEs then to calculate the gradient). 
	\item The approach is \emph{mesh independent}. The mesh independence implies that 
	the convergence behavior (e.g., convergence rate and number of iterations) of an infinite-dimensional method reflects the 
	behavior of properly discretized problems, when the mesh size is sufficiently small. 
\end{itemize}

Another method for solving inverse problems of PDEs is the Bayesian inverse methods studied in the current work. 
Similar to the PDE-constrained optimization methods, the Bayesian inverse methods also contain two typical approachs:
\begin{itemize}
	\item \emph{Discretize-then-Bayesianize}: The PDEs are initially discretized to approximate the original problem in
	some finite-dimensional space, and the reduced approximate problem is then solved by using the Bayes' method.
	\item \emph{Bayesianize-then-discretize}: The Bayes' formula and algorithms are initially constructed on infinite-dimensional space, and after the infinite dimensional algorithm is built, some finite-dimensional approximation is carried out.
\end{itemize}
\emph{Discretize-then-Bayesianize} and \emph{Bayesianize-then-discretize} are the finite- and in-\\
finite-dimensional approaches mentioned in the main contexts, respectively. 
Similar as the optimization case, these two approaches both have their own advantages and disadvantages, and also either could be suggested to be used dependent on the specific properties of the investigated inverse problems of PDEs. 
By our understanding, the advantages of \emph{Bayesianize-then-discretize} are similar as the case of \emph{Optimize-then-discretize}: 
\begin{itemize}
	\item It is important to have a better understanding of the function space structures in order to design optimal numerical schemes of PDEs, especially when the gradient information is employed. To design sampling algorithms, infinite-dimensional theories will be helpful to design appropriate discretization of probability measures. 
	\item \emph{Bayesianize-then-discretize} approach is \emph{mesh independent}. The sampling efficiency will not highly 
	depend on the dimension of the discretization, which is an important expected property for solving inverse problems of PDEs.  
\end{itemize}
The book \cite{Kaipio2004Book} provides a comprehensive discussions on the finite-dimensional approach, i.e., \emph{discretize-then-Bayesianize} approach. For the infinite-dimensional approach, we refer to 
\cite{Stuart2010AN,Cotter2009IP,Tan2013SISC,Beskos2015SC,Cotter2013SS} and the references there in. 

\subsection{A simple example}\label{compareFiniteInfinite2}
In Subsection \ref{compareFiniteInfinite1}, general discussions on finite- and infinite-dimensional approaches are given, 
which can hardly provide some intuitions on the differences of the numerical schemes. 
In the following, we consider a simple example that illustrates
the implementation differences between ``optimize (Bayesianize)-then-discretize'' and ``discretize-then-optimize (Bayesianize)'' approach. Let us consider the following equation:
\begin{align}\label{toyExample}
	\begin{split}
		\begin{cases}
			-0.1\Delta w + w = u, & \quad \text{in }\Omega,  \\
			u = 0, & \quad\text{on }\partial\Omega,
		\end{cases}
	\end{split}
\end{align}
where $\Omega = [0,1]^2$. 
Denote the forward operator $\mathcal{F}(u):=w$ and the measurement operator $\mathcal{M}(w) := (w(x_1),\ldots,w(x_{N_d}))^T$ 
where $\{x_i\}_{i=1}^{N_d}$ reside in $\Omega$ and $N_d$ is a positive integer. 
Define $\mathcal{G} := \mathcal{M}\circ\mathcal{F}$. We then have the following formulation:
\begin{align}\label{toy1}
	\bm{d} = \mathcal{G}(u) + \bm{\epsilon},
\end{align}
where $\bm{d}$ is the noisy data and $\bm{\epsilon}$ is the random noise. 
The simplest way for estimating $u$ from $\bm{d}$ is to solve the following minimization problem: 
\begin{align}\label{toy2} 
	\min_{u}F(u)
\end{align}
with $F(u) := \frac{1}{2}\|\mathcal{G}(u) - \bm{d}\|_{\ell^2}^2$. 
Now, we employ the finite-element method to discretize the above problem. Denote the finite element mass matrix by $\bm{M}$,
the stiffness matrix of equation (\ref{toyExample}) by $\bm{K}$, and the measurement matrix by $\bm{S}$. 
The forward operator $\mathcal{G}$ then has the following discretized form:
\begin{align}\label{toyDis}
	\bm{d} = \bm{S}\bm{K}^{-1}\bm{M}\bm{u} + \bm{\epsilon},
\end{align}
where $\bm{u}$ is the discretized vector of the function $u$.

\textbf{Discretize-then-optimize (Bayesianize)}: For using \emph{discretize-then-optimize (Bayesianize)} approach, 
we need to formulate the following discrete problem: 
\begin{align}\label{toyDis11} 
	\min_{\bm{u}}\frac{1}{2}\|\bm{S}\bm{K}^{-1}\bm{M}\bm{u} - \bm{d}\|_{\ell^2}^2.
\end{align}
Using the gradient descent method, we obtain the following iterative scheme: 
\begin{align}\label{toyDis12}
	\bm{u}_{k+1} = \bm{u}_k - \gamma (\bm{S}\bm{K}^{-1}\bm{M})^T(\bm{S}\bm{K}^{-1}\bm{M}\bm{u} - \bm{d}),
\end{align}
where $\gamma$ is the step size. 

\textbf{Optimize (Bayesianize)-then-discretize}: For using \emph{optimize (Bayesianize)-then-discretize} approach, 
we need to firstly formulate the infinite-dimensional problem:
\begin{align}\label{toyDis21} 
	\min_{u}\frac{1}{2}\|\mathcal{G}(u) - \bm{d}\|_{\ell^2}^2.
\end{align}
Then we derive the gradient descent iteration on infinite-dimensional space to obtain:
\begin{align}\label{toyDis22} 
	u_{k+1} = u_k - \gamma \mathcal{G}^*(\mathcal{G}(u) - \bm{d}),
\end{align}
where $\mathcal{G}^*$ is the adjoint-operator of $\mathcal{G}$. According to Subsection 3.3 of \cite{Tan2013SISC},
we may obtain the following iterative scheme on finite-dimensional space: 
\begin{align}\label{toyDis23}
	\bm{u}_{k+1} = \bm{u}_k - \gamma \bm{M}^{-1}(\bm{S}\bm{K}^{-1}\bm{M})^T(\bm{S}\bm{K}^{-1}\bm{M}\bm{u} - \bm{d}).
\end{align}
Comparing iterative schemes (\ref{toyDis12}) and (\ref{toyDis23}), we can see the difference. 
At a glance, this is a small difference. However, such a small difference leads to different behaviors of the two iterative schemes. 
We implement the two iterative schemes with different discretized dimensions $d = 20\times 20, 30\times 30, 40\times 40, 50\times 50, 60\times 60$ to visually see such different behavior. 
The step size $\gamma$ is set to be $0.01$ for all of the iterative schemes.
We define the step norm as follows: 
\begin{align}\label{step_norm}
	\text{The $k$-th step norm} = \|u_{k+1} - u_{k}\|_{L^2}.
\end{align}

\begin{figure}[t]
	\centering
	\includegraphics[width=0.8\textwidth]{dim_compare_toy.png}\\
	\vspace*{0.0em}
	\caption{{\small Left: Logarithm of the step norms computed by ``discretize-then-optimize (Bayesianize)'' approach with different discretized dimensions $d = 400, 900, 1600, 2500, 3600$. 
			Right: Logarithm of the step norms computed by ``optimize (Bayesianize)-then-discretize'' approach with different discretized dimensions $d = 400, 900, 1600, 2500, 3600$. }
	}\label{meshIndpendence1}
\end{figure}

In the left of Figure \ref{meshIndpendence1}, we draw the step norms of the iterative scheme (\ref{toyDis12}). We can see that the step norms
decay rapidly when the dimension grows. This indicate that the convergence speed of iterative scheme (\ref{toyDis12})
depends highly on the discretized dimension. 
In contrast, the step norms of the iterative scheme (\ref{toyDis23}) are almost the 
same for different discretizations. From this simple toy example, we can see that the ``discretize-then-optimize (Bayesianize)'' approach can hardly keep the infinite-dimensional natural. Hence, it usually lacks \emph{mesh independence} property. 
However, the ``optimize (Bayesianize)-then-discretize'' approach pushes the discretization implementation to the final step which makes it easier to catch the infinite-dimensional natural of the inverse problems of PDEs. The finally obtained algorithm usually 
has \emph{mesh independence} property, which is important for solving inverse problems of PDEs. 

\subsection{Mesh independence of iSVGD}\label{compareFiniteInfinite3}
In Subsection \ref{compareFiniteInfinite2}, we just provide a simple example. 
For the proposed iSVGD, we need more standard techniques that can be found in some typical literatures \cite{Tan2013SISC,Ghattas2021ActaNumerica,Reyes2015Book,Spantini2015SISC,Wang2018SISC,Thanh2016IPI}.
The lecture notes provided in ``https://uvilla.github.io/inverse17/'' are also beneficial for taking implementations. 

Now, let us illustrate that the proposed iSVGD algorithm possesses the \emph{mesh independence} property. 
That is to say, if the finite element mesh is refined, we indeed need more computational resources since the computations of each partial differential equation are more expensive. However, it might not need more iterations and particles when the finite element mesh is refined since discrete problems derived by refined mesh also approximate the infinite-dimensional formulation. 
For clearly illustrating this, we choose different discretized grids such that the dimensions of the function parameter are $d=20\times 20, 30\times 30, 40\times 40, 50\times 50, 60\times 60$. Using the same settings as in Section 4 of the main text, we only change the discretized dimension to see how discretized dimensions affect the behavior of the algorithm. In Figure \ref{meshIndpendence2}, we show the numerical results which demonstrate the \emph{mesh independence} as expected for \emph{Bayesianize-then-discretize} approach. 

Specifically speaking, we draw the variance functions with different discretized dimensions in the left of Figure \ref{meshIndpendence2}.
The variance functions are calculated by the iSVGDMPO with discretized dimensions $d = 400, 900, 1600, 2500, 3600$. 
When the algorithm generates the final particles, we calculate the variance functions and project the estimated variance functions 
on a mesh with dimension $d=20\times 20$. Then, we draw part of the grid point values of the variance functions calculated by different meshes. From the figure, it can be seen that the grid point values are similar.
This validates that the estimated variance function obtained by iSVGDMPO is not sensitive to one particular discretization.
Similar to other \emph{mesh independence} methods such as rMAP used for comparison in our numerical experiments, it may be difficult to obtain exactly the same values due to the quantities being evaluated approximately, especially the gradients and Hessian operators, are not evaluated accurately. For discretize-first type methods, we can calculate the gradients of the discretized system exactly. However, the gradients and Hessian operators defined on infinite-dimensional space could only be calculated approximately.   

\begin{figure}[t]
	\centering
	\includegraphics[width=0.8\textwidth]{dim_compare.png}\\
	\vspace*{0.0em}
	\caption{{\small Left: Pointwise sample variances computed by different discretized dimensions $d = 400, 900, 1600, 2500, 3600$
			(All of the variance functions are projected on a grid with $d=400$ for comparison). 
			Right: Decay of the averaged step norm $\frac{1}{m}\sum_{i=1}^m\|u_i^{\ell+1} - u_i^{\ell}\|_{L^2}$ w.r.t. the number
			of iterations for different discretized dimensions $d = 400, 900, 1600, 2500, 3600$. }
	}\label{meshIndpendence2}
\end{figure}

In the right of Figure \ref{meshIndpendence2}, we draw the averaged step norm defined as follows:
\begin{align}
	\frac{1}{m}\sum_{i=1}^{m}\|u_i^{\ell+1} - u_i^{\ell}\|_{L^2},
\end{align}
where $u_i^{\ell}$ stands for the $i$th particle at the $\ell$th iteration and $m$ is the number of particles. 
Obviously, the averaged step norms are similar for different discretized dimensions. 
It is evident that the curves under different discretized dimensions can hardly be distinguished, indicating that the algorithm has \emph{mesh independence} property. 
The convergence speed is not affected by discretized dimensions, which is not true for many algorithms developed 
under the finite-dimensional setting. 

At last, we should admit that more theoretical works are needed to ensure the \emph{mesh independence} property. 
Specifically speaking, we may need to do further research based on Subsection 3.4 in the main text on infinite-dimensional 
particle interacting system and the measure-valued evolution equation. The well-posedness of these complicated equations
should be proved and a theorem like Theorem 2.7 in \cite{Lu2019SIAM} needs to be established. 
Along this direction, we may consult to the studies on the semilinear Mckean--Vlasov stochastic evolution equation in 
Hilbert space \cite{Ahmed1995SPTA} and the theoretical analysis of the pCN algorithm \cite{Pillai2014SPDE}.  

\section{Numerical results for the Helmholtz equation}\label{Helmholtz}

In this section, we present numerical experiments for the Helmholtz equation
\begin{align}\label{equHelmholtz}
\begin{split}
-\Delta w - e^{2u}w & = 0 \text{ in }\Omega, \\
\frac{\partial w}{\partial\bm{n}} & = g \text{ in }\partial\Omega,
\end{split}
\end{align}
where $w$ is the acoustic field, $u$ is the logarithm of the distributed wave number field on $\Omega$ ($\Omega$ is a bounded domain), $\bm{n}$ is the unit outward normal on $\partial\Omega$, and $g$ is the prescribed Neumann source on the boundary.
The boundary value problem \eqref{equHelmholtz} may not have a unique solution due to possible resonances \cite{Colton2019Book}.  Hence, we can hardly verify Assumption 6 for this example.
However, this model was studied for the randomized maximum a posteriori (rMAP) method \cite{Wang2018SISC},
which is an approximate method used for our comparison in the main text.
From the proof in Subsection \ref{ProofOfThm17}, we may verify Assumption 6 under more suitable settings for the inverse medium scattering problem, e.g., Lemma 2.3  in \cite{Bao2005IP} gives a similar estimate to the Darcy flow model.

\begin{figure}[t]
	\centering
	\includegraphics[width=0.95\textwidth]{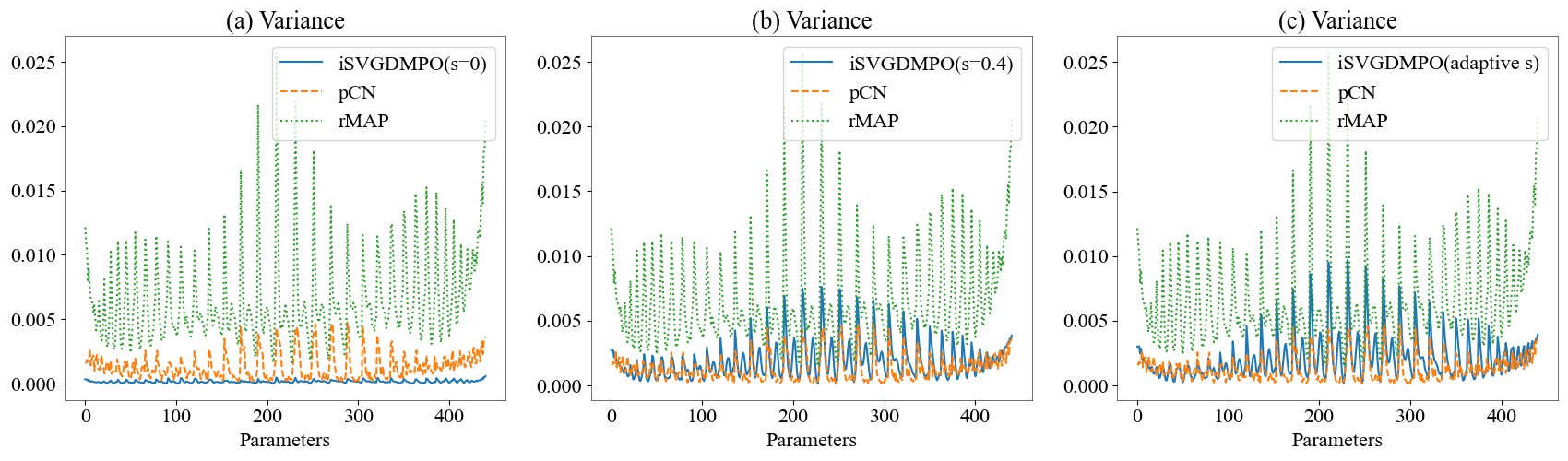}\\
	\vspace*{0.0em}
	\caption{{\small Comparison of the variances estimated by the pCN, rMAP, iSVGDMPO with different $s$ for the Helmholtz equation model.
		(a): Variances estimated by pCN, rMAP, and iSVGDMPO ($s=0$);
		(b): Variances estimated by pCN, rMAP, and iSVGDMPO ($s=0.4$);
		(c): Variances estimated by pCN, rMAP, and iSVGDMPO (adpative s).}
	}\label{Fig_compare_s_Helmholtz}
\end{figure}

\begin{figure}[t]
	\centering
	\includegraphics[width=0.75\textwidth]{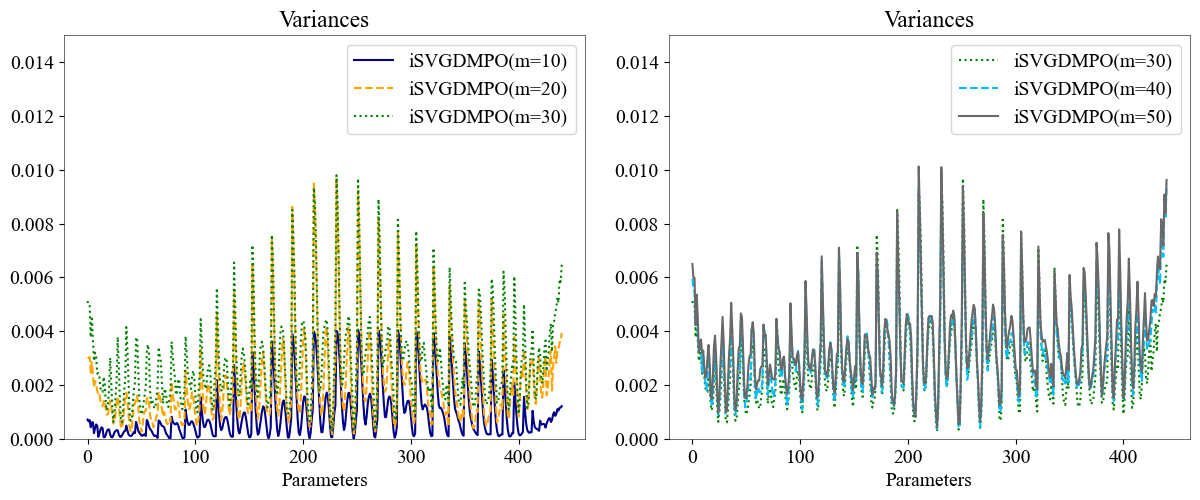}\\
	\vspace*{0.0em}
	\caption{{\small Comparison of the variances estimated by the iSVGDMPO when $s=10,20,30,40,50$ for the Helmholtz equation model.}
	}\label{Fig_compare_sample_number_Helmholtz}
\end{figure}

Basic settings and the finite-element discretization are similar to the Darcy flow model considered in the main text.
The only difference is that the measurement data are collected on the boundary of the domain, i.e.,
$x_i \in\partial\Omega$ for $i=1,\ldots,N_d$. Figure \ref{Fig_compare_s_Helmholtz} shows the estimates of the
variances obtained by the pCN, rMAP, and iSVGDMPO with parameter $s=0,0.4$, or choosing adaptively according to formula (68) in the main text. Similarly, the estimated variances are too small when $s=0$, which implies that the particles
are concentrated on a small set. When $s$ is taken as $0.4$ or chosen adaptively, we obtain similar estimates, which is more similar to the baseline provided by pCN compared with the estimates obtained by the rMAP. As in the main text, we use the empirical adaptive strategy to specify the parameter $s$ in the following.

As for the sample numbers, we also compare the estimated variances when the particle number $m$ equals to $10, 20, 30, 40$, and $50$. On the left and right in Figure \ref{Fig_compare_sample_number_Helmholtz}, we show the results obtained when $m=10, 20, 30$
and $m = 30, 40, 50$, respectively. Obviously, we find that $m=10$ is not enough to give reliable estimates and
the estimated variance functions are similar when $m = 30,40,50$. Hence, for the Helmholtz problem, it is enough to take $m=20$ or $30$ for our numerical examples, which attains a fine balance between efficiency and accuracy.

\begin{figure}[t]
	\centering
	\includegraphics[width=0.95\textwidth]{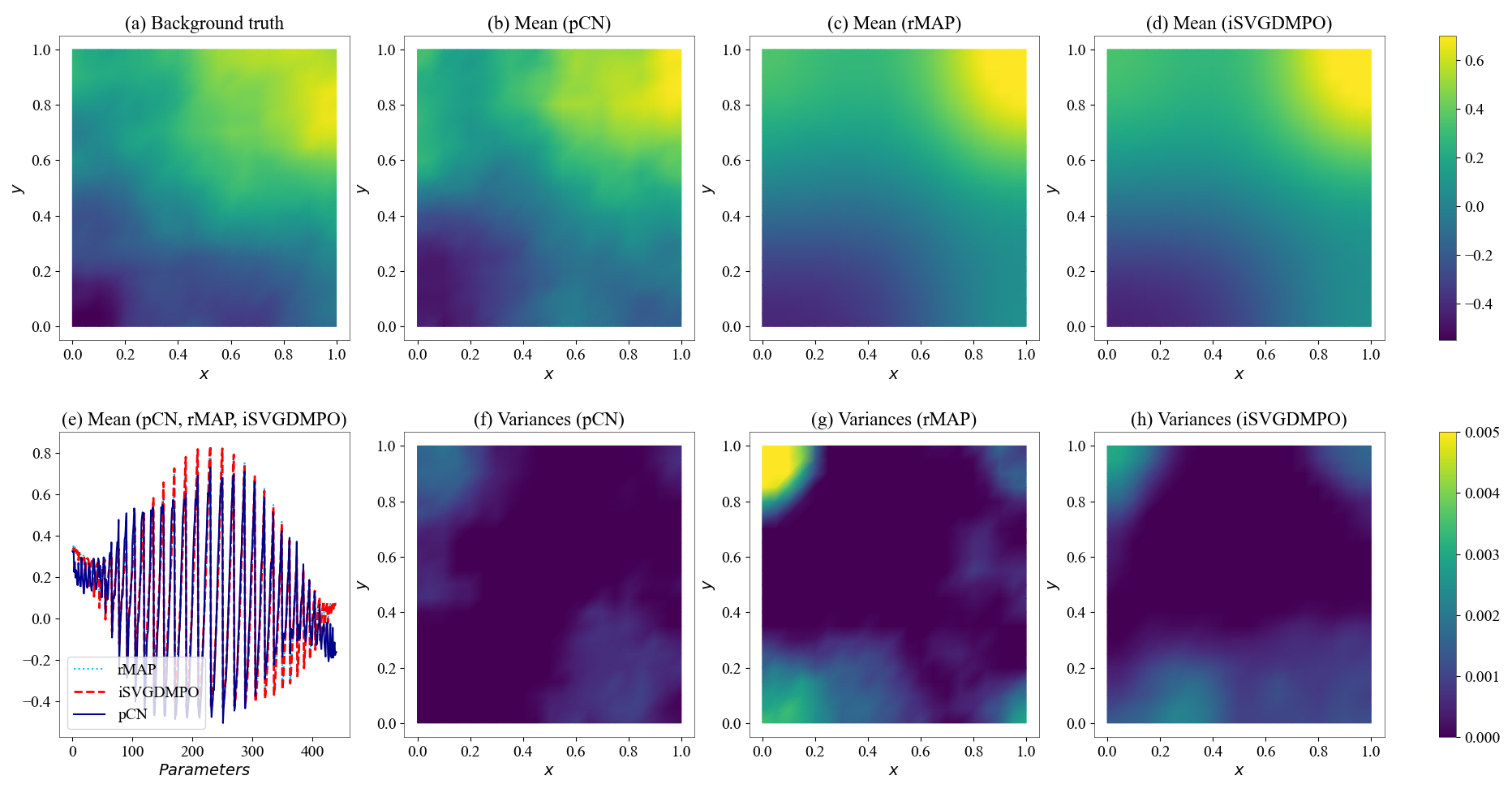}\\
	\vspace*{0.0em}
	\caption{{\small The background truth and estimated mean and variance functions by pCN, rMAP, and iSVGDMPO for the Helmholtz equation model.
		(a): Background truth; (b): Estimated mean function by pCN;
		(c): Estimated mean function by rMAP; (d): estimated mean function by iSVGDMPO;
		(e): Estimated mean function on mesh points by pCN (blue solid line), rMAP (light blue dotted line), and iSVGDMPO (red dashed line);
		(f): Estimated variances by pCN; (g): Estimated variances by rMAP;
		(h): Estimated variances by iSVGDMPO.}
	}\label{Fig_truth_mean_variances_Helm}
\end{figure}

\begin{figure}[htbp]
	\centering
	\includegraphics[width=0.95\textwidth]{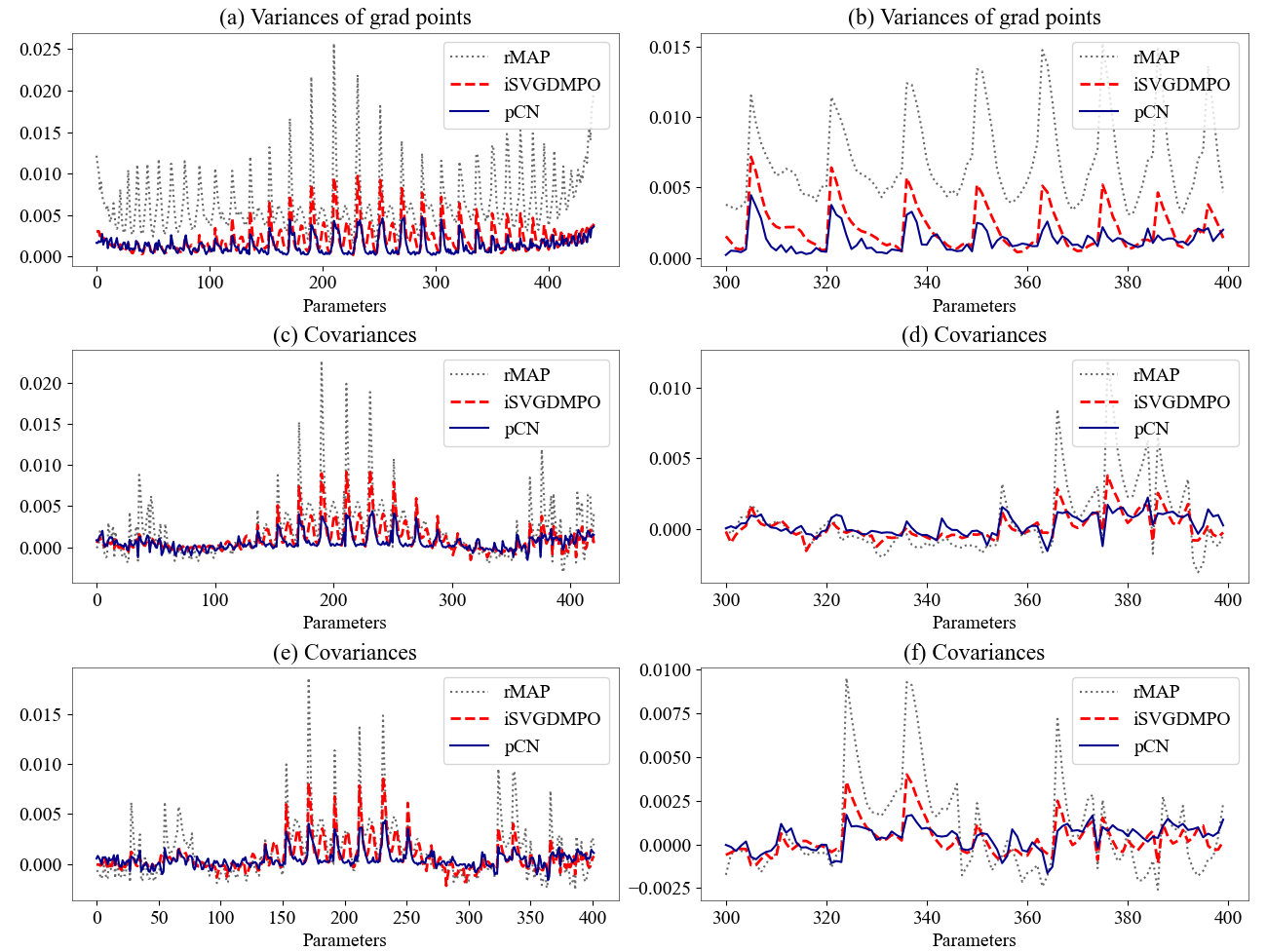}\\
	\vspace*{0.0em}
	\caption{{\small The estimated variances and covariances by the pCN (blue solid line), rMAP (gray dotted line), and iSVGDMPO (red dashed line).
		(a): Estimated variances $\{\text{var}_u(x_i)\}_{i=1}^{N_g}$ on all mesh points;
		(b): Estimated variances for mesh points with indexes from $300$ to $400$ (show details);
		(c): Estimated covariances $\{\text{cov}_u(x_i, x_{i+20})\}_{i=1}^{N_g-20}$ on mesh point pairs $\{(x_i,x_{i+20})\}_{i=1}^{N_g-20}$;
		(d): Estimated covariances shown in (c) with indexes from $300$ to $400$ (show details);
		(e): Estimated covariances $\{\text{cov}_u(x_i, x_{i+40})\}_{i=1}^{N_g-40}$ on mesh point pairs $\{(x_i,x_{i+40})\}_{i=1}^{N_g-40}$;
		(f): Estimated covariances shown in (e) with indexes from $300$ to $400$ (show details).}
	}\label{Fig_variance_covariance_Helm}
\end{figure}

For the following numerical experiments, we take $m = 30$ and set the parameter $s$ by the
empirical strategy (68) as presented in the main text.
In Figure \ref{Fig_truth_mean_variances_Helm}, we demonstrate the background truth and the estimated mean and variance functions
obtained by the pCN, rMAP, and iSVGDMPO, respectively. The iterative number of iSVGDMPO is set to be $30$.
The same observation can be made from the results. 
The mean functions obtained by the rMAP and iSVGDMPO are similar,
which are slightly smoother than the one obtained by the pCN algorithm.
Regarding the variance function, it can be seen from (f), (g), and (h) of the figure that the iSVGDMPO gives more reliable estimates than the rMAP does. 

Now, we provide some more comparisons of statistical quantities among the results obtained by the 
pCN, rMAP, and iSVGDMPO. Similarly, we compute variance and covariance functions on the mesh points
and exhibit the results in Figure \ref{Fig_variance_covariance_Helm}.
In all of the subfigures in Figure \ref{Fig_variance_covariance_Helm},
the estimates obtained by the pCN, rMAP, and iSVGDMPO are drawn in blue solid line, gray dotted line, and red dashed line, respectively. All the notations here are the same as those used in the main text.
We can also obtain the same conclusions from the results:
the estimates obtained by the iSVGDMPO are visually more similar to
the estimates provided by the pCN compared with the results obtained by the rMAP.

\begin{table}[htbp]\label{tab_Helmholtz}
	\caption{{\small The $\ell^2$-norm error of variance and covariance functions on mesh points for the rMAP and iSVGDMPO (estimates of the pCN are seen as the background truth)}}
	\begin{tabular}{cccc}
		\Xhline{1.1pt}
		& $\text{var}_u(x_i)$ & $\text{cov}_u(x_i,x_{i+5})$ & $\text{cov}_u(x_i,x_{i+10})$    \\
		\hline
		{\small rMAP} & $0.01525$ & $0.00155$ & $0.00237$    \\
		\hline
		{\small iSVGDMPO} & $0.00092$ & $0.00026$ & $0.00063$    \\
		\Xhline{1.1pt}
		& $\text{cov}_u(x_i,x_{i+15})$ & $\text{cov}_u(x_i,x_{i+20})$ & $\text{cov}_u(x_i,x_{i+25})$   \\
		\hline
		{\small rMAP} & $0.00295$ & $0.00353$ & $0.00153$   \\
		\hline
		{\small iSVGDMPO} & $0.00036$ & $0.00059$ & $0.00035$    \\
		\Xhline{1pt}	
		& $\text{cov}_u(x_i,x_{i+30})$ & $\text{cov}_u(x_i,x_{i+35})$ & $\text{cov}_u(x_i,x_{i+40})$   \\
		\hline
		{\small rMAP} & $0.00148$ & $0.00154$ & $0.00219$   \\
		\hline
		{\small iSVGDMPO} & $0.00055$ & $0.00037$ & $0.00048$   \\
		\Xhline{1.1pt}
	\end{tabular}
\end{table}

Besides these visual comparisons, a quantitative comparison of the differences among the pCN, rMAP, and iSVGDMPO are also given in Table 1. Again, all the notations have the same meaning as those used in the main text. It can be seen from Table 1 that
all the $\ell^2$-norm differences of the iSVGDMPO with the pCN are evidently smaller than the corresponding values of the rMAP.



\bibliographystyle{plain}
\bibliography{references}

%% file: ex_shared.tex
\usepackage{lipsum}
\usepackage{amsfonts}
\usepackage{graphicx}
\usepackage{epstopdf}
\usepackage{algorithmic}
\ifpdf
\DeclareGraphicsExtensions{.eps,.pdf,.png,.jpg}
\else
\DeclareGraphicsExtensions{.eps}
\fi

\newcommand{\TheTitle}{Stein variational gradient descent on infinite-dimensional space and applications to statistical inverse problems}
\newcommand{\TheAuthors}{J. Jia, P. Li, and D. Meng}
\newcommand{\TheRuningTitle}{SVGD on infinite-dimensional space}

\headers{\TheRuningTitle}{\TheAuthors}

\title{{\TheTitle}\thanks{
		The authors would like to thank associate professor Wei Gong for helpful discussions that greatly improved the manuscript.
		The first author was supported by the National Key R\&D Program of the Ministry of Science and Technology of China grant No. 2020YFA0713403 and the NSFC (Grant Nos. 11871392, 12090020, and 12090021). 
		The second author was supported in part by the NSF grant DMS-1912704. 
		The third author was supported by the NSFC (Grant Nos. 61721002 and U1811461), the Major Key Project of PCL (PCL2021A12) and the Macao Science and Technology Development Fund grant No. 061/2020/A2}}

\author{
	Junxiong Jia\thanks{School of Mathematics and Statistics,
		Xi'an Jiaotong University,
		Xi'an
		710049, China
		(\email{jjx323@xjtu.edu.cn}).}
	\and 
	Peijun Li\thanks{Department of Mathematics,
		Purdue University,
		West Lafayette, Indiana, 47907, USA
		(\email{lipeijun@math.purdue.edu}).}
	\and 
	Deyu Meng\thanks{Corresponding author. 
		School of Mathematics and Statistics and Ministry of  Education Key Lab of  Intelligent Networks and Network Security, 
		Xi’an Jiaotong University,  
		Xi’an, 710049, China, 
		and Macau Institute of Systems Engineering, 
		Macau University of Science and Technology, 
		Taipa, Macau
		(\email{dymeng@mail.xjtu.edu.cn}).}
}

\usepackage{amsopn}


%% file: ex_article_v3.bbl
\begin{thebibliography}{10}

\bibitem{Arridge2019AN}
{\sc S.~Arridge, P.~Maass, O.~\"{O}ktem, and C.-B. Sch\"{o}nlieb}, {\em Solving
  inverse problems using data-driven models}, Acta Numer., 28 (2019),
  pp.~1--174.

\bibitem{Beskos2015SC}
{\sc A.~Beskos, A.~Jasra, E.~A. Muzaffer, and A.~M. Stuart}, {\em Sequential
  {M}onte {C}arlo methods for {B}ayesian elliptic inverse problems}, Stat.
  Comput., 25 (2015), p.~727–737.

\bibitem{Bishop2006PRML}
{\sc C.~M. Bishop}, {\em Pattern Recognition and Machine Learning},
  Springer-Verlag, New York, NY, USA, 2006.

\bibitem{Tan2013SISC}
{\sc T.~Bui-Thanh, O.~Ghattas, J.~Martin, and G.~Stadler}, {\em A computational
  framework for infinite-dimensional {B}ayesian inverse problems part {I}:
  {T}he linearized case, with application to global seismic inversion}, SIAM J.
  Sci. Comput., 35 (2013), pp.~A2494--A2523.

\bibitem{Burger2014IP}
{\sc M.~Burger and F.~Lucka}, {\em Maximum a posteriori estimates in linear
  inverse problems with log-concave priors are proper {Bayes} estimators},
  Inverse Probl., 30 (2014), p.~114004.

\bibitem{Thanh2016IPI}
{\sc T.~But-Thanh and Q.~P. Nguyen}, {\em {FEM}-based discretization-invariant
  {MCMC} methods for {PDE}-constrained {B}ayesian inverse problems}, Inverse
  Probl. Imag., 10 (2016), pp.~943--975.

\bibitem{Carmeli2006AA}
{\sc C.~Carmeli, E.~D. Vito, and A.~Toigo}, {\em Vector valued reproducing
  kernel {H}ilbert spaces of integrable functions and {M}ercer theorem}, Anal.
  Appl., 4 (2006), pp.~377--408.

\bibitem{Carmeli2010AA}
{\sc C.~Carmeli, E.~D. Vito, and A.~Toigo}, {\em Vector-valued reproducing
  kernel {H}ilbert spaces and universality}, Anal. Appl., 8 (2010), pp.~19--61.

\bibitem{Carvalho2020CVPR}
{\sc E.~D.~C. Carvalho, R.~Clark, A.~Nicastro, and P.~H.~J. Kelly}, {\em
  Scalable uncertainty for computer vision with functional variational
  inference}, in CVPR, 2020, pp.~12003--12013.

\bibitem{Chen2020SISC}
{\sc P.~Chen and O.~Ghattas}, {\em Stein variational reduced basis {B}ayesian
  inversion}, SIAM J. Sci. Comput., 43 (2021), pp.~A1163--A1193.

\bibitem{Chen2019NIPS}
{\sc P.~Chen, K.~Wu, J.~Chen, T.~O'Leary-Roseberry, and O.~Ghattas}, {\em
  Projected {S}tein variational {N}ewton: a fast and scalable {B}ayesian
  inference method in high dimensions}, in NeurIPS, vol.~32, 2019.

\bibitem{Cotter2009IP}
{\sc S.~L. Cotter, M.~Dashti, J.~C. Robinson, and A.~M. Stuart}, {\em Bayesian
  inverse problems for functions and applications to fluid mechanics}, Inverse
  Probl., 25 (2009), p.~115008.

\bibitem{Cotter2013SS}
{\sc S.~L. Cotter, G.~O. Roberts, A.~M. Stuart, and D.~White}, {\em {MCMC}
  methods for functions: modifying old algorithms to make them faster}, Stat.
  Sci., 28 (2013), pp.~424--446.

\bibitem{Cui2016JCP}
{\sc T.~Cui, K.~J.~H. Law, and Y.~M. Marzouk}, {\em Dimension-independent
  likelihood-informed {MCMC}}, J. Comput. Phys., 304 (2016), pp.~109--137.

\bibitem{DaPrato1996Book}
{\sc G.~DaPrato and J.~Zabczyk}, {\em Stochastic Equations in Infinite
  Dimensions}, Cambridge University Press, Cambridge, 1992.

\bibitem{Dashti2017}
{\sc M.~Dashti and A.~M. Stuart}, {\em The {Bayesian} approach to inverse
  problems}, Handbook of Uncertainty Quantification,  (2017), pp.~311--428.

\bibitem{Detommaso2018NIPS}
{\sc G.~Detommaso, T.~Cui, A.~Spantini, and Y.~Marzouk}, {\em A {S}tein
  variational {N}ewton method}, in NeurIPS, vol.~32, 2018.

\bibitem{Duncan2019arXiv}
{\sc A.~Duncan, N.~N\"{u}sken, and L.~Szpruch}, {\em On the geometry of {S}tein
  variational gradient descent}.
\newblock arXiv:1912.00894, 2019.

\bibitem{Engl1996Book}
{\sc H.~W. Engl, M.~Hanke, and A.~Neubauer}, {\em Regularization of Inverse
  Problems}, Springer, Netherlands, 1996.

\bibitem{Feng2018SISC}
{\sc Z.~Feng and J.~Li}, {\em An adaptive independence sampler {MCMC} algorithm
  for {B}ayesian inferences of functions}, SIAM J. Sci. Comput., 40 (2018),
  pp.~A1310--A1321.

\bibitem{Fichtner2011Book}
{\sc A.~Fichtner}, {\em Full Seismic Waveform Modelling and Inversion},
  Springer, New York, 2011.

\bibitem{alfredo2020SIADDS}
{\sc A.~Garbuno-Inigo, F.~Hoffmann, W.~C. Li, and A.~M. Stuart}, {\em
  Interacting {L}angevin diffusions: gradient structure and ensemble {K}alman
  sampler}, SIAM J. Appl. Dyn. Syst., 19 (2020), pp.~412--441.

\bibitem{Guha2015JCP}
{\sc N.~Guha, X.~Wu, Y.~Efendiev, B.~Jin, and B.~K. Malick}, {\em A variational
  {B}ayesian approach for inverse problems with skew-t error distribution}, J.
  Comput. Phys., 301 (2015), pp.~377--393.

\bibitem{Helin2015IP}
{\sc T.~Helin and M.~Burger}, {\em Maximum a posteriori probability estimates
  in infinite-dimensional {B}ayesian inverse problems}, Inverse Probl., 31
  (2015), p.~085009.

\bibitem{Jia2020IPI}
{\sc J.~Jia, J.~Peng, and J.~Gao}, {\em Posterior contraction for empirical
  {B}ayesian approach to inverse problems under non-diagonal assumption},
  Inverse Probl. Imag., 15 (2020), pp.~201--228.

\bibitem{Jia2019IP}
{\sc J.~Jia, B.~Wu, J.~Peng, and J.~Gao}, {\em Recursive linearization method
  for inverse medium scattering problems with complex mixture {G}aussian error
  learning}, Inverse Probl., 35 (2019), p.~075003.

\bibitem{Jia2018JFA}
{\sc J.~Jia, S.~Yue, J.~Peng, and J.~Gao}, {\em Infinite-dimensional {B}ayesian
  approach for inverse scattering problems of a fractional {H}elmholtz
  equation}, J. Funct. Anal., 275 (2018), pp.~2299--2332.

\bibitem{Jia2020SISC}
{\sc J.~Jia, Q.~Zhao, D.~Meng, and Y.~Leung}, {\em Variational {B}ayes' method
  for functions with applications to some inverse problems}, SIAM J. Sci.
  Comput., 43 (2021), pp.~A355--A383.

\bibitem{Jin2012JCP}
{\sc B.~Jin}, {\em A variational {B}ayesian method to inverse problems with
  implusive noise}, J. Comput. Phys., 231 (2012), pp.~423--435.

\bibitem{Jin2010JCP}
{\sc B.~Jin and J.~Zou}, {\em Hierarchical {B}ayesian inference for ill-posed
  problems via variational method}, J. Comput. Phys., 229 (2010),
  pp.~7317--7343.

\bibitem{Kadri2016JMLR}
{\sc H.~Kadri, E.~Duflos, P.~Preus, S.~Canu, A.~Rakotomamonjy, and
  J.~Audiffren}, {\em Operator-valued kernels for learning from functional
  response data}, J. Mach. Learn. Res., 17 (2016), pp.~1--54.

\bibitem{Kaipio2004Book}
{\sc J.~Kaipio and E.~Somersalo}, {\em Statistical and Computational Inverse
  Problems}, Springer-Verlag, New York, 2005.

\bibitem{Korba2020NIPS}
{\sc A.~Korba, A.~Salim, M.~Arbel, G.~Luise, and A.~Gretton}, {\em A
  non-asymptotic analysis for {S}tein variational gradient descent}, in
  NeurIPS, vol.~33, 2020.

\bibitem{Lei2020Bernoulli}
{\sc J.~Lei}, {\em Convergence and concentraction of empirical measures under
  wasserstein distance in unbounded functional space}, Bernoulli, 26 (2020),
  pp.~767--798.

\bibitem{Levin2017Book}
{\sc D.~A. Levin, Y.~Peres, and E.~L. Wilmer}, {\em Markov Chains and Mixing
  Times}, American Mathematical Society, second~ed., 2017.

\bibitem{Li2021JMP}
{\sc W.~C. Li}, {\em Hessian metric via transport information geometry}, J.
  Math. Phys, 62 (2021), p.~033301.

\bibitem{pmlr-v97-liu19i}
{\sc C.~Liu, J.~Zhuo, P.~Cheng, R.~Zhang, and J.~Zhu}, {\em Understanding and
  accelerating particle-based variational inference}, in ICML, vol.~97, 2019,
  pp.~4082--4092.

\bibitem{Liu2017NIPS}
{\sc Q.~Liu}, {\em Stein variational gradient descent as gradient flow}, in
  NeurIPS, vol.~30.

\bibitem{Liu2016NIPS}
{\sc Q.~Liu and D.~Wang}, {\em Stein variational gradient descent: {A} general
  purpose {B}ayesian inference algorithm}, in NeurIPS, vol.~29, 2016.

\bibitem{Logg2012Book}
{\sc A.~Logg, K.~A. Mardal, and G.~N. Wells}, {\em Automated Solution of
  Differential Equations by the Finite Element Method}, Springer, United
  Kingdom, 2012.

\bibitem{Reyes2015Book}
{\sc J.~C.~D. los Reyes}, {\em Numerical PDE-Constrained Optimization},
  Springer, New York, 2015.

\bibitem{Lu2019SIAM}
{\sc J.~Lu, Y.~Lu, and J.~Nolen}, {\em Scaling limit of the {S}tein variational
  gradient descent: the mean field regime}, SIAM J. Math. Anal., 5 (2019),
  pp.~648--671.

\bibitem{Matthews2016PHD}
{\sc A.~G. D.~G. Matthews}, {\em Scalable {G}aussian process inference using
  variational methods}, PhD thesis, University of Cambridge, 9 2016.

\bibitem{Nickl2020JEMS}
{\sc R.~Nickl}, {\em Betnstein-von {M}ises theorem for statistical inverse
  problems {I}: {S}chr\"{o}dinger equation}, J. Eur. Math. Soc., 22 (2020),
  pp.~2697--2750.

\bibitem{Pinski2015SIAMSC}
{\sc F.~J. Pinski, G.~Simpson, A.~M. Stuart, and H.~Weber}, {\em Algorithms for
  {K}ullback-{L}eibler approximation of probability measures in infinite
  dimensions}, SIAM J. Sci. Comput., 37 (2015), pp.~A2733--A2757.

\bibitem{Pinski2015SIAMMA}
{\sc F.~J. Pinski, G.~Simpson, A.~M. Stuart, and H.~Weber}, {\em
  Kullback-{L}eibler approximation for probability measures on infinite
  dimensional space}, SIAM J. Math. Anal., 47 (2015), pp.~4091--4122.

\bibitem{Prato2004Book}
{\sc G.~D. Prato}, {\em Kolmogorov Equations for Stochastic PDEs},
  Birkh\"{a}user Verlag, Basel, 2004.

\bibitem{Prato2006IDAnalysis}
{\sc G.~D. Prato}, {\em An Introduction to Infinite-Dimensional Analysis},
  Springer-Verlag, Berlin, 2006.

\bibitem{Ramsay2005Book}
{\sc J.~O. Ramsay and B.~W. Silverman}, {\em Functional Data Analysis},
  Springer, New York, second~ed., 2005.

\bibitem{Reed1980FunctionalAnalysis}
{\sc M.~Reed and B.~Simon}, {\em Functional Analysis {I}: Methods of Modern
  Mathematical Physics}, Elsevier (Singapore) Pte Ltd, revised and
  enlarged~ed., 2003.

\bibitem{Spantini2015SISC}
{\sc A.~Spantini, A.~Solonen, T.~Cui, J.~Martin, L.~Tenorio, and Y.~Marzouk},
  {\em Optimal low-rank approximations of {B}ayesian linear inverse problems},
  SIAM J. Sci. Comput., 37 (2015), pp.~A2451--A2487.

\bibitem{SteinwartBook}
{\sc I.~Steinwart and A.~Christmann}, {\em Support Vector Machines}, Springer,
  Germany, 2006.

\bibitem{Stuart2010AN}
{\sc A.~M. Stuart}, {\em Inverse problems: A {B}ayesian perspective}, Acta
  Numer., 19 (2010), pp.~451--559.

\bibitem{Sun2019ICLR}
{\sc S.~Sun, G.~Zhang, J.~Shi, and R.~Grosse}, {\em Functional variational
  {B}ayesian neural networks}, in ICLR, 2019.

\bibitem{Tarantola2005Book}
{\sc A.~Tarantola}, {\em Inverse Problem Theory and Methods for Model Parameter
  Estimation}, SIAM, United States, 2005.

\bibitem{Tarantola1982JG}
{\sc A.~Tarantola and B.~Valette}, {\em Inverse problems = quset for
  information}, J. Geophys., 50 (1982), pp.~159--170.

\bibitem{Trillos2015CJM}
{\sc N.~G. Trillos and D.~Slep\v~cev}, {\em On the rate of convergence of
  empirical measures in $\infty$-transportation distance}, Canad. J. Math., 67
  (2015), pp.~1358--1383.

\bibitem{Wang2019NIPS}
{\sc D.~Wang, Z.~Tang, C.~Bajaj, and Q.~Liu}, {\em Stein variational gradient
  descent with matrix-valued kernels}, in NeurIPS, vol.~33, 2019.

\bibitem{Wang2018SISC}
{\sc K.~Wang, T.~Bui-Thanh, and O.~Ghattas}, {\em A randomized maximum a
  posteriori method for posterior sampling of high dimensional nonlinear
  {B}ayesian inverse problems}, SIAM J. Sci. Comput., 40 (2018),
  pp.~A142--A171.

\bibitem{Wang2020arXiv}
{\sc Y.~Wang and W.~C. Li}, {\em Accelerated information gradient flows}.
\newblock arXiv:1909.02102, 2020.

\bibitem{Wang2019ICLR}
{\sc Z.~Wang, T.~Ren, J.~Zhu, and B.~Zhang}, {\em Function space particle
  optimization for {B}ayesian neural networks}, in ICLR, 2019.

\bibitem{Zhang2018IEEE}
{\sc C.~Zhang, J.~Butepage, H.~Kjellstrom, and S.~Mandt}, {\em Advances in
  variational inference}, IEEE T. Pattern Anal., 41 (2018), pp.~2008--2026.

\bibitem{Zhao2015IEEE}
{\sc Q.~Zhao, D.~Meng, Z.~Xu, W.~Zuo, and Y.~Yan}, {\em $l_{1}$-norm low-rank
  matrix factorization by variational {B}ayesian method}, IEEE T. Neur. Net.
  Lear., 26 (2015), pp.~825--839.

\bibitem{Zhou2020SIIMS}
{\sc Q.~Zhou, T.~Yu, X.~Zhang, and J.~Li}, {\em Bayesian inference and
  uncertainty quantification for medical image reconstruction with poisson
  data}, SIAM J. Imaging Sci., 13 (2020), pp.~29--52.

\end{thebibliography}


\begin{thebibliography}{10}

\bibitem{Agapiou2013SPA}
S.~Agapiou, S.~Larsson, and A.~M. Stuart.
\newblock Posterior contraction rates for the {B}ayesian approach to linear
  ill-posed inverse problems.
\newblock {\em Stoch. Proc. Appl.}, 123(10):3828--3860, 2013.

\bibitem{Ahmed1995SPTA}
N.~U. Ahmed and X.~Ding.
\newblock A semilinear {M}ckean-{V}lasov stochastic evolution equation in
  {H}ilbert space.
\newblock {\em Stoch. Proc. Appl.}, 60:65--85, 1995.

\bibitem{Bao2005IP}
G.~Bao and P.~Li.
\newblock Inverse medium scattering for the {Helmholtz} equation at fixed
  frequency.
\newblock {\em Inverse Probl.}, 21(5):1621--1641, 2005.

\bibitem{Beskos2015SC}
A.~Beskos, A.~Jasra, E.~A. Muzaffer, and A.~M. Stuart.
\newblock Sequential {M}onte {C}arlo methods for {B}ayesian elliptic inverse
  problems.
\newblock {\em Stat. Comput.}, 25:727–737, 2015.

\bibitem{Tan2013SISC}
T.~Bui-Thanh, O.~Ghattas, J.~Martin, and G.~Stadler.
\newblock A computational framework for infinite-dimensional {B}ayesian inverse
  problems part {I}: {T}he linearized case, with application to global seismic
  inversion.
\newblock {\em SIAM J. Sci. Comput.}, 35(6):A2494--A2523, 2013.

\bibitem{Thanh2016IPI}
T.~But-Thanh and Q.~P. Nguyen.
\newblock {FEM}-based discretization-invariant {MCMC} methods for
  {PDE}-constrained {B}ayesian inverse problems.
\newblock {\em Inverse Probl. Imag.}, 10(4):943--975, 2016.

\bibitem{Colton2019Book}
D.~Colton and R.~Kress.
\newblock {\em Inverse Acoustic and Electromagnetic Scattering Theory}.
\newblock Springer, Cham, fourth edition, 2019.

\bibitem{Cotter2013SS}
S.~L. Cotter, G.~O. Roberts, A.~M. Stuart, and D.~White.
\newblock {MCMC} methods for functions: modifying old algorithms to make them
  faster.
\newblock {\em Stat. Sci.}, 28(3):424--446, 2013.

\bibitem{Cotter2009IP}
Simon~L Cotter, Massoumeh Dashti, James~Cooper Robinson, and Andrew~M Stuart.
\newblock Bayesian inverse problems for functions and applications to fluid
  mechanics.
\newblock {\em Inverse Probl.}, 25(11):115008, 2009.

\bibitem{Dashti2017}
M.~Dashti and A.~M. Stuart.
\newblock The {Bayesian} approach to inverse problems.
\newblock {\em Handbook of Uncertainty Quantification}, pages 311--428, 2017.

\bibitem{Ghattas2021ActaNumerica}
O.~Ghattas and K.~Willcox.
\newblock Learning physics-based models from data: perspectives from inverse
  problems and model reduction.
\newblock {\em Acta Numer.}, 30:445--554, 2021.

\bibitem{Hinze2009Book}
M.~Hinze, R.~Pinnau, M.~Ulbrich, and S.~Ulbrich.
\newblock {\em Optimization with PDE Constraints}.
\newblock Springer Netherlands, 2009.

\bibitem{Jia2020IPI}
J.~Jia, J.~Peng, and J.~Gao.
\newblock Posterior contraction for empirical {B}ayesian approach to inverse
  problems under non-diagonal assumption.
\newblock {\em Inverse Probl. Imag.}, 15(2):201--228, 2020.

\bibitem{Kaipio2004Book}
J.~Kaipio and E.~Somersalo.
\newblock {\em Statistical and Computational Inverse Problems}.
\newblock Springer-Verlag, New York, 2005.

\bibitem{Reyes2015Book}
Juan Carlos~De los Reyes.
\newblock {\em Numerical PDE-Constrained Optimization}.
\newblock Springer, New York, 2015.

\bibitem{Lu2019SIAM}
J.~Lu, Y.~Lu, and J.~Nolen.
\newblock Scaling limit of the {S}tein variational gradient descent: the mean
  field regime.
\newblock {\em SIAM J. Math. Anal.}, 5(2):648--671, 2019.

\bibitem{Pillai2014SPDE}
N.~S. Pillai, A.~M. Stuart, and A.~H. Thiery.
\newblock Noisy gradient flow from a random walk in {H}ilbert space.
\newblock {\em Stoch. Partial. Differ.}, 2(2):196--232, 2014.

\bibitem{Prato2004Book}
G.~D. Prato.
\newblock {\em Kolmogorov Equations for Stochastic PDEs}.
\newblock Birkh\"{a}user Verlag, Basel, 2004.

\bibitem{Spantini2015SISC}
A.~Spantini, A.~Solonen, T.~Cui, J.~Martin, L.~Tenorio, and Y.~Marzouk.
\newblock Optimal low-rank approximations of {B}ayesian linear inverse
  problems.
\newblock {\em SIAM J. Sci. Comput.}, 37(6):A2451--A2487, 2015.

\bibitem{Stuart2010AN}
A.~M. Stuart.
\newblock Inverse problems: A {B}ayesian perspective.
\newblock {\em Acta Numer.}, 19:451--559, 2010.

\bibitem{Wang2018SISC}
K.~Wang, T.~Bui-Thanh, and O.~Ghattas.
\newblock A randomized maximum a posteriori method for posterior sampling of
  high dimensional nonlinear {B}ayesian inverse problems.
\newblock {\em SIAM J. Sci. Comput.}, 40(1):A142--A171, 2018.

\end{thebibliography}
